\documentclass[review]{elsarticle}
\usepackage{hyperref}
\usepackage{amsmath}
\usepackage{latexsym}

\journal{Journal of \LaTeX\ Templates}

\begin{document}

\begin{frontmatter}

\title{Reduction approach to the dynamics of interacting front solutions in a bistable reaction-diffusion system and its application to heterogeneous media}

\author{Kei Nishi}
\address{Kyoto Sangyo University, Kyoto 603-8555, Japan}
\fntext[myfootnote]{knishi@cc.kyoto-su.ac.jp}

\author{Yasumasa Nishiura}
\address{Tohoku University, Sendai 980-8577, Japan}
\fntext[myfootnote]{nishiura@tohoku.ac.jp}

\author{Takashi Teramoto}
\address{Asahikawa Medical University, Asahikawa 078-8510, Japan}
\fntext[myfootnote]{teramoto@asahikawa-med.ac.jp}




%
\begin{abstract}
  The dynamics of pulse solutions in a bistable reaction-diffusion system are studied analytically by reducing partial differential equations (PDEs) to finite-dimensional ordinary differential equations (ODEs). 
  For the reduction, we apply the multiple-scales method to the mixed ODE-PDE system obtained by taking a singular limit of the PDEs.
  The reduced equations describe the interface motion of a pulse solution formed by two interacting front solutions. This motion is in qualitatively good agreement with that observed for the original PDE system. Furthermore, it is found that the reduction not only facilitates the analytical study of the pulse solution, especially the specification of the onset of local bifurcations, but also allows us to elucidate the global bifurcation structure behind the pulse behavior. As an application, the pulse dynamics in a heterogeneous bump-type medium are explored numerically and analytically. The reduced ODEs clarify the transition mechanisms between four pulse behaviors that occur at different parameter values.
\end{abstract}


\end{frontmatter}


%
\section{Introduction}
%
%
%
\subsection{Historical background to the present study}
Reaction-diffusion systems have been widely utilized to model various spatiotemporal patterns observed in nature, such as chemical reactions, animal skin patterns, and morphogenesis. A class of spatially localized patterns is one of the most fundamental objects observed in many dissipative systems \cite{Vanag_Epstein_2007}.  
Nishiura and Mimura \cite{Nishiura_Mimura_1989} studied a two-component bistable system in one-dimensional space of the form
\begin{eqnarray}
  \hspace{1.8cm}\left \{
  \begin{array}{lll}
    \tau \hspace{0.5mm}\epsilon \hspace{0.5mm}u_t &=&  \epsilon^2 \hspace{0.5mm}u_{xx} + f(u,v)\hspace{0.5mm}, \\
    \hspace{0.2cm}v_t  &=&  D \hspace{0.5mm}v_{xx} + g(u,v)\hspace{0.5mm}, \hspace{0.8cm}(\hspace{0.5mm}t > 0, \hspace{1.0mm} x \in \mathbf{R}\hspace{0.5mm})\\
  \end{array}
  \right.
  \label{eq_tau_epsilon}
\end{eqnarray}
where $0 < \epsilon \ll 1$ and $\tau,  \hspace{0.5mm}D > 0$. This system, often referred to as a $\tau$\hspace{0.5mm}-\hspace{0.5mm}$\epsilon$ system, includes the two key parameters $\tau$ and $\epsilon$, which control the reaction rate and diffusion rate of the $u$ component, respectively.
The difference in the diffusion coefficients $\epsilon^2 \ll D$ between the $u$ and $v$ components gives rise to the localization of solutions as well as the sharp interface of width $\mathcal{O}(\epsilon)$ in the profile of the $u$ component. In contrast, the $v$ component behaves moderately . 
In the presence of odd symmetry in the nonlinearities $f$ and $g$, namely, $f(-u,-v)=-f(u,v)$ and $g(-u,-v)=-g(u,v)$, it has been shown that a stable stationary front solution exists for large values of $\tau$, and two stable traveling front solutions appear via a supercritical pitchfork bifurcation as $\tau$ decreases \cite{Nishiura_Mimura_Ikeda_Fujii_1990}.
These mathematical aspects of front solutions in bistable systems have been extensively studied \cite{Hagberg_Meron_1998}\cite{Hagberg_Meron_Rubistein_Zaltzman_1996}\cite{Hagberg_Meron_Rubistein_Zaltzman_1997}\cite{Ikeda_1998}\cite{Ikeda_Ikeda_2000}\cite{Ikeda_Ikeda_Mimura_2000}\cite{McKay_Kolokolnikov_Muir_2012}\cite{Mimura_Tabata_Hosono_1980}\cite{Sakamoto_1990}. Singularly perturbed three-component reaction-diffusion systems have also been investigated \cite{Chirilus_Doelman_Heijster_Rademacher_2015}\cite{Doelman_Heijster_Kaper_2009}\cite{Heijster_Chen_Nishiura_Teramoto_2016}\cite{Heijster_Doelman_Kaper_2008}\cite{Heijster_Doelman_Kaper_Promislow_2010}. 

In this paper, we numerically and analytically study the pulse dynamics arising in the $\tau$\hspace{0.5mm}-\hspace{0.5mm}$\epsilon$ system (\ref{eq_tau_epsilon}) with bistable-type nonlinearity. For the analysis, we do not deal with the original PDEs of (\ref{eq_tau_epsilon}), but with finite-dimensional ODEs that describe the pulse motion.
Our study is inspired in part by the work of Ei et al. \cite{Ei_Ikeda_Kawana_2008}\cite{Kusaka_2007}. In \cite{Ei_Ikeda_Kawana_2008}, they employed the following nonlinearities for $f$ and $g$ of (\ref{eq_tau_epsilon}):
\begin{eqnarray}
  \begin{array}{ll}
   f(u,v) = (\hspace{0.5mm}1/2 - u \hspace{0.5mm})(\hspace{0.5mm}u + v/2\hspace{0.5mm})(\hspace{0.5mm}u - v/2\hspace{0.5mm})\hspace{0.5mm}, \hspace{0.3cm} g(u,v) = u-v\hspace{0.5mm}.
  \end{array}  \label{eq_reaction_terms1}
\end{eqnarray}
System (\ref{eq_tau_epsilon}) with (\ref{eq_reaction_terms1}) is of bistable type and has the three spatially uniform solutions $P_{\pm} \equiv(\hspace{0.5mm}\pm 1/2, \hspace{0.5mm}\pm 1/2\hspace{0.5mm})$, and $P_0 \equiv (\hspace{0.5mm}0,\hspace{0.5mm}0\hspace{0.5mm})$, where $P_{\pm}$ are stable and $P_0$ is unstable.
By the center manifold reduction, Ei et al. showed the existence and stability of a front solution in (\ref{eq_tau_epsilon}) with (\ref{eq_reaction_terms1}), in which a stable stationary front solution becomes a pair of stable traveling front solutions via a pitchfork bifurcation. In particular, the dynamics near the bifurcation point are governed by the two-dimensional ODEs
\begin{eqnarray}
  \left \{
  \begin{array}{ll}
    \dot{l} = r\hspace{0.5mm}, \\
    \dot{r} = -M_1 \hspace{0.2mm}r^3 + M_2 \hspace{0.2mm}(\hspace{0.5mm}\tau_c - \tau \hspace{0.5mm}) \hspace{0.2mm}r\hspace{0.5mm}, \\
  \end{array}
  \right. \label{eq_Ei_ODE1}
\end{eqnarray}
where the variables $l$ and $r$ represent the location of the front interface and its velocity, respectively. The overdot $\hspace{0.5mm}\left< \hspace{1.0mm}\dot{} \hspace{1.0mm}\right>\hspace{0.5mm}$ denotes the time derivative $d/dt$, and $\tau_c$ is the pitchfork bifurcation point of the stationary front solution.
To analytically calculate the coefficients in Eq.~(\ref{eq_Ei_ODE1}), they employed matched asymptotic expansions to obtain the functions necessary for the reduction, such as the profiles of the stationary front solutions and the eigenfunctions. Thus, to the lowest order of $\epsilon$, they obtained the values $\tau_c = 1/(\hspace{0.2mm}4\sqrt{2\hspace{0.2mm}D}\hspace{0.2mm})$, $M_1 = 1/(\hspace{0.2mm}6 \hspace{0.2mm}D\hspace{0.2mm})$, $M_2 = 16 \sqrt{2\hspace{0.2mm}D}/3$.
Next, in \cite{Kusaka_2007}, Ei and Kusaka considered the problem of two interacting fronts in a bistable system in which the two front solutions are glued together to form a pulse shape.
Note, however, that system (\ref{eq_tau_epsilon}) with (\ref{eq_reaction_terms1}) does not have a pulse solution, because the interaction between the front solutions is repulsive for the nonlinearities in (\ref{eq_reaction_terms1}). Hence, the interfaces monotonically repel each other because of the odd symmetry in $f$ and $g$, and no pulse solution of finite width is formed. Such a repelling front motion has been extensively analyzed \cite{Car_Pego_1989}\cite{Ei_Ohta_1994}\cite{Fusco_Hale_1989}\cite{Kawasaki_Ohta_1982}.
Thus, to construct a pulse solution with a finite width, Ei and Kusaka broke the odd symmetry by adding perturbation terms to Eq.~(\ref{eq_reaction_terms1}):
\begin{eqnarray}
  \begin{array}{ll}
    f(u,v) = (\hspace{0.5mm}1/2 - u\hspace{0.5mm})(\hspace{0.5mm}u + v/2\hspace{0.5mm})(\hspace{0.5mm}u - v/2\hspace{0.5mm})+\delta_0 \hspace{0.5mm}h_1(u,v)\hspace{0.2mm}, \vspace{0.1cm}\\
    \hspace{0.0cm}g(u,v) = u-v+\delta_0 \hspace{0.2mm}h_2(u,v) \hspace{0.5mm},
  \end{array}  \label{eq_reaction_terms2}
\end{eqnarray}
where $0<\delta_0 \ll 1$. The perturbation terms play the role of gluing the two fronts together. In fact, for particular choices of $h_1$ and $h_2$, they numerically observed a stable pulse solution for system (\ref{eq_tau_epsilon}) with (\ref{eq_reaction_terms2}) that exhibited both traveling and oscillatory motion.
Using the center manifold reduction together with weak interaction theory \cite{Ei_2002}\cite{Ei_Mimura_Nagayama_2002}, Ei and coworkers obtained the equations of motion of the pulse behavior as
\begin{eqnarray} \left \{
	\begin{array}{rcl}
	\dot{l}_2 & = & r_2 + \tilde{M_0} \hspace{0.2mm}e^{-h/\sqrt{D}} + \delta_0 \hspace{0.5mm}\beta_1\hspace{0.5mm}, \\
	\dot{l}_1 & = & r_1 - \tilde{M_0} \hspace{0.2mm}e^{-h/\sqrt{D}} - \delta_0 \hspace{0.5mm}\beta_1\hspace{0.5mm}, \\
	\dot{r}_2 & = & -M_1 \hspace{0.2mm}r_2^3 + M_2 \hspace{0.2mm}(\hspace{0.5mm}\tau_c - \tau\hspace{0.5mm}) \hspace{0.5mm}r_2 + M_0 \hspace{0.2mm}e^{-h/\sqrt{D}} - \delta_0 \hspace{0.2mm}\beta_2\hspace{0.5mm}, \\
	\dot{r}_1 & = & -M_1 \hspace{0.2mm}r_1^3 + M_2 \hspace{0.2mm}(\hspace{0.5mm}\tau_c - \tau\hspace{0.5mm}) \hspace{0.5mm}r_1 - M_0 \hspace{0.2mm}e^{-h/\sqrt{D}} + \delta_0 \hspace{0.2mm}\beta_2\hspace{0.5mm}.
	\end{array}
	\right. \label{eq_Ei_ODE2}
\end{eqnarray}
The variables $l_2$ and $l_1$ represent the locations of the front interfaces, and $r_2$ and $r_1$ correspond to their velocities. Note that the third and fourth equations for $r_2$ and $r_1$ consist of three parts: the first two terms on the right-hand side describe the dynamics near the pitchfork bifurcation of each stationary front solution, which is exactly the same as in Eq.~(\ref{eq_Ei_ODE1}). The terms $\pm \hspace{0.2mm}\delta_0 \hspace{0.2mm}\beta_2$ come from adding the perturbation terms in Eq.~(\ref{eq_reaction_terms2}), and the terms $\pm \hspace{0.2mm}M_0 \hspace{0.2mm}e^{-h/\sqrt{D}}$ correspond to the exponentially weak interaction between the front solutions, with $h := l_2 - l_1$ being the width of the pulse solution. 
Each coefficient in Eq.~(\ref{eq_Ei_ODE2}) was obtained explicitly. In the particular case of $h_1(u,v)=0$ and $h_2(u,v)=u+v+1$, for instance, they were computed to the lowest order of $\epsilon$ as $M_0 = 8 \sqrt{D}/3$, $\tilde{M}_0 = 20 \sqrt{D}/9$, $\beta_1 = 8 \sqrt{D}/ 9$, and $\beta_2 = 16 \sqrt{D}/ 3$, where the values of $\tau_c$, $M_1$, and $M_2$ were given earlier.

By analyzing the resulting ODEs (\ref{eq_Ei_ODE2}), Ei et al. investigated the stability of a stationary pulse solution (i.e., standing pulse solution), and found that it underwent pitchfork and Hopf bifurcations as the parameter $\tau$ varied. 
As $\tau$ decreased, the standing pulse solution, which was stable for large values of $\tau$,
was destabilized by a pitchfork bifurcation at $\tau=\tau_{\rm p}$ and a stable traveling pulse solution appeared instead.
As $\tau$ decreased further, a Hopf bifurcation occurred as a secondary bifurcation at $\tau=\tau_{\rm \hspace{0.3mm}\scalebox{0.65}{H}}$. As a result, a standing breather solution, which exhibited in-phase oscillation of the two interfaces, remained unstable after appearing via the Hopf bifurcation point.
Regarding the order of the bifurcation points, the relation $\tau_{\rm \hspace{0.3mm}\scalebox{0.65}{H}} < \tau_{\rm p}$ was found analytically for ODE system (\ref{eq_Ei_ODE2}).
However, these analytic results disagreed with the pulse behavior observed by the direct numerical simulation of the original PDEs (\ref{eq_tau_epsilon}) with (\ref{eq_reaction_terms2}), which strongly suggested that the standing pulse solution underwent a Hopf bifurcation first and then a pitchfork bifurcation, implying that $\tau_{\rm p} < \tau_{\rm \hspace{0.3mm}\scalebox{0.65}{H}}$. This was supported by a more elaborate numerical approach that clarified the global bifurcation structure of a pulse solution for a bistable reaction-diffusion system similar to ours \cite{Nagayama_Ueda_Yadome_2010}.
Therefore, it seems most likely that the resulting ODEs (\ref{eq_Ei_ODE2}) fail to properly reproduce the pulse dynamics for the original PDE system, and hence the underlying bifurcation structure.

One of the main goals of this paper is to derive equations for the interacting fronts that properly describe the pulse motion, thus resolving the discrepancy between the pulse dynamics for the original PDE system and those for the reduced ODE system. For the reduction, we start with the limit system in (\ref{eq_hybrid}), which we introduced in our previous paper \cite{Nishi_Nishiura_Teramoto_2013}, and apply the multiple-scales method \cite{Hagberg_Meron_1998}\cite{Hagberg_Meron_Rubistein_Zaltzman_1996}\cite{Hagberg_Meron_Rubistein_Zaltzman_1997} rather than the center manifold reduction. This approach yields the four-dimensional ODEs
\begin{eqnarray}
\left \{
\begin{array}{l}
\displaystyle \; \; \; \; \dot{l}_2 = r_2, \quad \; \; \dot{l}_1 = r_1, \vspace{0.1cm} \\ 
~m_{0} \hspace{0.5mm}\dot{r}_2 = \sqrt{2} \hspace{0.5mm}(\tau_c-\tau) \hspace{0.5mm}r_2 - g_{3} \hspace{0.5mm}r^{3}_2
+(\hspace{0.5mm}G_0-G_1 \hspace{0.5mm}r_1\hspace{0.5mm}) \hspace{0.5mm}e^{-\frac{\;\;r_1+\phi(r_1)}{2 D}h}-\delta_0, \\
~m_{0} \hspace{0.5mm}\dot{r}_1 = \sqrt{2} \hspace{0.5mm}(\tau_c-\tau) \hspace{0.5mm}r_1 - g_{3} \hspace{0.5mm}r^{3}_1
-(\hspace{0.5mm}G_0+G_1 \hspace{0.5mm}r_2\hspace{0.5mm}) \hspace{0.5mm}e^{-\frac{-r_2+\phi(r_2)}{2 D}h}+\delta_0.
\end{array}
\right. \nonumber
\end{eqnarray}
which are the same as in (\ref{eq_ODE_renormalized_expanded2_0}) in Section 2.3.
The resulting equations are similar to (\ref{eq_Ei_ODE2}) in appearance, but they successfully reproduce the pulse dynamics observed for the original reaction-diffusion system, including the aforementioned order of the Hopf and pitchfork bifurcation points.
A key ingredient in remedying this discrepancy lies in the interaction terms $+(\hspace{0.5mm}G_0-G_1 \hspace{0.5mm}r_1\hspace{0.5mm}) \hspace{0.5mm}e^{-\frac{\;\;r_1+\phi(r_1)}{2 D}h}$ and $-(\hspace{0.5mm}G_0+G_1 \hspace{0.5mm}r_2\hspace{0.5mm}) \hspace{0.5mm}e^{-\frac{-r_2+\phi(r_2)}{2 D}h}$, which include higher-order terms that have largely been neglected \cite{Ei_2002}\cite{Ei_Mimura_Nagayama_2002}\cite{Ei_Ohta_1994}. 

\subsection{Pulse behavior in heterogeneous media}
A second issue covered by this paper is the pulse dynamics in heterogeneous media.
Recently, the effect of heterogeneity on pattern formation has attracted such attention that numerous experiments and numerical simulations have been performed by imposing various types of heterogeneity, such as parameter variations in space and time, and changes in system geometry \cite{Mikhailov_Showalter_2006}\cite{Vanag_Epstein_2008}. However, the mechanism for such heterogeneity-induced dynamics is still poorly understood considering the increasing number of experimental and numerical results \cite{Xin_2000}.
In this paper, we deal with a bump-type spatial heterogeneity, for which a system parameter changes from one value to another in a certain interval of $x$.
It is shown that our reduction method can be extended to the spatially heterogeneous case.
First, we numerically obtain a complete phase diagram of the pulse dynamics, both for the original PDE system and the reduced ODE system. This diagram suggests four different kinds of behavior as the height and width of the bump are varied as bifurcation parameters.
Based on the reduced ODEs, the transition mechanism for the pulse behavior is clarified from the perspective of a dynamical system. In particular, we characterize the difference between two kinds of behavior called pulse decomposition, whereby the pulse decomposes into two counter-propagating front solutions after encountering the heterogeneity.
We find that the transition is caused by a change in the behavior of either pulse interface, which is unique to a pulse consisting of two interacting front solutions.

The remainder of this paper is organized as follows.
In Section 2, we introduce the model system employed throughout the paper, namely, a system of two-component reaction-diffusion equations and its limit system of mixed ODE-PDE equations in the presence of a spatial heterogeneity.
Next, we give a brief description of the derivation of the reduced four-dimensional ODEs for the motion of two interfaces of the pulse solution. The details of the derivation are given in Appendices A and B. 
In Section 3, we investigate the reaction-diffusion system and the reduced ODE system for the pulse dynamics in the bump-type heterogeneous medium, and study the transition mechanism for the pulse behavior from the perspective of a dynamical system.
We conclude the paper with a summary and discussion in Section 4. 
%

%
\section{Equations of motion for two interfaces}
\subsection{Bistable reaction-diffusion system and its limit system}
\begin{figure}[h]
   \centering
 \includegraphics[width=12cm]{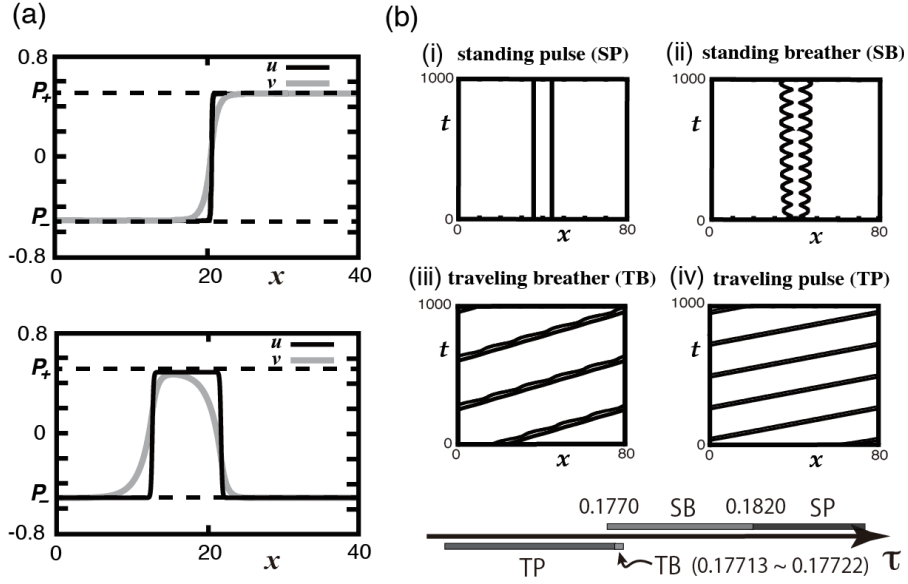}
  \caption{ 
(a) Profiles for front (top) and pulse (bottom) solutions observed numerically for Eq.~(\ref{eq_PDE}). The front solutions tend asymptotically to $P_{\pm}=(\hspace{0.5mm}\pm 1/2,\pm1/2+\delta_0\hspace{0.5mm})$ at $x \to \pm \infty$, whereas the pulse solutions tend to $P_{-}$.
(b) Four kinds of pulse behavior numerically observed for Eq.~(\ref{eq_PDE}) as the parameter $\tau$ varies: (i) standing pulse (\hspace{0.5mm}$\tau=0.1850$\hspace{0.5mm}), (ii) standing breather (\hspace{0.5mm}$\tau=0.1810$\hspace{0.5mm}), (iii)  traveling breather (\hspace{0.5mm}$\tau=0.17721$\hspace{0.5mm}), (iv) traveling pulse (\hspace{0.5mm}$\tau=0.1750$\hspace{0.5mm}). The horizontal and vertical axes denote $x$ and $t$, respectively, and the location of the interfaces for the $u$ component is shown.
The horizontal bar on the bottom schematically represents the parameter regime where the four kinds of pulse behavior were observed.
For the numerical simulation of the pulse solution, a finite difference method was employed (\hspace{0.5mm}$\Delta x=0.025$, $\Delta t=2.5 \times 10^{-3}$\hspace{0.5mm}) with a periodic boundary condition. The other parameters were fixed to $\epsilon=0.05$, $D=1.0$, and $\delta_0=0.001$.
}
 \label{pulse_dynamics_pde_homo}
\end{figure}

We consider the two-component reaction-diffusion system in one-dimensional space:
\begin{equation} 
\left\{
\begin{array}{l@{\,}ll} 
\displaystyle \tau \hspace{0.5mm} \epsilon \hspace{0.5mm} u_{t} \hspace{0.1cm} &=& \epsilon^{2} u_{xx} + \biggl(u+\frac{1}{2}\biggr)\biggl(\frac{1}{2}-u\biggr)\biggl(u-\frac{1}{2} v\biggr)\hspace{0.2mm}, \\
\displaystyle \hspace{0.2cm}v_{t} \hspace{0.1cm} &=& D v_{xx} + u - v +\delta(x)\hspace{0.2mm}, \\
\end{array}
\right. \label{eq_PDE}
\end{equation}

\noindent where $u(x,t)$ and $v(x,t)$ depend on space $x \in \mathbf{R}$ and $t>0$. The function $\delta(x)$ represents a given spatial heterogeneity. For now, we assume that the system is homogeneous with $\delta(x) \equiv \delta_0$, where $\delta_0$ is some positive constant. In this case, we find that the system is bistable, with the two stable spatially uniform solutions $(u_{\pm}(x,t),v_{\pm}(x,t)) \equiv (\pm 1/2, \pm 1/2+\delta_0)=: P_{\pm}$ and one unstable uniform solution $P_0$. 
As shown in Figure \ref{pulse_dynamics_pde_homo}(a), Eq.~(\ref{eq_PDE}) not only has a front solution connecting the two spatially uniform solutions $(u_{\pm},v_{\pm})$, but also a pulse solution that exhibits four kinds of behavior as the parameter $\tau$ varies: (i) standing pulse (SP), (ii) standing breather (SB), (iii) traveling breather (TB), (iv) traveling pulse (TP) (Fig.~\ref{pulse_dynamics_pde_homo}(b)).
%

In our previous paper \cite{Nishi_Nishiura_Teramoto_2013}, we took advantage of the fact that $\epsilon \ll D$ and considered the limit system obtained as $\epsilon \to 0$, for which the equation for $u$ is replaced by those for the two interfaces of the pulse solution located at $x=l_2, \hspace{0.5mm}l_1$:
\begin{equation}
\left\{
\begin{array}{l@{\,}ll} 
  \displaystyle{ \dot{l}_2 = - \frac{v(l_2,t)}{\sqrt{2} \tau},\hspace{0.5cm}\dot{l}_1 = \frac{v(l_1,t)}{\sqrt{2} \tau }, } \\
\displaystyle v_{\hspace{0.5mm}t}= D v_{\hspace{0.5mm}xx} + u(x\hspace{0.5mm};l_2,l_1) - v +\delta(x),
\end{array}
\right. \label{eq_hybrid}
\end{equation}

\noindent where $x=l_2(t),\hspace{0.5mm}l_1(t)$ are the locations of the interfaces ($\hspace{0.5mm}l_2>l_1$), the overdot denotes differentiation with respect to $t$, and the profile of the $u$ component $u(\hspace{0.2mm}x\hspace{0.2mm};l_2,\hspace{0.2mm}l_1\hspace{0.2mm})$ is given by 
\begin{equation}
u(x\hspace{0.5mm};l_2,l_1):= F(\hspace{0.2mm}x-l_1\hspace{0.2mm}) - F(\hspace{0.2mm}x - l_2\hspace{0.2mm}) - 1/\hspace{0.2mm}2, \label{eq_function_u}
\end{equation}
using a piecewise constant function
$F(x)=\hspace{0.2mm}1\hspace{0.2mm}/\hspace{0.2mm}2 \hspace{3.0mm} (\hspace{0.2mm}x \le 0\hspace{0.2mm}), \hspace{3.0mm} -1\hspace{0.2mm}/\hspace{0.2mm}2 \hspace{3.0mm}(\hspace{0.2mm}x > 0\hspace{0.2mm})$.
This type of mixed ODE-PDE system, which we call a hybrid system hereafter, typically arises in free boundary problems, such as the Stefan problem for phase transition \cite{Cannon_1984}, and has been used to study the dynamics of localized patterns in bistable reaction-diffusion systems \cite{Hilhorst_Nishiura_Mimura_1991}\cite{Ikeda_Ikeda_Mimura_2000}\cite{Ohta_Mimura_Kobayashi_1989}.
In our previous study \cite{Ei_Nishi_Nishiura_Teramoto_2015}\cite{Nishi_Nishiura_Teramoto_2013}, we examined the pulse dynamics for the hybrid system (\ref{eq_hybrid}) both numerically and analytically. In particular, we followed the method described in \cite{Ikeda_Ikeda_Mimura_2000} to investigate the eigenvalue problem for the stability of standing and traveling pulse solutions, confirming that the order of the Hopf and pitchfork bifurcations was consistent with that for the original PDEs, namely, $\tau_{\rm p} < \tau_{\rm \hspace{0.3mm}\scalebox{0.65}{H}}$.

Numerical results indicated that the hybrid system (\ref{eq_hybrid}) gave a good qualitative reproduction of the pulse dynamics observed in the original PDE system (\ref{eq_PDE}). Comparing systems (\ref{eq_PDE}) and (\ref{eq_hybrid}), we find that the equation for the $u$ component in Eq. (\ref{eq_PDE}) is replaced by two ODEs for the interfaces, whereas the $u$ term in the $v$ component becomes piecewise constant. This allowed us to analytically show the existence and stability of the time-independent pulse solutions for both the heterogeneous and homogeneous cases, and hence clarify several mechanisms for the pulse dynamics observed in a jump-type heterogeneous medium \cite{Nishi_Nishiura_Teramoto_2013}.
However, the reduced hybrid system (\ref{eq_hybrid}) still includes a PDE for the $v$ component. This has hindered the analytical study of more complicated behaviors, such as the standing breather and traveling breather shown in Figure~\ref{pulse_dynamics_pde_homo}(b--ii) and (b--iii). 

We take a further step to reduce the hybrid system to finite-dimensional ODEs by perturbatively solving the equation for the $v$ component.
To this end, we employ the multiple-scales method, which was applied to a similar hybrid system by Hagberg et al. to study the front motion in a bistable reaction-diffusion system in both one- and two-dimensional space \cite{Hagberg_Meron_1998}\cite{Hagberg_Meron_Rubistein_Zaltzman_1996}\cite{Hagberg_Meron_Rubistein_Zaltzman_1997}.
In the present paper, we extend their approach to the dynamics of two interacting front solutions, and derive a system of ODEs that describes the interface motion observed in the original PDE system (\ref{eq_PDE}).
%

%
\subsection{Reduction to ODE system}
In this section, a reduction from the hybrid system (\ref{eq_hybrid}) to ODEs is performed by perturbatively solving the PDE for the $v$ component in Eq.~(\ref{eq_hybrid}). For the reduction, a multiple-scales method is utilized, which proceeds as follows.
We introduce a slow time scale $T=\mu\hspace{0.5mm}t$, where $\mu$ is some infinitesimally small constant $0<\mu \ll 1$, and assume that the variables $l_2$, $l_1$, and $v$ in Eq.~(\ref{eq_hybrid}) have two timescales, $t$ and $T$. However, we also assume that, after a sufficiently long time, these variables are described by the slow timescale $T$ alone and become independent of the fast time scale $t$ as $t \to \infty$. Now, we expand the solution $v$ in power series of $\mu$ as
\begin{eqnarray}
v(x,T)=v_0(x,T)+\mu \hspace{0.5mm}v_1(x,T)+\mu^2 \hspace{0.5mm}v_2(x,T)+\mu^3 \hspace{0.5mm}v_3(x,T)+\cdots.
\end{eqnarray}
Substituting the above ansatz into the third equation in (\ref{eq_hybrid}) and collecting terms with equal powers of $\mu$, we can write equations for the functions $v_i(x,t)\hspace{0.2cm}(\hspace{0.5mm}i=0,1,2,\cdots\hspace{0.5mm})$ as
\begin{eqnarray}
\left \{
\begin{array}{l}
\displaystyle{ 
\mathcal{O}(\; \mu^0 \;): \quad \left( D \frac{\partial^2}{\partial x^2} -1 \right) v_0(x,T)+ u(\hspace{0.5mm}x;l_2,l_1\hspace{0.5mm})+\delta(x)=0 } \vspace{0.2cm} \\
\displaystyle{
\mathcal{O}(\; \mu^1 \;): \quad \left( D \frac{\partial^2}{\partial x^2} -1 \right) v_1(x,T)=\frac{\partial v_0}{\partial T} } \vspace{0.2cm} \\
\displaystyle{
\mathcal{O}(\; \mu^2 \;): \quad \left( D \frac{\partial^2}{\partial x^2} -1 \right) v_2(x,T)=\frac{\partial v_1}{\partial T} } \vspace{0.2cm} \\
\hspace{5.0cm} \vdots
\end{array}
\right. \label{eq_regular_perturbation}
\end{eqnarray}
where the order of the heterogeneity function $\delta(x)$ is assumed to be $\mathcal{O}(\mu^0)$, which does not affect the final result of the reduction. In general, the following relation holds for $n \ge 1$:
\begin{eqnarray}
\mathcal{O}(\mu^n): \quad \; \; \left( D \frac{\partial^2}{\partial x^2} -1 \right) v_n(x,T)=\frac{\partial v_{n-1}}{\partial T}
\end{eqnarray}
Each equation in (\ref{eq_regular_perturbation}) is a linear ODE with respect to $x$ and can be solved iteratively to give explicit solutions $v_n(x,T)$.
Substituting these $v(x,t)$ into the first two equations in (\ref{eq_hybrid}) yields closed-form ODEs for $l_1$ and $l_2$ (\hspace{0.5mm}see Appendix A for the details\hspace{0.5mm}).
Neglecting the higher-order terms, we obtain the following four-dimensional ODEs for $(\hspace{0.5mm}l_2(t), l_1(t), r_2(t), r_1(t)\hspace{0.5mm})$: 
\begin{eqnarray}
\left \{
\begin{array}{l}
\displaystyle \; \; \; \; \; \; \dot{l}_2 = r_2, \quad \; \; \dot{l}_1 = r_1, \vspace{0.1cm} \nonumber \\ 
\displaystyle m(r_2)\hspace{0.5mm}\dot{r}_2-M(r_1,h)\hspace{0.5mm}e^{-\frac{\hspace{1.5mm}r_1+\phi(r_1)}{2D} h} \hspace{1.0mm} \dot{r}_1 = g(r_2)+G(-r_1) \hspace{0.5mm}e^{-\frac{r_1+\phi(r_1)}{2D} h}-\Delta_0(l_2), \vspace{0.2cm} \nonumber\\
\displaystyle m(r_1)\hspace{0.5mm}\dot{r}_1 -M(r_2,h) \hspace{0.5mm}e^{-\frac{-r_2+\phi(r_2)}{2D} h} \hspace{1.0mm} \dot{r}_2= g(r_1)-G(r_2) \hspace{0.5mm}e^{-\frac{-r_2+\phi(r_2)}{2D} h}+\Delta_0(l_1),
\end{array} 
\right.
\end{eqnarray}
\begin{equation}
\hspace{5.0cm} \label{eq_ODE_renormalized_full}
\end{equation}
where $h(t):=l_2(t)-l_1(t)$ and $\phi(r):=\sqrt{r^2+4D}$. 
The function $\Delta_0(x)$ originates from the heterogeneity term $\delta(x)$ in Eq.~(\ref{eq_PDE}), which solves 
\begin{equation}
\left(\hspace{0.5mm} D \frac{d^2}{dx^2}-1 \hspace{0.5mm} \right) \Delta_0(x)+\delta(x)=0. \label{eq_delta}
\end{equation}
As $l_2$ and $l_1$ represent the locations of the two interfaces, the first two equations in (\ref{eq_ODE_renormalized_full}) indicate that the variables $r_2$ and $r_1$ correspond to the velocities of the interfaces. The prefactor functions are defined as
\begin{eqnarray}
\begin{array}{cc}
M(r,h) := \displaystyle \left( \hspace{0.5mm} \frac{6D^2}{\phi(r)^5} + \frac{3Dh}{\phi(r)^4} + \frac{h^2}{2\phi(r)^3} \hspace{0.5mm}\right), \vspace{0.2cm} \\
\displaystyle \hspace{1.0cm} m(r) := \frac{6 D^2}{\phi(r)^5}, \hspace{0.5cm} g(r) :=  -\sqrt{2}\hspace{0.2mm} \tau \hspace{0.2mm}r + \frac{r}{2 \phi(r)}, \hspace{0.5cm} G(r) := \frac{r + \phi(r)}{2 \phi(r)}.
\end{array} \label{eq_prefactor_function}
\end{eqnarray}
Note that Ohta et al. \cite{Ohta_Kiyose_Mimura_1997} and Mimura et al. \cite{Mimura_Nagayama_Ikeda_Ikeda_2000}\cite{Mimura_Nagayama_Ohta_2002} also studied a hybrid system similar to Eq.~(\ref{eq_hybrid}) with two ODEs for the interface motion and one PDE for the diffusion field, and reduced this to a four-dimensional system. The resulting ODEs were similar to Eq.~(\ref{eq_ODE_renormalized_full}) in appearance, but their reduction method was different from the aforementioned multiple-timescale technique.
They first applied a Fourier transform to the PDE and solved for $v(x,t)$ in an integral form. Next, to evaluate the integral, they expanded the interface locations $l_1(t)$ and $l_2(t)$ in series of $t$ up to the third order, thus obtaining closed-form second-order ODEs for $l_1(t)$ and $l_2(t)$.
They also utilized AUTO \cite{Doedel_2008} to study the bifurcation structure of the resulting second-order ODEs numerically, which will be mentioned in Section 2.3.

In this paper, we do not deal with the reduced ODE system (\ref{eq_ODE_renormalized_full}), but make further simplifications by assuming the interfaces move very slowly and are sufficiently far apart.\vspace{0.2cm}\\
\noindent
{\bf Proposition 1.}
{\it Assume that $|\hspace{0.3mm}r_2\hspace{0.3mm}|, \hspace{0.2mm}|\hspace{0.3mm}r_1\hspace{0.3mm}| \ll 1$, $|\hspace{0.3mm}\dot{r}_2\hspace{0.3mm}|, \hspace{0.2mm}|\hspace{0.3mm}\dot{r}_1\hspace{0.3mm}| \ll 1$, and $h \gg 1$. Then, the principal part of (\ref{eq_ODE_renormalized_full}) is given by
\begin{eqnarray}
\left \{
\begin{array}{l}
\; \; \; \; \; \; \dot{l}_2 = r_2, \quad \; \; \; \; \dot{l}_1 = r_1, \\ 
~m_{0} \hspace{0.5mm}\dot{r}_2 = \sqrt{2} \hspace{0.5mm}(\tau_c-\tau) \hspace{0.5mm}r_2 - g_{3} \hspace{0.5mm}r^{3}_2
+(\hspace{0.5mm}G_0-G_1 \hspace{0.5mm}r_1\hspace{0.5mm}) \hspace{0.5mm}e^{-\frac{\;\;r_1+\phi(r_1)}{2 D}h}-\Delta_0(l_2), \\
~m_{0} \hspace{0.5mm}\dot{r}_1 = \sqrt{2} \hspace{0.5mm}(\tau_c-\tau) \hspace{0.5mm}r_1 - g_{3} \hspace{0.5mm}r^{3}_1
-(\hspace{0.5mm}G_0+G_1 \hspace{0.5mm}r_2\hspace{0.5mm}) \hspace{0.5mm}e^{-\frac{-r_2+\phi(r_2)}{2 D}h}+\Delta_0(l_1),
\end{array}
\right. \label{eq_ODE_renormalized_expanded1}
\end{eqnarray}
where $m_0$, $\tau_c$, $g_3$, $G_0$, and $G_1$ are positive constants independent of $t$.} \vspace{0.2cm}\\
{\bf Proof.}\hspace{0.1cm}
{\it Expand the functions $m(r)$, $g(r)$, and $G(r)$ in (\ref{eq_prefactor_function}) as power series of $r$:
\begin{eqnarray}
\begin{array}{lll}
\displaystyle m(r) & = & m_0 + \mathcal{O}(r^2),  \vspace{0.2cm}\nonumber \\
\displaystyle g(r) & = &-\sqrt{2} \hspace{0.5mm}(\tau - \tau_c) \hspace{0.5mm}r - g_3\hspace{0.5mm}r^3 + \mathcal{O}(r^5), \vspace{0.2cm}\nonumber \\
\displaystyle G(r) & = & G_0 + G_1 \hspace{0.5mm}r + \mathcal{O}(r^3),
\end{array}
\end{eqnarray}
where the coefficients are given by $m_0=3/(\hspace{0.5mm}16\sqrt{D}\hspace{0.5mm})$, $\tau_c=1/(\hspace{0.5mm}4\sqrt{2D}\hspace{0.5mm})$, $g_3=1/(\hspace{0.5mm}32(\sqrt{D})^3\hspace{0.5mm})$, $G_0=1/2$, and $G_1=1/(\hspace{0.5mm}4\sqrt{D}\hspace{0.5mm})$.
Substituting these series into (\ref{eq_ODE_renormalized_full}) yields
\begin{eqnarray}
\left \{
\begin{array}{l}
\; \; \; \; \; \; \dot{l}_2 = r_2, \quad \; \; \; \; \dot{l}_1 = r_1, \vspace{0.2cm}\\ 
~(\hspace{0.5mm}m_{0}+\mathcal{O}(r_2^2)\hspace{0.5mm})\hspace{0.5mm}\dot{r}_2 -M(r_1,h)\hspace{0.5mm}e^{-\frac{\hspace{1.5mm}r_1+\phi(r_1)}{2D} h} \hspace{1.0mm} \dot{r}_1 \vspace{0.1cm}\\
\hspace{1.5cm}= \sqrt{2} \hspace{0.5mm}(\tau_c-\tau) \hspace{0.5mm}r_2 - g_{3} \hspace{0.5mm}r^{3}_2 + \mathcal{O}(r_2^5) \vspace{0.1cm}\\
\hspace{4.2cm}+(\hspace{0.5mm}G_0-G_1 \hspace{0.5mm}r_1+\mathcal{O}(r_1^3)\hspace{0.5mm}) \hspace{0.5mm}e^{-\frac{\;\;r_1+\phi(r_1)}{2 D}h}-\Delta_0(l_2), \vspace{0.2cm}\\
~(\hspace{0.5mm}m_{0}+\mathcal{O}(r_1^2)\hspace{0.5mm})\hspace{0.5mm}\dot{r}_1 -M(r_2,h) \hspace{0.5mm}e^{-\frac{-r_2+\phi(r_2)}{2D} h} \hspace{1.0mm} \dot{r}_2 \vspace{0.1cm}\\
\hspace{1.5cm}= \sqrt{2} \hspace{0.5mm}(\tau_c-\tau) \hspace{0.5mm}r_1 - g_{3} \hspace{0.5mm}r^{3}_1 + \mathcal{O}(r_1^5) \vspace{0.1cm}\\
\hspace{4.2cm}-(\hspace{0.5mm}G_0+G_1 \hspace{0.5mm}r_2+\mathcal{O}(r_2^3)\hspace{0.5mm}) \hspace{0.5mm}e^{-\frac{-r_2+\phi(r_2)}{2 D}h}+\Delta_0(l_1),
\end{array}
\right. \nonumber
\end{eqnarray}
Truncating the higher-order terms including $\mathcal{O}(r_i^3)\hspace{0.5mm}\times\hspace{0.5mm}e^{-\frac{\pm\hspace{0.5mm}r_i+\phi(r_i)}{2D} h}$ and $M(r_i,h)\hspace{0.5mm}\times\hspace{0.5mm}$ $\hspace{0.5mm}e^{-\frac{\pm\hspace{0.5mm}r_i+\phi(r_i)}{2D} h} \hspace{1.0mm} \dot{r}_i$ \hspace{0.2cm} $(\hspace{0.5mm}i=1,2\hspace{0.5mm})$, we have the equations in (\ref{eq_ODE_renormalized_expanded1}).} \hspace{2.5cm} \scalebox{1.0}{$\Box$} \vspace{0.2cm}

The four-dimensional ODEs (\ref{eq_ODE_renormalized_expanded1}) describe the motion of the two interfaces of the pulse solution observed in the PDEs (\ref{eq_PDE}), with the terms $+(\hspace{0.5mm}G_0-G_1 \hspace{0.5mm}r_1\hspace{0.5mm}) \hspace{0.5mm}e^{-\frac{\;\;r_1+\phi(r_1)}{2 D}h}$ and $-(\hspace{0.5mm}G_0+G_1 \hspace{0.5mm}r_2\hspace{0.5mm}) \hspace{0.5mm}e^{-\frac{-r_2+\phi(r_2)}{2 D}h}$ serving as exponentially weak interactions between the interfaces. 
In the limit $h \to \infty$, the interaction terms vanish and the system reduces to equations for the separate motion of two front solutions, for which, in the absence of the heterogeneity terms, $l=l_2=l_1$ and $r=r_2=r_1$ obey $\dot{l} = r\hspace{0.5mm}$, $m_{0} \hspace{0.5mm}\dot{r} = \sqrt{2} \hspace{0.5mm}(\tau_c-\tau) \hspace{0.5mm}r - g_{3} \hspace{0.5mm}r^{3}$ or, equivalently,
\begin{eqnarray}
  \left \{
  \begin{array}{ll}
    \dot{l} = r\hspace{0.5mm}, \\
    \dot{r} = -\tilde{M}_1 \hspace{0.2mm}r^3 + \tilde{M}_2 \hspace{0.2mm}(\hspace{0.5mm}\tau_c - \tau \hspace{0.5mm}) \hspace{0.2mm}r\hspace{0.5mm}, \\
  \end{array} \nonumber
  \right.
\end{eqnarray}
where $\tilde{M}_1 = g_3/\hspace{0.2mm}m_0= 1/(\hspace{0.2mm}6 \hspace{0.2mm}D\hspace{0.2mm})$ and $\tilde{M}_2 = \sqrt{2}/\hspace{0.2mm}m_0= 16 \sqrt{2\hspace{0.2mm}D}/3$. The equations are exactly the same as Eq.~(\ref{eq_Ei_ODE1}), which was derived by center manifold reduction around a drift bifurcation point of a front solution \cite{Ei_Ikeda_Kawana_2008}.

Furthermore, the interaction terms simplify to $\pm \hspace{0.5mm}G_0 \hspace{0.5mm}e^{-h/\sqrt{D}}$ when both $r_2$ and $r_1$ are set to zero. In weak interaction theory \cite{Ei_2002}, reduced ODEs usually have this simple form of interaction term.
In fact, the simple interaction term also appears in the ODEs of (\ref{eq_Ei_ODE2}), which are quite similar to those in (\ref{eq_ODE_renormalized_expanded1}) in appearance, except for the additional terms $\pm \hspace{0.5mm}\tilde{M_0} \hspace{0.2mm}e^{-h/\sqrt{D}}$ and $\pm \hspace{0.5mm}\delta_0 \hspace{0.5mm}\beta_1$, and the form of the interaction terms.
However, the form of the interaction term qualitatively changes the behavior of the solution to the reduced ODE system, and hence its bifurcation structure, as we shall see in the next section.

\subsection{Analysis of ODE system}
Section 2.1 presented numerical results for the homogeneous case in which the PDE system (\ref{eq_PDE}) has $\delta_0(x)\equiv\delta_0$. The function $\Delta_0(l)\equiv\delta_0$ in the reduced ODEs (\ref{eq_ODE_renormalized_expanded1}) can be found by solving Eq.~(\ref{eq_delta}), yielding $\Delta_0(l)\equiv\delta_0$.
In this section, we investigate the resulting ODEs (\ref{eq_ODE_renormalized_expanded1}) for the homogeneous case $\Delta_0(l)\equiv\delta_0$:
\begin{eqnarray}
\left \{
\begin{array}{l}
\displaystyle \; \; \; \; \dot{l}_2 = r_2, \quad \; \; \dot{l}_1 = r_1, \vspace{0.1cm} \\ 
~m_{0} \hspace{0.5mm}\dot{r}_2 = \sqrt{2} \hspace{0.5mm}(\tau_c-\tau) \hspace{0.5mm}r_2 - g_{3} \hspace{0.5mm}r^{3}_2
+(\hspace{0.5mm}G_0-G_1 \hspace{0.5mm}r_1\hspace{0.5mm}) \hspace{0.5mm}e^{-\frac{\;\;r_1+\phi(r_1)}{2 D}h}-\delta_0, \\
~m_{0} \hspace{0.5mm}\dot{r}_1 = \sqrt{2} \hspace{0.5mm}(\tau_c-\tau) \hspace{0.5mm}r_1 - g_{3} \hspace{0.5mm}r^{3}_1
-(\hspace{0.5mm}G_0+G_1 \hspace{0.5mm}r_2\hspace{0.5mm}) \hspace{0.5mm}e^{-\frac{-r_2+\phi(r_2)}{2 D}h}+\delta_0.
\end{array}
\right. \label{eq_ODE_renormalized_expanded2_0}
\end{eqnarray}

Using the relation $h=l_2-l_1$, this reduces to the three-dimensional system
\begin{eqnarray}
\left \{
\begin{array}{l}
\; \; \; \; \; \; \dot{h} = r_2-r_1, \\ 
~m_{0} \hspace{0.5mm}\dot{r}_2 = \sqrt{2} \hspace{0.5mm}(\tau_c-\tau) \hspace{0.5mm}r_2 - g_{3} \hspace{0.5mm}r^{3}_2
+(\hspace{0.5mm}G_0-G_1 \hspace{0.5mm}r_1\hspace{0.5mm}) \hspace{0.5mm}e^{-\frac{\;\;r_1+\phi(r_1)}{2 D}h}-\delta_0, \\
~m_{0} \hspace{0.5mm}\dot{r}_1 = \sqrt{2} \hspace{0.5mm}(\tau_c-\tau) \hspace{0.5mm}r_1 - g_{3} \hspace{0.5mm}r^{3}_1
-(\hspace{0.5mm}G_0+G_1 \hspace{0.5mm}r_2\hspace{0.5mm}) \hspace{0.5mm}e^{-\frac{-r_2+\phi(r_2)}{2 D}h}+\delta_0.
\end{array}
\right. \label{eq_ODE_renormalized_expanded2}
\end{eqnarray}
Figure \ref{pulse_dynamics_ode_homo}(a) shows a numerical simulation of the ODEs (\ref{eq_ODE_renormalized_expanded2}) for various values of the bifurcation parameter $\tau$. Compared with the numerical results in Figure \ref{pulse_dynamics_pde_homo}(b), the four kinds of pulse behavior (SP, SB, TB, and TP) are successfully reproduced, and are in qualitatively very good agreement with the pulse dynamics of the original PDE system.

Now that the system has been reduced to ODEs, the continuation and bifurcation software AUTO \cite{Doedel_2008} can be utilized to investigate the global bifurcation structure for the four kinds of solutions (Fig. \ref{pulse_dynamics_ode_homo}(b)).
The SP solution is stable for large values of $\tau$, but loses its stability
 via a Hopf bifurcation as $\tau$ decreases. From the Hopf bifurcation point,
 a stable SB solution branch appears. As $\tau$ decreases further, the SP solution undergoes a pitchfork bifurcation, at which point an unstable TP solution appears. The TP solution recovers its stability via a Hopf bifurcation, and a stable TB solution appears instead. The TB solution branch turns around at a fold bifurcation.
 Note that both the SB and TB solutions appear via Hopf bifurcations and their maximal $h$ values seem to diverge as the parameter $\tau$ approaches a particular value. However, the TB solution is only stable before the fold bifurcation, whereas the SB solution remains stable after the Hopf bifurcation, yielding a small coexistence regime of SB and TB solutions, as seen in Figure \ref{pulse_dynamics_ode_homo}(a).

As mentioned at the end of Section 2.2, Mimura et al. \cite{Mimura_Nagayama_Ikeda_Ikeda_2000} studied second-order ODEs similar to Eq.~(\ref{eq_ODE_renormalized_full}) using AUTO. They also obtained a qualitatively identical bifurcation structure for the four solution (SP, SB, TP, and TB) by varying a parameter that corresponds to $\tau$ in our ODEs. They observed that the pitchfork bifurcation on the SP solution branch became subcritical as another parameter varied, producing a saddle-node point on the TP solution branch. This was confirmed for Eq.~(\ref{eq_ODE_renormalized_expanded2}) with a larger value of $\delta_0$.
Nagayama et al. \cite{Nagayama_Ueda_Yadome_2010} numerically studied a two-component $\tau$\hspace{0.5mm}-\hspace{0.5mm}$\epsilon$ system with a bistable nonlinearity to clarify the bifurcation structure of the pulse solution. They obtained a bifurcation structure that is qualitatively the same as in Figure~\ref{pulse_dynamics_ode_homo}(b), except that the SB solution branch, unlike in the case of the ODEs (\ref{eq_ODE_renormalized_expanded2}), did not remain stable after the Hopf bifurcation point, but was destabilized via a saddle-node bifurcation. Considering that the PDE system investigated by Nagayama et al. was quite similar to Eq.~(\ref{eq_PDE}), it seems plausible to assume that the SB solution also undergoes a saddle-node bifurcation in our PDE system (\ref{eq_PDE}), implying that the reduced ODEs (\ref{eq_ODE_renormalized_expanded2}) do not completely reproduce the pulse dynamics of the original PDE system.
Remarkably, similar pulse dynamics have been found in an excitatory--inhibitory neural field model \cite{Folias_2017}. The author explored the behavior of traveling pulse solutions both numerically and analytically for a wide range of parameters, and revealed the bifurcation properties of stationary pulse solutions. 
The neural field model has a non-local interaction across a wide spatial region, in contrast to our system of diffusion-induced local interactions, and the similarities in the pulse behavior and associated bifurcation properties seem to indicate that the bistable reaction-diffusion system is a good approximation of the neural model system, at least for some parameter regimes.

The stationary solution $(\hspace{0.5mm}h,r_2,r_1\hspace{0.5mm})=(\hspace{0.5mm}h^*,r^*,r^*\hspace{0.5mm})$ to Eq.~(\ref{eq_ODE_renormalized_expanded2}) corresponds to either the TP solution ($r^* \neq 0$) or the SP solution ($r^* = 0$). The values of $\tau$ at which the SP solution undergoes the Hopf and pitchfork bifurcations can also be found analytically. \\

\noindent
{\bf Proposition 2.}
{\it For $0 < \delta_0 < 1/\hspace{0.5mm}2$,
\begin{equation}
(\hspace{0.5mm}h,r_2,r_1\hspace{0.5mm})=(\hspace{0.5mm}h^*,r^*,r^*\hspace{0.5mm})=(\hspace{0.5mm}-\sqrt{D}\hspace{0.5mm}\log(\hspace{0.5mm}2\hspace{0.5mm}\delta_0),\hspace{0.5mm}0,\hspace{0.5mm}0\hspace{0.5mm})\hspace{0.5mm},\label{eq_stationary_solution1}
\end{equation}
is a stationary solution to Eq.~(\ref{eq_ODE_renormalized_expanded2}). Furthermore, the stationary solution undergoes a pitchfork bifurcation at $\tau=\tau_d$ and a Hopf bifurcation at $\tau=\tau_H$ as $\tau$ varies. These points are explicitly given by
\begin{equation}
\displaystyle \tau_d=\tau_c-\sqrt{2}\hspace{0.5mm}(\hspace{0.5mm}G_1+ G_0\hspace{0.5mm}\phi_1\hspace{0.5mm}h^*\hspace{0.5mm})\hspace{0.5mm}\delta_0, \hspace{0.5cm}\tau_H=\tau_c+\sqrt{2}\hspace{0.5mm}(\hspace{0.5mm}G_1+ G_0\hspace{0.5mm}\phi_1\hspace{0.5mm}h^*\hspace{0.5mm})\hspace{0.5mm}\delta_0. \label{eq_bifurcation_points}
\end{equation}
where $\phi_1=1/(\hspace{0.5mm}2\hspace{0.5mm}D\hspace{0.5mm})$. \vspace{0.2cm}\\
{\bf Proof.}\hspace{0.1cm}
{\it The stationary solution $(\hspace{0.5mm}h,r_2,r_1\hspace{0.5mm})=(\hspace{0.5mm}h^*,r^*,r^*\hspace{0.5mm})$ to (\ref{eq_ODE_renormalized_expanded2}) satisfies
\begin{eqnarray}
\left \{
\begin{array}{l}
\displaystyle 0 = \sqrt{2} \hspace{0.5mm}(\tau_c-\tau) \hspace{0.5mm}r^* - g_{3} \hspace{0.5mm}r^{*3}
+(\hspace{0.5mm}G_0-G_1 r^*\hspace{0.5mm}) \hspace{0.5mm}e^{-\frac{\;\;r^*+\phi(r^*)}{2 D}h^*}-\delta_0, \\
\displaystyle 0 = \sqrt{2} \hspace{0.5mm}(\tau_c-\tau) \hspace{0.5mm}r^* - g_{3} \hspace{0.5mm}r^{*3}
-(\hspace{0.5mm}G_0+G_1 r^*\hspace{0.5mm}) \hspace{0.5mm}e^{-\frac{-r^*+\phi(r^*)}{2 D}h^*}+\delta_0.
\end{array}
\right. \label{eq_stationary_solution}
\end{eqnarray}
Substituting $r^* = 0$ into (\ref{eq_stationary_solution}), we find that $\displaystyle e^{-\frac{h^*}{\sqrt{D}}}=\hspace{0.5mm}2\hspace{0.5mm}\delta_0$ or, equivalently, $h^*=-\sqrt{D}\hspace{0.5mm}\log(\hspace{0.5mm}2\hspace{0.5mm}\delta_0\hspace{0.5mm})$. Therefore, we have that $(\hspace{0.5mm}h,r_2,r_1\hspace{0.5mm})=(\hspace{0.5mm}-\sqrt{D}\hspace{0.5mm}\log(\hspace{0.5mm}2\hspace{0.5mm}\delta_0),\hspace{0.5mm}0,\hspace{0.5mm}0\hspace{0.5mm})$ is a stationary solution to (\ref{eq_ODE_renormalized_expanded2}).
Linearizing Eq.~(\ref{eq_ODE_renormalized_expanded2}) around this stationary solution, we obtain the Jacobi matrix
\begin{eqnarray}
J_p=\left (
\begin{array}{ccc}
  0 & 1 & -1 \\
  -P & Q & -R \\
  P & -R & Q 
\end{array}
\right ) \label{eq_Jacobi_matrix}
\end{eqnarray}
where $P=2 \hspace{0.5mm}G_0\hspace{0.5mm}\phi_0\hspace{0.5mm}\delta_0/m_0$, $Q=\sqrt{2} \hspace{0.5mm}(\hspace{0.5mm}\tau_c-\tau\hspace{0.5mm})/m_0$, and $R=2\hspace{0.5mm}(\hspace{0.5mm}G_1+ G_0\hspace{0.5mm}\phi_1\hspace{0.5mm}h^*\hspace{0.5mm})/m_0$.
The constants $\phi_0$ and $\phi_1$ are the coefficients in the series expansion of
\begin{equation}
\displaystyle \Phi(r):=\hspace{0.5mm}\frac{r+\phi(r)}{2\hspace{0.5mm}D}\hspace{0.5mm}=\hspace{0.5mm}\phi_0 + \phi_1 \hspace{0.5mm}r + \mathcal{O}(r^2)\hspace{0.2mm}, \nonumber
\end{equation}
where $\phi(r)=\sqrt{r^2+4D}$, and are given by $\phi_0=1/\sqrt{D}$ and $\phi_1=1/(\hspace{0.5mm}2\hspace{0.5mm}D\hspace{0.5mm})$. The eigenvalue $\lambda$ of the matrix $J_p$ satisfies
\begin{equation}
  \displaystyle \lambda \hspace{0.5mm}(\hspace{0.5mm}\lambda -Q\hspace{0.5mm})^2+2\hspace{0.5mm}P\hspace{0.5mm}(\hspace{0.5mm}\lambda-Q\hspace{0.5mm})-R^2\hspace{0.5mm}\lambda+2\hspace{0.5mm}P\hspace{0.5mm}R=0 \label{eq_eigenvalue}
\end{equation}
The cubic equation for $\lambda$ in (\ref{eq_eigenvalue}) has one real root $\lambda_R$ and a pair of imaginary roots $\lambda_I$, $\bar{\lambda}_I$ for some appropriate values of $P$, $Q$, and $R$, which depend on the parameters in the ODEs (\ref{eq_ODE_renormalized_expanded2}). In particular, the eigenvalues $\lambda_R$ and $\lambda_I$ change depending on the parameter $\tau$ through $Q$.
Suppose that $\lambda_R=0$ for $\tau=\tau_d$, and $\lambda_I=\hspace{0.5mm}i\hspace{0.5mm}k\hspace{1.5mm}(\hspace{0.5mm}i=\sqrt{-1},\hspace{0.5mm}k \in \mathbf{R},\hspace{0.5mm}k>0\hspace{0.5mm})$ for $\tau=\tau_H$. Setting $\lambda=0$ in (\ref{eq_eigenvalue}) yields $Q=R$ or, equivalently, $\tau_d=\tau_c-\sqrt{2}\hspace{0.5mm}(\hspace{0.5mm}G_1+ G_0\hspace{0.5mm}\phi_1\hspace{0.5mm}h^*\hspace{0.5mm})\hspace{0.5mm}\delta_0$. Similarly, setting $\lambda=\hspace{0.5mm}i\hspace{0.5mm}k$ in (\ref{eq_eigenvalue}) yields $Q^3+(\hspace{0.5mm}P-R^2\hspace{0.5mm})\hspace{0.5mm}Q+PR=0$, which can be solved by $Q=-R$ and $k=\sqrt{2\hspace{0.5mm}P}$. Therefore, $\tau_H=\tau_c+\sqrt{2}\hspace{0.5mm}(\hspace{0.5mm}G_1+ G_0\hspace{0.5mm}\phi_1\hspace{0.5mm}h^*\hspace{0.5mm})\hspace{0.5mm}\delta_0$.
} \hspace{2.0cm} \scalebox{1.0}{$\Box$} \vspace{0.2cm}

\rm \noindent
The existence of the TP solution near the pitchfork bifurcation point $\tau_d$ can also be shown by the regular perturbation method.\\

\noindent
{\bf Proposition 3.}
{\it For $0 < \delta_0 < 1/\hspace{0.5mm}2$ and $0<\tau_d-\tau \ll 1$, the stationary solution $(\hspace{0.5mm}h,r_2,r_1\hspace{0.5mm})=(\hspace{0.5mm}h^*,r^*,r^*\hspace{0.5mm})$ to Eq.~(\ref{eq_ODE_renormalized_expanded2}) with $r^* \neq 0$ is approximately given by
\begin{eqnarray}
\begin{array}{l}
\displaystyle h^* = \hspace{0.5mm}h^{(0)} + h^{(2)} \hspace{0.5mm}(\hspace{0.5mm}\tau_d-\tau\hspace{0.5mm}) + \mathcal{O}(\hspace{0.5mm}\sqrt{\tau_d-\tau}^{\hspace{0.5mm}3}\hspace{0.5mm})\hspace{0.5mm},\vspace{0.1cm}\\
\displaystyle r^* = \hspace{0.5mm}r^{(1)} \hspace{0.5mm}\sqrt{\tau_d-\tau} + \mathcal{O}(\hspace{0.5mm}\sqrt{\tau_d-\tau} \hspace{0.8mm}{}^3\hspace{0.5mm})\hspace{0.5mm},\\
\end{array}
\label{eq_stationary_solution2_1}
\end{eqnarray}
where 
\begin{eqnarray}
\begin{array}{l}
  \displaystyle h^{(0)} = -\sqrt{D}\hspace{0.5mm}\log(\hspace{0.5mm}2\hspace{0.5mm}\delta_0)\hspace{0.5mm}, \hspace{0.5cm}r^{(1)} = \pm \hspace{0.5mm} \sqrt{ \frac{ \sqrt{2} }{ g_3-\tilde{g}_3 } }, \vspace{0.2cm}\\
  \displaystyle h^{(2)} =\hspace{1.0mm} 2\hspace{0.5mm}\frac{ -G_0 \phi_2 h^{(0)} +G_1 \phi_1 h^{(0)} + \frac{1}{2} \hspace{0.5mm}G_0 (\phi_1 h^{(0)})^2 }{ G_0 \hspace{0.5mm}\phi_2 \hspace{0.5mm}(\hspace{0.5mm}g_3-\tilde{g}_3\hspace{0.5mm})}. \vspace{0.1cm}\\
\end{array}
\label{eq_stationary_solution2_2}
\end{eqnarray}
The constants $\phi_0$, $\phi_1$, and $\phi_2$ are the coefficients in the series expansion of
\begin{equation}
\displaystyle \Phi(r):=\hspace{0.5mm}\frac{r+\phi(r)}{2\hspace{0.5mm}D}\hspace{0.5mm}=\hspace{0.5mm}\phi_0 + \phi_1 \hspace{0.5mm}r + \phi_2 \hspace{0.5mm}r^2 + \mathcal{O}(r^3)\hspace{0.2mm}, \nonumber
\end{equation}
with $\phi(r)=\sqrt{r^2+4D}$, which are given by $\phi_0=1/\sqrt{D}$, $\phi_1=1/(\hspace{0.5mm}2\hspace{0.5mm}D\hspace{0.5mm})$, and $\phi_2=1/(\hspace{0.5mm}8\hspace{0.5mm}D\hspace{0.5mm}\sqrt{D}\hspace{0.5mm})$. The constant $\tilde{g}_3$ is defined as \\
\begin{eqnarray}
  \begin{array}{l}
     \displaystyle \hspace{1.5cm}\tilde{g}_3  = \displaystyle \delta_0\hspace{0.5mm} \left \{ \left( 1+\frac{G^2_1}{G^2_0}- \frac{G_1}{G_0}\hspace{0.5mm}\frac{\phi_1}{\phi_2} \right) \phi_1\hspace{0.5mm}h^{(0)} \right.  \vspace{0.18cm}\\
     \hspace{4.0cm} \displaystyle \left. + \left( \frac{G_1}{G_0}-\frac{1}{2}\hspace{0.5mm}\frac{\phi_1}{\phi_2} \right) (\hspace{0.5mm}\phi_1\hspace{0.5mm}h^{(0)}\hspace{0.5mm})^2 + \frac{1}{3}\hspace{0.5mm}(\hspace{0.5mm}\phi_1\hspace{0.5mm}h^{(0)}\hspace{0.5mm})^3\hspace{0.5mm} \right \}\hspace{0.2mm}.
  \end{array} \nonumber
\end{eqnarray}
\vspace{0.2cm}\\
{\bf Proof.}\hspace{0.1cm}
{\it The stationary solution $(\hspace{0.5mm}h,r_2,r_1\hspace{0.5mm})=(\hspace{0.5mm}h^*,r^*,r^*\hspace{0.5mm})$ to (\ref{eq_ODE_renormalized_expanded2}) satisfies
\begin{eqnarray}
\left \{
\begin{array}{l}
\displaystyle 0 = \sqrt{2} \hspace{0.5mm}(\tau_c-\tau) \hspace{0.5mm}r^* - g_{3} \hspace{0.5mm}r^{*3}
+(\hspace{0.5mm}G_0-G_1 r^*\hspace{0.5mm}) \hspace{0.5mm}e^{-\frac{\;\;r^*+\phi(r^*)}{2 D}h^*}-\delta_0, \\
\displaystyle 0 = \sqrt{2} \hspace{0.5mm}(\tau_c-\tau) \hspace{0.5mm}r^* - g_{3} \hspace{0.5mm}r^{*3}
-(\hspace{0.5mm}G_0+G_1 r^*\hspace{0.5mm}) \hspace{0.5mm}e^{-\frac{-r^*+\phi(r^*)}{2 D}h^*}+\delta_0.
\end{array}
\right. \label{eq_stationary_solution2}
\end{eqnarray}
We now introduce a small parameter $0 < \eta \ll 1$, and assume that the solution has the form
\begin{eqnarray}
 \begin{array}{l}
   \displaystyle \hspace{2.0cm}\tau_c - \tau = \hspace{0.5mm}\tau_0 + \eta^2 \hspace{0.5mm},\vspace{0.1cm}\\
  \displaystyle h^* = \hspace{0.5mm}h^{(0)} + h^{(1)} \hspace{0.5mm} \eta + h^{(2)} \hspace{0.5mm} \eta^2 + h^{(3)} \hspace{0.5mm} \eta^3 + \mathcal{O}(\hspace{0.5mm}\eta^4\hspace{0.5mm})\hspace{0.5mm},\vspace{0.1cm}\\
\displaystyle r^* = \hspace{0.5mm}r^{(1)} \hspace{0.5mm}\eta + r^{(2)} \hspace{0.5mm}\eta^2 + r^{(3)} \hspace{0.5mm}\eta^3 + \mathcal{O}(\hspace{0.5mm}\eta^4\hspace{0.5mm})\hspace{0.5mm},\\
\end{array}
\label{eq_stationary_solution_ansatz}
\end{eqnarray}
where \hspace{0.5mm}$\tau_0$, $h^{(i)}$, and \hspace{0.5mm}$r^{(i)}$ \hspace{1.0mm} {\rm (}\hspace{0.5mm}$i=0,1,2,3$\hspace{0.5mm} {\rm )} are unknown constants to be determined. Note that rescaling $\eta$ allows us to set the coefficient of $\eta^2$ in the first equation to $1$. Substituting the ansatz (\ref{eq_stationary_solution_ansatz}) into the first equation in (\ref{eq_stationary_solution2}) and collecting powers of $\eta$, we obtain 
\begin{eqnarray}
  \left \{
  \begin{array}{ll}
    \mathcal{O}(\eta^{0}):
    \displaystyle \hspace{0.2cm}G_0 Q_0 - \delta_0 = 0 \hspace{0.5mm}, \vspace{0.1cm}\\
    \mathcal{O}(\eta^1): 
    \displaystyle \hspace{0.2cm}G_0 Q_1 - G_1 r^{(1)} Q_0 + \sqrt{2}  \hspace{0.5mm}\tau_0 r^{(1)} = 0 \hspace{0.5mm}, \vspace{0.1cm}\\
    \mathcal{O}(\eta^2):
    \displaystyle \hspace{0.2cm}G_0 Q_2 - G_1 r^{(1)} Q_1 - G_1 r^{(2)} Q_0 + \sqrt{2}  \hspace{0.5mm}\tau_0 r^{(2)} = 0 \hspace{0.5mm}, \vspace{0.1cm}\\
    \mathcal{O}(\eta^3):
    \displaystyle \hspace{0.2cm}\sqrt{2} \hspace{0.5mm}\tau_0 r^{(3)} + \sqrt{2} \hspace{0.5mm}\tau_1 r^{(1)} - g_3 (r^{(1)})^3 \\
    \hspace{2.0cm} + (G_0 Q_3 - G_1 r^{(1)} Q_2 - G_1 r^{(2)} Q_1 - G_1 r^{(3)} Q_0)= 0 \hspace{0.5mm}, \vspace{0.1cm}\\
  \end{array}
\right . \label{eq_stationary_solution_orders}
\end{eqnarray}
where the constants $Q_i$ \hspace{0.1cm}$(\hspace{0.5mm}i=0,\cdots,3\hspace{0.5mm})$ are defined as
\begin{eqnarray}
  \left \{
  \begin{array}{ll}
    \displaystyle  \hspace{0.1cm}Q_0 = e^{- \phi_0 h^{(0)}} \hspace{0.5mm}, \vspace{0.1cm}\\
    \displaystyle  \hspace{0.1cm}Q_1 = - e^{- \phi_0 h^{(0)}} \hspace{0.5mm}p_1\hspace{0.5mm}, \vspace{0.1cm}\\
    \displaystyle  \hspace{0.1cm}Q_2 = e^{- \phi_0 h^{(0)}} \hspace{0.5mm}\left (-p_2 + \frac{1}{2}\hspace{0.5mm}p^2_1\right)\hspace{0.5mm}, \vspace{0.1cm}\\
    \displaystyle  \hspace{0.1cm}Q_3 = e^{- \phi_0 h^{(0)}} \hspace{0.5mm}\left (-p_3 + p_1\hspace{0.5mm}p_2 - \frac{1}{6}\hspace{0.5mm}p^3_1\right)\hspace{0.5mm}, \vspace{0.1cm}\\
  \end{array}
 \right. \label{eq_stationary_solution_Q}
\end{eqnarray}
with
\begin{eqnarray}
  \left \{
  \begin{array}{ll}
    \displaystyle \hspace{0.1cm}p_1 = \phi_0 h^{(1)} + \phi_1 h^{(0)} r^{(1)} \hspace{0.5mm}, \vspace{0.1cm}\\
    \displaystyle \hspace{0.1cm}p_2 = \phi_2 h^{(2)} + \phi_1 h^{(1)} r^{(1)} + h^{(0)} \hspace{0.5mm} \left \{ \phi_1 r^{(2)} + \phi_2 (r^{(1)})^2 \right \} \hspace{0.5mm}, \vspace{0.1cm}\\
    \displaystyle \hspace{0.1cm}p_3 = \phi_0 h^{(3)} + \phi_1 h^{(2)} r^{(1)} + h^{(1)} \hspace{0.5mm} \left \{ \phi_1 r^{(2)} + \phi_2 (r^{(1)})^2 \right \} \vspace{0.1cm}\\
    \hspace{4.3cm} + h^{(0)} \hspace{0.5mm} \left \{ \phi_1 r^{(3)} + 2 \phi_2 r^{(1)} r^{(2)} \right \} \hspace{0.5mm}, \vspace{0.1cm}\\    
  \end{array}
 \right. \label{eq_stationary_solution_q}
\end{eqnarray}
The first equation in (\ref{eq_stationary_solution_orders}) yields $e^{- \phi_0 h^{(0)}}=\delta_0 / G_0$ or, equivalently,
\begin{eqnarray}
h^{(0)}=-\sqrt{D}\hspace{0.5mm}\log(\hspace{0.5mm}2\hspace{0.5mm}\delta_0\hspace{0.5mm}). \label{eq_stationary_solution_order0}
\end{eqnarray}
Substituting (\ref{eq_stationary_solution_Q}) and (\ref{eq_stationary_solution_q}) into the second equation for $\mathcal{O}(\eta^1)$ in (\ref{eq_stationary_solution_orders}) yields
\begin{eqnarray}
  \begin{array}{ll}
    \displaystyle G_0 \hspace{0.5mm}\phi_0 \hspace{0.5mm}h^{(1)} + (G_0 \hspace{0.5mm}\phi_1 \hspace{0.5mm}h^{(0)}+G_1-\hspace{0.5mm}\sqrt{2}\tau_0\hspace{0.5mm}e^{\phi_0 h^{(0)}})\hspace{0.5mm}r^{(1)}=0\hspace{0.2mm}.
  \end{array} \label{eq_stationary_solution_order1-1}
\end{eqnarray}
By substituting the ansatz (\ref{eq_stationary_solution_ansatz}) into the second equation in (\ref{eq_stationary_solution2}), we have a similar equation:
\begin{eqnarray}
  \begin{array}{ll}
    \displaystyle - G_0 \hspace{0.5mm}\phi_0 \hspace{0.5mm}h^{(1)} + (G_0 \hspace{0.5mm}\phi_1 \hspace{0.5mm}h^{(0)}+G_1-\hspace{0.5mm}\sqrt{2}\tau_0\hspace{0.5mm}e^{\phi_0 h^{(0)}})\hspace{0.5mm}r^{(1)}=0\hspace{0.2mm},
  \end{array} \label{eq_stationary_solution_order1-2}
\end{eqnarray}
which is also obtained by the change in sign
\begin{eqnarray}
  \begin{array}{ll}
    \displaystyle G_0 \to - G_0 \hspace{0.5mm}, \hspace{0.3cm}\phi_1 \to - \phi_1 \hspace{0.5mm}, \hspace{0.3cm}\delta_0 \to - \delta_0 \hspace{0.5mm}.
  \end{array} \label{eq_stationary_solution_sign_change}
\end{eqnarray}
Equations (\ref{eq_stationary_solution_order1-1}) and (\ref{eq_stationary_solution_order1-2}) for $h^{(1)}$ and $r^{(1)}$ lead to
\begin{eqnarray}
  \begin{array}{ll}
    \displaystyle \hspace{0.1cm}h^{(1)}=0, \hspace{0.4cm}\tau_0 = (\hspace{0.5mm}G_1+ G_0\hspace{0.5mm}\phi_1\hspace{0.5mm}h^*\hspace{0.5mm})\hspace{0.5mm}\delta_0/( \sqrt{2}\hspace{0.5mm}G_0).
  \end{array} \label{eq_stationary_solution_order1-3}
\end{eqnarray}
Similarly, from the equation for $\mathcal{O}(\eta^2)$ in (\ref{eq_stationary_solution_orders}) together with the sign inversion (\ref{eq_stationary_solution_sign_change}), we obtain
\begin{eqnarray}
  \displaystyle h^{(2)} = \frac{ -G_0 \hspace{0.5mm}\phi_2 \hspace{0.5mm}h^{(0)}  +G_1 \hspace{0.5mm}\phi_1 \hspace{0.5mm}h^{(0)} + \frac{1}{2} \hspace{0.5mm}G_0 \hspace{0.5mm}(\phi_1 h^{(0)})^2  }{ G_0 \hspace{0.5mm}\phi_2 } \hspace{0.5mm}(r^{(1)})^2\hspace{0.2mm}. \label{eq_stationary_solution_order2-1}
\end{eqnarray}
Finally, from the equation for $\mathcal{O}(\eta^3)$ in (\ref{eq_stationary_solution_orders}), we obtain
\begin{eqnarray}
  \left \{
  \begin{array}{ll}
     \displaystyle \sqrt{2} \hspace{0.5mm}r^{(1)} - g_3 (r^{(1)})^3 + \tilde{g}_3 (r^{(1)})^3 = 0\hspace{0.2mm},\vspace{0.1cm}\\
  \displaystyle h^{(3)} = \frac{ -2 \hspace{0.5mm}G_0 \hspace{0.5mm}\phi_2 \hspace{0.5mm}h^{(0)} + 2 \hspace{0.5mm}G_1 \hspace{0.5mm}\phi_1 \hspace{0.5mm}h^{(0)} +G_0 \hspace{0.5mm}(\phi_1 \hspace{0.5mm}h^{(0)})^2 }{ G_0 \hspace{0.5mm}\phi_0 } \hspace{0.5mm}r^{(1)}\hspace{0.5mm}r^{(2)}\hspace{0.2mm},
  \end{array} \label{eq_stationary_solution_order3-1}
 \right. 
\end{eqnarray}
where $\tilde{g}_3$ is a constant defined by
\begin{eqnarray}
   \begin{array}{ll}
     \displaystyle \tilde{g}_3  = \displaystyle \delta_0\hspace{0.5mm} \left \{ \left( 1+\frac{G^2_1}{G^2_0}- \frac{G_1}{G_0}\hspace{0.5mm}\frac{\phi_1}{\phi_2} \right) \phi_1\hspace{0.5mm}h^{(0)} + \left( \frac{G_1}{G_0}-\frac{1}{2}\hspace{0.5mm}\frac{\phi_1}{\phi_2} \right) (\hspace{0.5mm}\phi_1\hspace{0.5mm}h^{(0)}\hspace{0.5mm})^2 \right.  \vspace{0.15cm}\\
\hspace{8.5cm} \displaystyle \left.+ \frac{1}{3}\hspace{0.5mm}(\hspace{0.5mm}\phi_1\hspace{0.5mm}h^{(0)}\hspace{0.5mm})^3\hspace{0.5mm} \right \}\hspace{0.2mm}.
  \end{array} \label{eq_stationary_solution_order3-2}
\end{eqnarray}
From the first equation in (\ref{eq_stationary_solution_order3-1}), the value of $r^{(1)} \neq 0$ is determined as
\begin{eqnarray}
  \begin{array}{ll}
    \displaystyle r^{(1)} =\hspace{0.2cm}\pm \hspace{0.5mm}\sqrt{\frac{\sqrt{2}}{g_3-\tilde{g}_3}} \hspace{0.4mm}.
  \end{array} \label{eq_stationary_solution_order3-4}
\end{eqnarray}
Now that we have obtained $r^{(1)}$ and $h^{(i)}$ \hspace{0.1cm}$(\hspace{0.5mm}i=0,1,2\hspace{0.5mm})$, we construct the approximate solution in (\ref{eq_stationary_solution_ansatz}). The first equation in (\ref{eq_stationary_solution_ansatz}) leads to $\eta=\sqrt{(\hspace{0.2mm}\tau_c-\tau_0\hspace{0.2mm})-\tau}= \sqrt{\tau_d-\tau}$, where $\tau_d:=\tau_c-\tau_0$; this coincides with $\tau_d$ in Eq.~(\ref{eq_bifurcation_points}) of Proposition 2. The second and third equations in (\ref{eq_stationary_solution_ansatz}) yield
\begin{eqnarray}
 \begin{array}{l}
   \displaystyle h^* = \hspace{0.5mm}h^{(0)} + h^{(2)} \hspace{0.5mm}(\hspace{0.5mm}\tau_d-\tau\hspace{0.5mm}) + \mathcal{O}(\hspace{0.5mm}\sqrt{\tau_d-\tau}^{\hspace{0.5mm}3}\hspace{0.5mm})\hspace{0.5mm},\vspace{0.1cm}\\
   \displaystyle r^* = \hspace{0.5mm}r^{(1)} \hspace{0.5mm}\sqrt{\tau_d-\tau} + \mathcal{O}(\hspace{0.5mm}\sqrt{\tau_d-\tau}^{\hspace{0.5mm}2}\hspace{0.5mm})\hspace{0.5mm},\\
\end{array}
\label{eq_stationary_solution_3-1}
\end{eqnarray}
where $h^{(0)}$ and $h^{(2)}$ are defined in Eqs.~(\ref{eq_stationary_solution_order0}) and (\ref{eq_stationary_solution_order2-1}), respectively.
} \hspace{0.9cm} \scalebox{1.0}{$\Box$} \vspace{0.4cm}\\
\rm
For later use, let us define the SP and TP solutions more rigorously based on the above two propositions.\\

\noindent
{\bf Definition 1.} 
{\it The stationary solution $(\hspace{0.5mm}h,r_2,r_1\hspace{0.5mm})=(\hspace{0.5mm}h^*,r^*,r^*\hspace{0.5mm})$ to Eq.~(\ref{eq_ODE_renormalized_expanded2}) with $r^* = 0$ is defined as ${\rm SP}$, which is explicitly given by Eq.~(\ref{eq_stationary_solution1}). The stationary solution $(\hspace{0.5mm}h,r_2,r_1\hspace{0.5mm})=(\hspace{0.5mm}h^*,r^*,r^*\hspace{0.5mm})$ to Eq.~(\ref{eq_ODE_renormalized_expanded2}) with either $r^* > 0$ or $r^* < 0$ is defined as \hspace{0.8mm}${\rm TP}^{\hspace{0.5mm}+}$\hspace{0.8mm}or \hspace{0.8mm}${\rm TP}^{\hspace{0.5mm}-}$, respectively, which is approximately given by Eq.~(\ref{eq_stationary_solution2_1}) near $\tau = \tau_d$.
}\\

\rm
Proposition 3 analytically guarantees the existence of the left- and right-moving traveling pulse solutions near the pitchfork bifurcation point of the SP solution.
Numerical simulations indicate that these solutions are initially unstable, but recover their stability via a Hopf bifurcation (Fig.~\ref{pulse_dynamics_ode_homo}(b)).
As for the SP stability, Proposition 2 shows that the SP is destabilized by a Hopf bifurcation. In particular, as the constants on the right-hand side of (\ref{eq_bifurcation_points}) are all positive, we find that $\tau_H>\tau_d$, which agrees with the numerical results shown in Figure~\ref{pulse_dynamics_ode_homo}(b).
As mentioned at the end of the previous section, an interaction term of the form $G_0 \hspace{0.5mm}e^{-\frac{h}{\sqrt{D}}}$ usually appears in the weak interaction theory. In fact, the interaction terms $+(\hspace{0.5mm}G_0-G_1 \hspace{0.5mm}r_1\hspace{0.5mm}) \hspace{0.5mm}e^{-\frac{\;\;r_1+\phi(r_1)}{2 D}h}$ and $-(\hspace{0.5mm}G_0+G_1 \hspace{0.5mm}r_2\hspace{0.5mm}) \hspace{0.5mm}e^{-\frac{-r_2+\phi(r_2)}{2 D}h}$ in Eq.~(\ref{eq_ODE_renormalized_expanded2}) also reduce to $\pm \hspace{0.5mm}G_0 \hspace{0.5mm}e^{-\frac{h}{\sqrt{D}}}$ when $r_1=r_2=0$. In this case, both $G_1$ and $\phi_1$ vanish in (\ref{eq_bifurcation_points}), and we are led to the following corollary. \vspace{0.2cm}

\noindent
{\bf Corollary 1.}
{\it
  If $G_1=0$ and $\phi_1=0$ in (\ref{eq_bifurcation_points}), then $\tau_d=\tau_H=\tau_c$.
}\vspace{0.2cm}

Therefore, the two bifurcation points coincide if the interaction terms are of the form $\pm \hspace{0.5mm}G_0 \hspace{0.5mm}e^{-\frac{h}{\sqrt{D}}}$, resulting in a discrepancy in the pulse dynamics of the PDE system and the reduced ODE system.
The interaction terms in our reduced ODEs (\ref{eq_ODE_renormalized_expanded1}) include, in a sense, higher-order interaction terms, with $G_0 \hspace{0.5mm}e^{-\frac{h}{\sqrt{D}}}$ being the lowest-order term.
In fact, dividing both sides of the second and third equations in (\ref{eq_ODE_renormalized_expanded2}) by $m_0$, we find that $G_0/\hspace{0.5mm}m_0= 8 \sqrt{D}/3$, which coincides with the coefficient $M_0$ of the interaction terms in (\ref{eq_Ei_ODE2}).

Thus, the above analysis suggests that higher-order interaction terms must be taken into account for our particular case of the front--back pulse dynamics observed in PDE system (\ref{eq_PDE}). 
\begin{figure}[ht!]
   \centering
 \includegraphics[width=10cm]{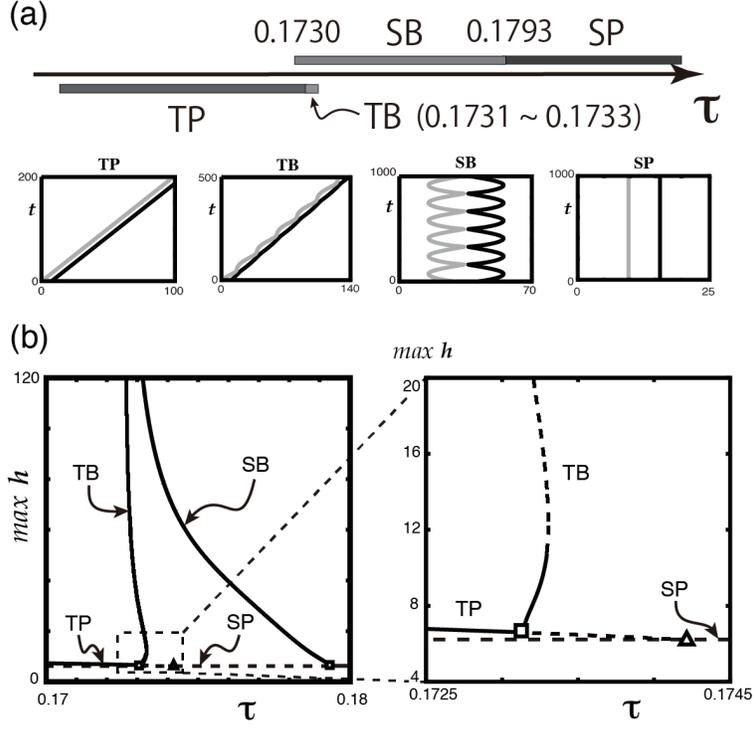}
 \caption{(a) Schematic phase diagram obtained by numerical simulation of the ODE system (\ref{eq_ODE_renormalized_expanded1}) and typical spatiotemporal plots of $l_2(t)$ (black line) and $l_1(t)$ (gray line) for the four kinds of behavior. For the spatiotemporal plots, values of $\tau=0.1820$ (SP), $\tau=0.1760$ (SB), $\tau=0.1733$ (TB), and $\tau=0.1700$ (TP) were selected.
(b) Bifurcation diagram for the four behaviors obtained numerically for the ODE system (\ref{eq_ODE_renormalized_expanded1}) (left) and its magnification (right). The horizontal axis denotes the parameter $\tau$ and the vertical axis denotes the maximal value of the pulse width $h=l_2-l_1$.
The solid (resp. broken) line represents stable (resp. unstable) solutions, and the squares and triangles on each solution branch denote Hopf and pitchfork bifurcation points, respectively.
For (a) and (b), the other parameters are fixed to $D=1.0$ and $\delta_0=0.001$.
}
 \label{pulse_dynamics_ode_homo}
\end{figure}
%

%
\section{Application to dynamics in bump-type heterogeneous media}
\subsection{Numerical simulation of PDE system}
So far, we have assumed that the function $\delta(x)$ in Eq.~(\ref{eq_PDE}) is constant with $\delta(x) \equiv \delta_0$, and hence $\Delta(l) \equiv \delta_0$ in Eq.~(\ref{eq_ODE_renormalized_expanded1}). In this section, we consider the heterogeneous case where $\delta(x)$ varies with $x$.
Although the function $\delta(x)$ may have an arbitrary form depending on the situation, the spatial heterogeneities of jump-, bump-, and periodic-type dynamics are ubiquitous in nature. Hence, they are often employed for numerical and analytical studies
\cite{Bode_1997}\cite{Ikeda_Ei_2010}\cite{Ikeda_Mimura_1989}\cite{Li_2003}\cite{Prat_Li_2003}\cite{Schutz_Bode_Purwins_1995}. 
In \cite{Nishiura_Teramoto_Yuan_Ueda_2007}\cite{Nishiura_Teramoto_Yuan_2012}\cite{Teramoto_Yuan_Baer_Nishiura_2009}\cite{Yadome_Nishiura_Teramoto_2014}\cite{Yuan_Teramoto_Nishiura_2007}, for instance, we considered pulse dynamics in excitable media with jump- and bump-type heterogeneities and clarified their mechanism from a bifurcation point of view.
Compared to the research on excitable systems, however, relatively few theoretical studies have examined the influence of heterogeneity on pulse behavior in bistable systems \cite{Heijster_Chen_Nishiura_Teramoto_2017}\cite{Heijster_Doelman_Kaper_Nishiura_Ueda_2011}. The geometric singular perturbation theory has been effectively used to clarify the existence and stability of stationary pulse solutions induced by heterogeneity, and in our previous paper \cite{Nishi_Nishiura_Teramoto_2013}, we considered a mixed ODE-PDE system  (\ref{eq_hybrid}) and applied a center manifold reduction to derive ODEs describing the pulse dynamics in the presence of a jump-type heterogeneity. 

In the present paper, we deal with a bump-type heterogeneity of the form
\begin{eqnarray}
  \displaystyle{ \delta(x)=\delta_0+\frac{\epsilon_0}{1+e^{\gamma (x-x_c+d_0/2)}}-\frac{\epsilon_0}{1+e^{-\gamma (x-x_c-d_0/2)}},
  }
  \label{eq_bump_pde}
\end{eqnarray}
and consider the situation in which a traveling pulse coming from the left infinity encounters this bump-type heterogeneity (Fig.~\ref{phase_diagram_bump}(a)).
In Eq.~(\ref{eq_bump_pde}), the parameter $\epsilon_0$ represents the bump height, $d_0$ is the bump width, and $\gamma$ is the steepness of the bump interface.
Through the numerical simulation of the PDEs in (\ref{eq_PDE}) with Eq.~(\ref{eq_bump_pde}),
we find four kinds of pulse behavior as $d_0$ and $\epsilon_0$ vary as bifurcation parameters: penetration (PEN), decomposition 1 (DEC1), rebound (REB), and decomposition 2 (DEC2). 
These are illustrated in Figure~\ref{phase_diagram_bump}(b). 
For PEN, the traveling pulse passes the bump region and travels on (Fig.~\ref{phase_diagram_bump}(b-i)). As the bump height is increased, the pulse exhibits either DEC1 or REB behavior, depending on the bump width $d_0$.
When $d_0$ is relatively small, DEC1 is observed, whereby the interfaces of the pulse separate after encountering the bump (Fig.~\ref{phase_diagram_bump}(b-ii)).
When $d_0$ is sufficiently large, the DEC1 behavior is no longer observed. Instead, the REB regime appears, whereby the pulse traveling to the right encounters the bump and turns back to the left (Fig.~\ref{phase_diagram_bump}(b-iii)).
For negative values of $\epsilon_0$, DEC2 is observed, whereby, as for DEC1, the two interfaces move away from each other and the pulse width diverges (Fig.~\ref{phase_diagram_bump}(b-iv)).

Note the difference between the two types of pulse decomposition behavior, DEC1 and DEC2. In the case of DEC1, the right interface can pass the bump region, whereas the left one quickly rebounds around the left edge of the bump region, resulting in pulse decomposition.
In the case of DEC2, the left interface enters the bump region, but cannot reach the right edge. Therefore, the left interface changes its direction of propagation after passing the left edge of the bump region.

\subsection{Mechanism for PEN-DEC and DEC-REB transitions}
The reduced four-dimensional ODEs shown in Section 2.2 are
\begin{eqnarray}
\left \{
\begin{array}{l}
\; \; \; \; \; \; \dot{l}_2 = r_2, \quad \; \; \; \; \dot{l}_1 = r_1, \\ 
~m_{0} \hspace{0.5mm}\dot{r}_2 = \sqrt{2} \hspace{0.5mm}(\tau_c-\tau) \hspace{0.5mm}r_2 - g_{3} \hspace{0.5mm}r^{3}_2
+(\hspace{0.5mm}G_0-G_1 \hspace{0.5mm}r_1\hspace{0.5mm}) \hspace{0.5mm}e^{-\frac{\;\;r_1+\phi(r_1)}{2 D}h}-\Delta_0(l_2), \\
~m_{0} \hspace{0.5mm}\dot{r}_1 = \sqrt{2} \hspace{0.5mm}(\tau_c-\tau) \hspace{0.5mm}r_1 - g_{3} \hspace{0.5mm}r^{3}_1
-(\hspace{0.5mm}G_0+G_1 \hspace{0.5mm}r_2\hspace{0.5mm}) \hspace{0.5mm}e^{-\frac{-r_2+\phi(r_2)}{2 D}h}+\Delta_0(l_1),
\end{array}
\right. \label{eq_ODE_renormalized_expanded1_2}
\end{eqnarray}
and these can also be applied to analytically investigate the pulse dynamics of the PDEs, as we saw in the previous section.
For the sake of analysis, we assume that the bump interfaces are sufficiently steep and take the limit $\gamma \to \infty$, so that the function $\delta_0(x)$ simplifies to
\begin{eqnarray}
\displaystyle{ 
\delta(x)=\left \{
\begin{array}{l}
\delta_0+\epsilon_0, \quad (\hspace{0.5mm}-d_0/2 \hspace{0.5mm}\hspace{0.5mm}\le x \le d_0/2\hspace{0.5mm})\\
\quad \delta_0 \quad  \quad \quad \quad {\rm otherwise} \hspace{0.5mm}
\end{array}
\right. }
\label{eq_bump_ode}
\end{eqnarray}
For the piecewise constant function (\ref{eq_bump_ode}), (\ref{eq_delta}) is solved analytically as follows. \vspace{0.2cm}

\noindent
{\bf Lemma 1.}
{\it Suppose that the function $\delta(x)$ is given by the piecewise constant function in (\ref{eq_bump_ode}). Then, a bounded $C^1$ solution $\Delta_0(x)$ to (\ref{eq_delta}) is given by
\begin{eqnarray}
\displaystyle{ 
\Delta_0(x)=\left \{
\begin{array}{l}
\displaystyle{ 
\quad \quad \delta_0+\frac{\epsilon_0}{2}\hspace{0.2mm} \left( \hspace{0.5mm}e^{\frac{d_0}{2\sqrt{D}}}-e^{-\frac{d_0}{2\sqrt{D}}} \hspace{0.5mm}\right) \hspace{0.5mm}e^{\frac{x}{\sqrt{D}}}, \quad \quad \left(\hspace{0.5mm}x \le -d_0/2 \hspace{0.5mm}\right) }  \vspace{0.1cm}\\
\displaystyle{ \delta_0+\epsilon_0-\frac{\epsilon_0}{2} \hspace{0.5mm}e^{-\frac{d_0}{2\sqrt{D}}} \left( \hspace{0.5mm}e^{\frac{x}{\sqrt{D}}}+ e^{-\frac{x}{\sqrt{D}}} \hspace{0.5mm}\right), \quad \; \left(\hspace{0.5mm}-d_0/2 < x \le d_0/2 \hspace{0.5mm}\right) }  \vspace{0.1cm}\\
\displaystyle{ 
\quad \quad \delta_0+\frac{\epsilon_0}{2}\hspace{0.2mm} \left(\hspace{0.5mm} e^{\frac{d_0}{2\sqrt{D}}}-e^{-\frac{d_0}{2\sqrt{D}}} \hspace{0.5mm}\right) \hspace{0.5mm}e^{-\frac{x}{\sqrt{D}}}. \; \quad \left(\hspace{0.5mm}d_0/2 \le x \hspace{0.5mm}\right) }  \vspace{0.1cm}\\
\end{array}
\right. }
\label{eq_Delta0}
\end{eqnarray}
\vspace{0.2cm}\\
{\bf Proof.}\hspace{0.1cm}
{\it As (\ref{eq_delta}) with (\ref{eq_bump_ode}) is an ODE with constant coefficients, it is readily solved as
\begin{eqnarray}
\displaystyle{ 
\Delta_0(x)=\left \{
\begin{array}{l}
\quad \quad \delta_0+ C_1 e^{\frac{x}{\sqrt{D}}}, \; \quad \quad \quad \quad \quad \quad (\hspace{0.5mm}x \le -d_0/2 \hspace{0.5mm})\\ 
\delta_0+\epsilon_0+C_2 e^{\frac{x}{\sqrt{D}}}+C_3 e^{-\frac{x}{\sqrt{D}}}, \quad (\hspace{0.5mm}-d_0/2 \hspace{0.5mm}\hspace{0.5mm}\le x \le d_0/2\hspace{0.5mm})\\
\quad \quad \delta_0+ C_4 e^{-\frac{x}{\sqrt{D}}}, \; \; \quad \quad \quad \quad \quad (\hspace{0.5mm}d_0/2 \le x \hspace{0.5mm})
\end{array}
\right. }\nonumber
\end{eqnarray}
where the unknown constants $C_i$ \hspace{0.1cm}(\hspace{0.5mm}i=1, 2, 3, 4\hspace{0.5mm}) are determined by the condition that both $\Delta_0(x)$ and $\Delta_0'(x)$ are continuous at $x=\pm \hspace{0.5mm}d_0/2$. Thus, we have $\displaystyle C_1=C_4=\epsilon_0\hspace{0.5mm}(\hspace{0.5mm}e^{\frac{d_0}{\sqrt{2D}}}-e^{-\frac{d_0}{\sqrt{2D}}}\hspace{0.5mm})/2$ and $\displaystyle C_2=C_3=-\epsilon_0\hspace{0.5mm}e^{-\frac{d_0}{\sqrt{2D}}}\hspace{0.5mm}/2$, which leads to Eq.~(\ref{eq_Delta0}).
} \hspace{0.0cm} \scalebox{1.0}{$\Box$} \vspace{0.2cm}

\rm
Now, we define penetration, rebound, and decomposition in terms of the solution orbit of the ODE system.\\

\noindent
{\bf Definition 2.} {\it When the orbit of (\ref{eq_ODE_renormalized_expanded1_2}) with (\ref{eq_Delta0}) starting from ${\rm TP}^{\hspace{0.5mm}+}$ converges to ${\rm TP}^{\hspace{0.5mm}+}$ as $t \to +\infty$, we refer to the dynamics as penetration.}\\
\rm \noindent
{\bf Definition 3.} {\it When the orbit of (\ref{eq_ODE_renormalized_expanded1_2}) with (\ref{eq_Delta0}) starting from ${\rm TP}^{\hspace{0.5mm}+}$ converges to ${\rm TP}^{\hspace{0.5mm}-}$ as $t \to +\infty$, we refer to the dynamics as rebound.}\\
\rm \noindent
{\bf Definition 4.} {\it When the orbit of (\ref{eq_ODE_renormalized_expanded1_2}) with (\ref{eq_Delta0}) starting from ${\rm TP}^{\hspace{0.5mm}+}$ converges to the solution consisting of ${\rm TF}_{2}^{\hspace{0.5mm}+}$ and ${\rm TF}_{1}^{\hspace{0.5mm}-}$ as $t \to +\infty$, we refer to the dynamics as decomposition.}\\

\rm \noindent
{\bf Remark}
{\it For the definitions of ${\rm TF}_{2}^{\hspace{0.5mm}+}$ and ${\rm TF}_{1}^{\hspace{0.5mm}-}$ corresponding to the right- and left-moving traveling front solutions, see Definitions 5 and 6, respectively. Note that, for the decomposition behavior, the distance between the two interfaces diverges as $t \to +\infty$, and hence $h \to \infty$. As mentioned at the end of Section 2.2, the limit $h \to \infty$ decouples the four-dimensional ODEs (\ref{eq_ODE_renormalized_expanded1_2}) into two ODE systems, each describing the separate motions of the left and right interfaces. ${\rm TF}_{2}^{\hspace{0.5mm}+}$ and ${\rm TF}_{1}^{\hspace{0.5mm}-}$ in Definition 4 are solutions to each of the two ODE systems, which correspond to traveling front solutions in the original PDE system.}
\\

\rm
Numerical simulations of Eq.~(\ref{eq_ODE_renormalized_expanded1_2}) with Eq.~(\ref{eq_Delta0}) indicate that the reduced ODEs provide a good qualitative reproduction of the pulse dynamics in the original PDE system (Fig.~\ref{phase_diagram_bump}(c)). In the following, we utilize the ODEs to investigate the mechanism for the PEN-DEC1, DEC1-REB, and PEN-DEC2 transitions from the viewpoint of bifurcation theory.
In particular, we characterize the difference between DEC1 and DEC2, which exhibit similar asymptotic behavior as shown in Figures~\ref{phase_diagram_bump}(c-ii) and \ref{phase_diagram_bump}(c-iv).

To this end, the parameters were set very close to the phase boundaries; $d_0$ was fixed and $\epsilon_0$ was varied (Fig.~\ref{phase_diagram_bump_scattor}(a)). The results reveal that the pulse interfaces exhibit different kinds of transient behavior right after the pulse encounters the heterogeneity.
For parameters close to the PEN-DEC1 boundary, the left interface bounces off the bump heterogeneity and transiently moves with almost a constant velocity to the left before eventually traveling either to the left or to the right, resulting in the DEC1 or PEN behavior, respectively (Fig.~\ref{phase_diagram_bump_scattor}(a-i)).
For those close to the DEC1-REB boundary, however, the right interface transiently moves to the right (Fig.~\ref{phase_diagram_bump_scattor}(a-ii)).
In contrast, for the PEN-DEC2 transition, the left interface stays near the right edge of the bump heterogeneity for a certain period of time (Fig.~\ref{phase_diagram_bump_scattor}(a-iii)).

These numerical findings lead us to infer that one of the interfaces transiently approaches some unstable traveling front solution for the PEN-DEC1 and DEC1-REB transitions, whereas it approaches some unstable stationary front solution for the PEN-DEC2 transition. This can be confirmed analytically using the ODEs in (\ref{eq_ODE_renormalized_expanded1_2}).
First, for the PEN-DEC1 transition, the spatiotemporal plot in Figure \ref{phase_diagram_bump_scattor}(a-i) suggests that the left interface moves away from both the right interface and the heterogeneity immediately after it bounces off, and that the influence of these two is negligible. Therefore, we consider the following ODEs for the motion of the left interface by setting $h \to \infty$ and $\Delta_0(l_1)=\delta_0$ in Eq.~(\ref{eq_ODE_renormalized_expanded1_2}):
\begin{eqnarray}
\displaystyle{ 
	\left \{
\begin{array}{l}
\quad \quad \hspace{0.5mm}\dot{l}_1 = r_1,  \\ 
~m_{0}\hspace{0.5mm}\dot{r}_1 = \sqrt{2} \hspace{0.5mm}(\tau_c-\tau) ~r_1 - g_{3} \hspace{0.5mm}r^{3}_1 \hspace{0.5mm}+\delta_0\hspace{0.5mm}.
\end{array}
\right. }
\label{eq_ODE_left1}
\end{eqnarray}

The solution with $r_1= const$ in (\ref{eq_ODE_left1}) corresponds to the traveling front solution in PDE system (\ref{eq_PDE}).\\

\noindent
{\bf Definition 5.}
{\it  The solution {\rm (}\hspace{0.5mm}$l_1,\hspace{0.5mm}r_1$\hspace{0.5mm}{\rm )} to Eq.~(\ref{eq_ODE_left1}) with $\dot{r}_1\hspace{0.2mm}=\hspace{0.2mm}0$ is defined as ${\rm TF}^{\hspace{0.5mm}+}_{1} \hspace{0.8mm}{\rm (} \hspace{0.5mm}r_1=\bar{r}^{\hspace{0.3mm}+}_1\hspace{0.5mm}{\rm )}$, ${\rm TF}^{\hspace{0.5mm}-}_{1} \hspace{0.8mm}{\rm (} \hspace{0.5mm}r_1=\bar{r}^{\hspace{0.3mm}-}_1\hspace{0.5mm}{\rm )}$, and\hspace{1.5mm}${\rm TF}^{\hspace{0.5mm}0}_{1} \hspace{0.8mm}{\rm (} \hspace{0.5mm}r_1=\bar{r}^{\hspace{0.3mm}0}_1\hspace{0.5mm}{\rm )}$, where $\bar{r}^{\hspace{0.3mm}-}_1<\bar{r}^{\hspace{0.3mm}0}_1<0<\bar{r}^{\hspace{0.3mm}+}_1$.}\\

The equation for $\dot{r}_1=0$ in the second equation in (\ref{eq_ODE_left1}) is solved as follows. \vspace{0.2cm}

\noindent
{\bf Lemma 2.}
{\it
Assume that $0< \hspace{0.5mm}\delta_0 \hspace{0.5mm} \ll \hspace{0.5mm} 1$,\hspace{1.0mm}$g_3>0$, and $\tau < \tau_c$. Then, the cubic equation for $r_1$,
\begin{eqnarray}
   \sqrt{2} \hspace{0.5mm}(\tau_c-\tau) ~r_1 - g_{3} \hspace{0.5mm}r^{3}_1 \hspace{0.5mm}+\delta_0\hspace{0.5mm}=0, \label{eq_cubic}
\end{eqnarray}
has three roots, $r_1=r^{\pm}_1$ and $\hspace{0.5mm}r^{0}_1$. These are approximately given by
\begin{eqnarray}
\displaystyle{ 
	\left \{
\begin{array}{l}
\hspace{0.3cm}\displaystyle{ r^{\pm}_1=\pm \sqrt{\frac{\sqrt{2}(\tau_c-\tau)}{g_3}}+\frac{1}{2\sqrt{2}(\tau_c-\tau)} ~\delta_0 +\mathcal{O}(\delta^2_0), } \vspace{0.1cm}\\
\hspace{0.3cm}\displaystyle{ r^{0}_1=-\frac{1}{\sqrt{2}(\tau_c-\tau)}~\delta_0+\mathcal{O}(\delta^3_0).}
\end{array}
\right. }
\label{eq_front_velocity}
\end{eqnarray}
}
\vspace{0.2cm}\\
{\bf Proof.}\hspace{0.1cm}
{\it By assuming a solution of the form $r_1\hspace{0.5mm}=\hspace{0.5mm}a_0\hspace{0.5mm}+\hspace{0.5mm}a_1\hspace{0.5mm}\delta_0+a_2\hspace{0.5mm}\delta^2_0+ \cdots$ and substituting this into (\ref{eq_cubic}), we have
\begin{eqnarray}
\begin{array}{lll}
\mathcal{O}(\; \delta_0^0 \;): \quad \sqrt{2}\hspace{0.5mm}(\tau_c-\tau)\hspace{0.5mm}a_0-g_3\hspace{0.5mm}a^3_0=0\hspace{0.5mm}, \vspace{0.2cm} \\
\mathcal{O}(\; \delta_0^1 \;): \quad \sqrt{2}\hspace{0.5mm}(\tau_c-\tau)\hspace{0.5mm}a_1-3\hspace{0.5mm}g_3\hspace{0.5mm}a^2_0\hspace{0.5mm}a_1+1=0\hspace{0.5mm}, \vspace{0.2cm} \\
\mathcal{O}(\; \delta_0^2 \;): \quad \sqrt{2}\hspace{0.5mm}(\tau_c-\tau)\hspace{0.5mm}a_2-3\hspace{0.5mm}g_3 \left(\hspace{0.5mm} a_0\hspace{0.5mm}a^2_1+a_0^2\hspace{0.5mm}a_2\hspace{0.5mm}\right)=0\hspace{0.5mm}. 
\end{array} \nonumber
\end{eqnarray}
From the first equation, $a_0$ is solved as $a_0=0,\hspace{0.5mm}\pm \hspace{0.5mm}\sqrt{\sqrt{2}\hspace{0.5mm}(\tau_c-\tau)/g_3}$. For $a_0=0$, the second and the third equations lead to $a_1=-1/\{ \hspace{0.5mm}\sqrt{2}\hspace{0.5mm}(\tau_c-\tau) \}$ and $a_2=0$. Similarly, for $a_0=\pm \hspace{0.5mm}\sqrt{\sqrt{2}\hspace{0.5mm}(\tau_c-\tau)/g_3}$, $a_1$ and $a_2$ are solved as $a_1=1/\{ \hspace{0.5mm}2\sqrt{2}\hspace{0.5mm}(\tau_c-\tau) \}$ and $a_2=\mp \hspace{0.5mm}3\hspace{0.5mm}g_3/\{ \hspace{0.5mm}16\sqrt{2}\hspace{0.5mm}(\tau_c-\tau)^3 \}\hspace{0.5mm}\sqrt{\sqrt{2}\hspace{0.5mm}(\tau_c-\tau)/g_3}$. Collecting the $a_i$ terms up to $\mathcal{O}(\delta_0^2)$, we have the three solutions in (\ref{eq_front_velocity}).
} \hspace{0.0cm} \scalebox{1.0}{$\Box$} \vspace{0.2cm}

The solutions $r_1=\hspace{0.5mm}r^{+}_1, \hspace{0.5mm}r^{-}_1 \hspace{0.5mm}$, and $\hspace{1.5mm}r^{0}_1$ approximate the velocities of ${\rm TF}^{\hspace{0.5mm}+}_{1}$, ${\rm TF}^{\hspace{0.5mm}-}_{1}$, and ${\rm TF}^{\hspace{0.5mm}0}_{1}$, respectively. Similarly, the solution that corresponds to the traveling front for the right interface, whose motion is governed by
\begin{eqnarray}
\displaystyle{ 
	\left \{
\begin{array}{l}
\quad \quad \hspace{0.5mm}\dot{l}_2 = r_2,  \\ 
~m_{0}\hspace{0.5mm}\dot{r}_2 = \sqrt{2} \hspace{0.5mm}(\tau_c-\tau) ~r_2 - g_{3} \hspace{0.5mm}r^{3}_2 \hspace{0.5mm}-\delta_0\hspace{0.5mm},
\end{array}
\right. }
\label{eq_ODE_left2_0}
\end{eqnarray}
is defined as follows.\\

\noindent
{\bf Definition 6.}
{\it The solution {\rm (}\hspace{0.5mm}$l_2,\hspace{0.5mm}r_2$\hspace{0.5mm}{\rm )} to Eq.~(\ref{eq_ODE_left2_0}) with $\dot{r}_2\hspace{0.2mm}=\hspace{0.2mm}0$ is defined as ${\rm TF}^{\hspace{0.5mm}+}_{2} \hspace{0.8mm}{\rm (} \hspace{0.5mm}r_2=\bar{r}^{\hspace{0.3mm}+}_2\hspace{0.5mm}{\rm )}$, ${\rm TF}^{\hspace{0.5mm}-}_{2} \hspace{0.8mm}{\rm (} \hspace{0.5mm}r_2=\bar{r}^{\hspace{0.3mm}-}_2\hspace{0.5mm}{\rm )}$, and\hspace{1.5mm}${\rm TF}^{\hspace{0.5mm}0}_{2} \hspace{0.8mm}{\rm (} \hspace{0.5mm}r_2=\bar{r}^{\hspace{0.3mm}0}_2\hspace{0.5mm}{\rm )}$, where $\bar{r}^{\hspace{0.3mm}-}_2<0<\bar{r}^{\hspace{0.3mm}0}_2<\bar{r}^{\hspace{0.3mm}+}_2$.
}\\

\rm \noindent
The approximate value of $r_2$ in Definition 6 is obtained by inverting the sign of $\delta_0$ to $-\delta_0$ in Eq.~(\ref{eq_ODE_left1}).\\

\noindent
{\bf Corollary 2.}
{\it
Assume that $0< \hspace{0.5mm}\delta_0 \hspace{0.5mm} \ll \hspace{0.5mm} 1$,\hspace{1.0mm}$g_3>0$, and $\tau < \tau_c$. Then, the cubic equation for $r_2$,
\begin{eqnarray}
   \sqrt{2} \hspace{0.5mm}(\tau_c-\tau) ~r_2 - g_{3} \hspace{0.5mm}r^{3}_2 \hspace{0.5mm}-\delta_0\hspace{0.5mm}=0, \label{eq_cubic2}
\end{eqnarray}
has three roots, $r_2=r^{\pm}_2$ and $\hspace{0.5mm}r^{0}_2$. These are approximately given by
\begin{eqnarray}
\displaystyle{ 
	\left \{
\begin{array}{l}
\hspace{0.3cm}\displaystyle{ r^{\pm}_2=\pm \sqrt{\frac{\sqrt{2}(\tau_c-\tau)}{g_3}}-\frac{1}{2\sqrt{2}(\tau_c-\tau)} ~\delta_0 +\mathcal{O}(\delta^2_0), } \vspace{0.1cm}\\
\hspace{0.3cm}\displaystyle{ r^{0}_2=\frac{1}{\sqrt{2}(\tau_c-\tau)}~\delta_0+\mathcal{O}(\delta^3_0).}
\end{array}
\right. }
\label{eq_front_velocity2}
\end{eqnarray}
}\\

\rm \noindent
Figure \ref{phase_diagram_bump_scattor}(b-i) shows the numerical behavior of $r_1(t)$ for Eq.~(\ref{eq_ODE_left1}) with the parameters corresponding to the PEN-DEC1 phase boundary in Figure~\ref{phase_diagram_bump_scattor}(a-i). The propagation velocity $r_1$ of the left interface remains almost constant at $r_1^0$ for a certain period of time, before eventually converging to either $r_1^+$ or $r_1^-$, depending on the parameter change. This corresponds to the situation where the solution orbit of (\ref{eq_ODE_left1}) approaches ${\rm TF}^{\hspace{0.5mm}0}_{1}$ and converges to either ${\rm TF}^{\hspace{0.5mm}+}_{1}$ or ${\rm TF}^{\hspace{0.5mm}-}_{1}$.
Therefore, it may be concluded that the unstable traveling front solution with $r_1^0$ (i.e., ${\rm TF}^{\hspace{0.5mm}0}_{1}$) plays a crucial role in changing the destination of the left interface, leading to the two different outputs of PEN and DEC1.
In our previous work \cite{Nishiura_Teramoto_Ueda_2003}\cite{Nishiura_Teramoto_Ueda_2005}, we showed that unstable patterns, named scattors, and the local dynamics around them underlie a variety of outputs in the collision process of localized moving patterns. The orbital flow immediately after the collision travels along one of the scattors' unstable directions.
A similar argument holds for the DEC1-REB transition, and an unstable traveling front solution ${\rm TF}^{\hspace{0.5mm}0}_{2}$ plays the role of the scattor in this case, too.
The above discussions lead us to infer the following plausible scenario for the PEN-DEC1 and DEC1-REB transitions for the bump-up case $\epsilon_0 > 0$.\\

\noindent
{\bf Proposition 4.} There exist positive constants $\epsilon^{(1)}_0$ and $\epsilon^{(2)}_0$ (\hspace{0.5mm}$\epsilon^{(2)}_0 > \epsilon^{(1)}_0 > 0$\hspace{0.5mm}) such that the solution orbit of (\ref{eq_ODE_renormalized_expanded1_2}) with (\ref{eq_Delta0}) starting from ${\rm TP}^{\hspace{0.5mm}+}$ converges to the solution consisting of ${\rm TF}^{\hspace{0.5mm}0}_{1}$ and ${\rm TF}^{\hspace{0.5mm}+}_{2}$ at $\epsilon_0=\epsilon^{(1)}_0$, and to that consisting of ${\rm TF}^{\hspace{0.5mm}-}_{1}$ and ${\rm TF}^{\hspace{0.5mm}0}_{2}$ at $\epsilon_0=\epsilon^{(1)}_0$.\\

\noindent
Note that the constants $\epsilon^{(1)}_0$ and $\epsilon^{(2)}_0$ in Proposition 4 correspond to the values of $\epsilon_0$ for the PEN-DEC1 and DEC1-REB boundaries in Figure~\ref{phase_diagram_bump_scattor}(a-i,ii), respectively. The unstable manifolds of the unstable solutions ${\rm TF}^{\hspace{0.5mm}0}_{1}$ and ${\rm TF}^{\hspace{0.5mm}0}_{2}$ near $\epsilon_0=\epsilon^{(1)}_0$ and $\epsilon_0=\epsilon^{(2)}_0$ determine the destination of the solution orbit after encountering the bump region, resulting in penetration, decomposition, or rebound behavior.

For the PEN-DEC2 transition for the bump-down case $\epsilon_0 < 0$, Figure~\ref{phase_diagram_bump_scattor}(a-iii) shows that the left interface transiently stays close to the right edge of the bump, meaning that the influence from the heterogeneity may not be neglected in this case. Therefore, we consider the following ODEs for the motion of the left interface:
\begin{eqnarray}
\displaystyle{ 
	\left \{
\begin{array}{l}
\quad \quad \hspace{0.5mm}\dot{l}_1 = r_1,  \\ 
~m_{0}\hspace{0.5mm}\dot{r}_1 = \sqrt{2} \hspace{0.5mm}(\tau_c-\tau) ~r_1 - g_{3} \hspace{0.5mm}r^{3}_1 \hspace{0.5mm}+\Delta_0(l_1)\hspace{0.5mm}.
\end{array}
\right. }
\label{eq_ODE_left2}
\end{eqnarray}
For the explicit form of $\Delta_0(x)$ in (\ref{eq_Delta0}), the stationary solution to Eq.~(\ref{eq_ODE_left2}) is solved as follows.\vspace{0.2cm}

\noindent
{\bf Lemma 3.}
{\it
  For $\delta_0>0$ and $\epsilon_0 \le -\delta_0/\hspace{0.5mm}(\hspace{0.5mm}1-e^{-\frac{d_0}{2\sqrt{D}}}\hspace{0.5mm})$, a stationary solution $(\hspace{0.5mm}l^*_1,\hspace{0.5mm}r^*_1\hspace{0.5mm})$\hspace{0.3cm}$(\hspace{0.5mm}l^*_1 \ge 0\hspace{0.5mm})$ to Eq.~(\ref{eq_ODE_left2}) is given by $r^*_1=0$ and
\begin{eqnarray}
\displaystyle{ 
\hspace{0.3cm} l^*_1=\left \{
\begin{array}{l}
\hspace{0.5cm}\sqrt{D}\hspace{0.5mm} \log \left\{\hspace{0.5mm} \displaystyle{ 
  \frac{\delta_0+\epsilon_0}{\epsilon_0 z_0}-\sqrt{\left( \frac{\delta_0+\epsilon_0}{\delta_0 z_0}\right)^2-1} } \hspace{0.5mm}\right\}, \vspace{0.1cm}\\
\hspace{2.9cm} \left(\hspace{1.0mm} -2\hspace{0.2mm}\delta_0/\hspace{0.5mm}(\hspace{0.5mm}1-z^2_0\hspace{0.5mm}) \hspace{0.5mm} \le \hspace{0.5mm}\epsilon_0 \hspace{0.5mm}\le \hspace{0.5mm}-\delta_0/\hspace{0.5mm}(\hspace{0.5mm}1-z_0\hspace{1.0mm}) \hspace{0.5mm}\right) \vspace{0.5cm}\\
\hspace{0.8cm}-\sqrt{D}\hspace{0.5mm} \log \left\{\hspace{0.5mm} \displaystyle{ -\frac{2 z_0 \delta_0}{\epsilon_0 (1-z^2_0)} } \hspace{0.5mm}\right\}, \vspace{0.1cm} \hspace{0.8cm} \left(\hspace{1.0mm} \epsilon_0 \hspace{0.5mm}\le \hspace{0.5mm} -2\hspace{0.2mm}\delta_0/\hspace{0.5mm}(\hspace{0.5mm}1-z^2_0\hspace{1.0mm}) \hspace{0.5mm}\right) \\
\end{array}
\right. } \label{eq_front_position}
\end{eqnarray}
where $z_0=e^{-\frac{d_0}{2\sqrt{D}}}$. Furthermore, the stationary solution is unstable for $\tau<\tau_c$.
}
\vspace{0.2cm}\\
{\bf Proof.}\hspace{0.1cm}
{\it
  Setting $\dot{l}_1=\dot{r}_1=0$ and $(\hspace{0.5mm}l_1,\hspace{0.5mm}r_1\hspace{0.5mm})=(\hspace{0.5mm}l^*_1,\hspace{0.5mm}r^*_1\hspace{0.5mm})$ yield $r^*_1=0$ and $\Delta_0(\hspace{0.5mm}l^*_1\hspace{0.5mm})=0$. As $\Delta_0(\hspace{0.5mm}x\hspace{0.5mm})$ in (\ref{eq_Delta0}) is symmetric about $x=0$, we assume $l^*_1 \ge 0$. Then, $\Delta_0(\hspace{0.5mm}l^*_1\hspace{0.5mm})=0 \hspace{2.5mm}(\hspace{0.5mm}l^*_1 \ge 0\hspace{0.5mm})$ is solved as  
\begin{eqnarray}
\displaystyle{ 
\hspace{0.5cm} z^*_1=\left \{
\begin{array}{l}
\hspace{0.3cm}\hspace{0.5mm} \displaystyle{ \frac{\delta_0+\epsilon_0}{\delta_0 z_0}-\sqrt{\left( \frac{\delta_0+\epsilon_0}{\epsilon_0 z_0}\right)^2-1} } \hspace{0.5mm}, \hspace{0.5cm} \left(\hspace{0.5mm}z_0 \le z^*_1 \le 1 \hspace{0.5mm}\right) \vspace{0.3cm}\\
\hspace{1.2cm} \displaystyle{ -\frac{2 z_0 \delta_0}{\epsilon_0 (1-z^2_0)} } \hspace{0.5mm}, \hspace{2.0cm} \left(\hspace{0.5mm} 0 < z^*_1 \le z_0 \hspace{0.5mm}\right) \\
\end{array}
\right. } \label{eq_front_position2}
\end{eqnarray}
where $\displaystyle z^*_1:=e^{-\frac{l^*_1}{\sqrt{D}}}$ and $\displaystyle z_0:=e^{-\frac{d_0}{2\sqrt{D}}}$. Note that the conditions $\delta_0>0$ and $\epsilon_0 \le -\delta_0/\hspace{0.5mm}(\hspace{0.5mm}1-z_0 \hspace{0.5mm})$ ensure the inequality
\begin{eqnarray}
\displaystyle \left( \frac{\delta_0+\epsilon_0}{\delta_0 z_0} \right)^2-1 \hspace{0.5mm}>\hspace{0.5mm}0 \nonumber
\end{eqnarray}
holds in the first equation in (\ref{eq_front_position2}). Furthermore, the interval $\hspace{0.5mm}z_0 \le z^*_1 \le 1 \hspace{0.5mm}$ with
\begin{eqnarray}
\displaystyle
 \hspace{0.5cm} z^*_1= \frac{\delta_0+\epsilon_0}{\delta_0 z_0}-\sqrt{\left( \frac{\delta_0+\epsilon_0}{\delta_0 z_0}\right)^2-1} \hspace{0.5mm} \nonumber
\end{eqnarray}
in the first equation can be rewritten in terms of $\epsilon_0$ as $\epsilon_0 \hspace{0.5mm}\ge \hspace{0.5mm}-2\hspace{0.2mm}\delta_0/\hspace{0.5mm}(\hspace{0.5mm}1-z^2_0\hspace{0.5mm})$. Similarly, the interval $\hspace{0.5mm} z^*_1 \le z_0 \hspace{0.5mm}$ in the second equation is rewritten as $\epsilon_0 \hspace{0.5mm}\le \hspace{0.5mm}-2\hspace{0.2mm}\delta_0/\hspace{0.5mm}(\hspace{0.5mm}1-z^2_0\hspace{0.5mm})$. Finally, by the relation $l^*_1=\hspace{0.5mm}-\sqrt{D}\hspace{0.5mm} \log z^*_1$, we obtain the equations for $l^*_1$.

The stability of the stationary solution is determined by the eigenvalues of the Jacobi matrix
\begin{eqnarray}
  \displaystyle
J_f=\left (
\begin{array}{cc}
0 & 1  \\
p_0  &  \sqrt{2} \hspace{0.5mm}(\hspace{0.5mm}\tau_c-\tau\hspace{0.5mm})\hspace{0.5mm}/\hspace{0.5mm}m_0\\
\end{array}
\right ) \nonumber
\end{eqnarray}
where 
\begin{eqnarray}
\displaystyle{ 
\hspace{0.1cm} p_0 :=\Delta'_0(l^*_1)/\hspace{0.5mm}m_0=\left \{
\begin{array}{l}
\hspace{0.2cm}\displaystyle{ -\frac{1}{m_0\sqrt{D}}\hspace{0.5mm}\sqrt{(\delta_0+\epsilon_0)^2-(\epsilon_0 z_0)^2}, \hspace{1.0cm} \left(\hspace{0.5mm}z_0 \le z^*_1 \le 1 \hspace{0.5mm}\right) } \vspace{0.1cm}\\
\hspace{1.8cm}\displaystyle{ -\frac{1}{m_0\sqrt{D}} \hspace{0.8mm}\delta_0. \hspace{2.7cm}\left(\hspace{0.5mm} 0 < z^*_1 \le z_0 \hspace{0.5mm}\right) }
\end{array}
\right. } \nonumber
\end{eqnarray}
As ${\rm tr} \hspace{0.5mm} J_f=\sqrt{2} \hspace{0.5mm}(\hspace{0.5mm}\tau_c-\tau\hspace{0.5mm})/m_0>0$ and ${\rm det} \hspace{0.5mm} J_f=-p_0>0$, we find that the eigenvalues are both positive for $\tau<\tau_c$, and hence the stationary solution is unstable.
} \hspace{0.0cm} \scalebox{1.0}{$\Box$} \vspace{0.2cm}

Note that $l_1=l^*_1$ corresponds to the location of the stationary front solution. This is unstable for $\tau<\tau_c$, the parameter range we now consider.\\

\noindent
{\bf Definition 7.} {\it
  The stationary solution ${\rm (}\hspace{0.5mm}l_1,\hspace{0.5mm}r_1\hspace{0.5mm}{\rm )}={\rm (}\hspace{0.5mm}l^*_1,\hspace{0.5mm}0\hspace{0.5mm}{\rm )}$ to Eq.~(\ref{eq_ODE_left2}) with (\ref{eq_Delta0}) is defined as ${\rm SF}_{1}$, where $l^*_1$ is given in Eq.~(\ref{eq_front_position}).}\\

\noindent
Figure \ref{phase_diagram_bump_scattor}(b-ii) shows the numerical behavior of $l_1(t)$ for Eq.~(\ref{eq_ODE_left1}) with parameters corresponding to the PEN-DEC2 phase boundary in Figure~\ref{phase_diagram_bump_scattor}(a-iii). The location $l_1$ of the left interface remains almost constant around $l_1^*$ for a certain period of time, before eventually moving either to the left or right, depending on the parameter change.
Therefore, in the PEN-DEC2 transition, the unstable stationary front solution with $l_1=l^*_1$ plays the role of the scattor. Similar to Proposition 4, we propose the following plausible scenario for the PEN-DEC2 transition for the bump-down case $\epsilon_0 < 0$.\\

\noindent
{\bf Proposition 5.} {\it There exists a constant $\epsilon^{(3)}_0 <0$ such that, at $\epsilon_0=\epsilon^{(3)}_0$, the solution orbit of (\ref{eq_ODE_renormalized_expanded1_2}) with (\ref{eq_Delta0}) starting from ${\rm TP}^{\hspace{0.5mm}+}$ converges to the solution consisting of ${\rm SF_{1}}$ and ${\rm TF}^{\hspace{0.5mm}+}_{2}$.}\\

\rm
The constant $\epsilon^{(3)}_0$ in Proposition 5 corresponds to the value of $\epsilon_0$ for the PEN-DEC2 boundary in Figure~\ref{phase_diagram_bump_scattor}(a-iii), and the unstable manifold of ${\rm SF}_{1}$ near $\epsilon_0=\epsilon^{(3)}_0$ determines the destination of the solution orbit after encountering the bump region, resulting in either penetration or decomposition behavior.

In summary, all the transitions observed in Figure \ref{phase_diagram_bump}(c) can be understood as changes in the orbital behavior of the unstable solutions.
In the framework of the four-dimensional dynamical system (\ref{eq_ODE_renormalized_expanded1_2}), it seems plausible that the unstable manifolds emanating from the unstable traveling and stationary front solutions are connected to the stable traveling front solutions that propagate to $x \to \pm \infty$. This could be confirmed by more careful numerical analysis, as for the collision dynamics of localized patterns \cite{Nishiura_Teramoto_Ueda_2003}\cite{Nishiura_Teramoto_Ueda_2005}.
\begin{figure}[ht!]
   \centering
 \includegraphics[width=10cm]{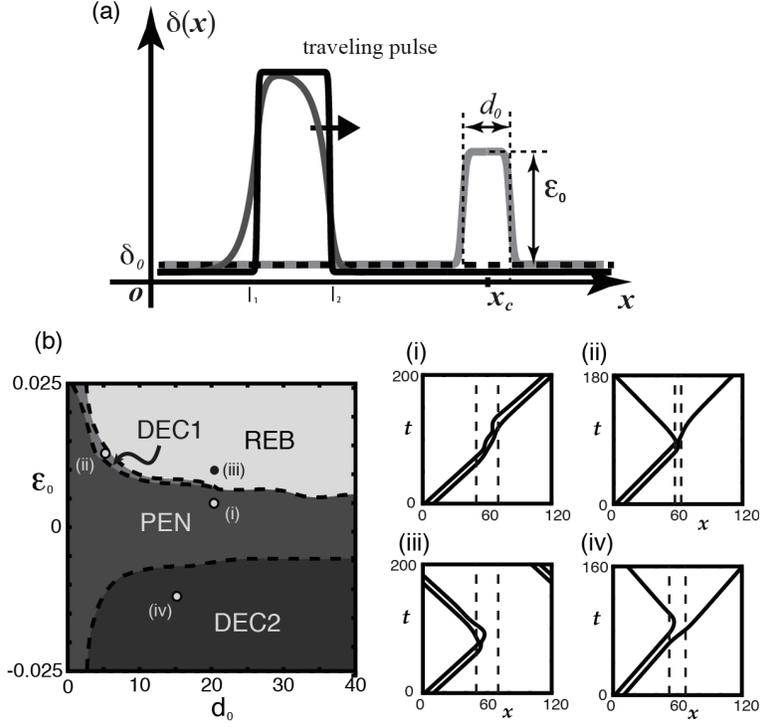}
 \caption{
   (a) Schematic figure for the pulse dynamics in a bump-type heterogeneous medium. The traveling pulse coming from the left infinity encounters the bump-type heterogeneity of width $d_0$ and height $\epsilon_0$. The black and gray lines represent the profiles of the $u$ and $v$ components, respectively.
   (b) Phase diagram obtained numerically for the PDEs in (\ref{eq_PDE}) with the heterogeneity function (\ref{eq_bump_pde}) by varying the bump width $d_0$ and bump height $\epsilon_0$. The right side shows typical spatiotemporal plots for the four kinds of behavior where the locations of the pulse interfaces are plotted, and 
   the two vertical broken lines denote the edges of the bump heterogeneity at $x=x_c+\pm \hspace{0.5mm}d_0/2$ with $x_c=60$. The parameters for each spatiotemporal plot correspond to the points labeled (i)--(iv) in the phase diagram.
   For the numerical simulation, we set $\tau=0.1700$ and the other parameters were the same as in Figure~\ref{pulse_dynamics_pde_homo}.
 }
   \label{phase_diagram_bump}
\end{figure}
\setcounter{figure}{2}
\begin{figure}[ht!]
   \centering
 \includegraphics[width=10cm]{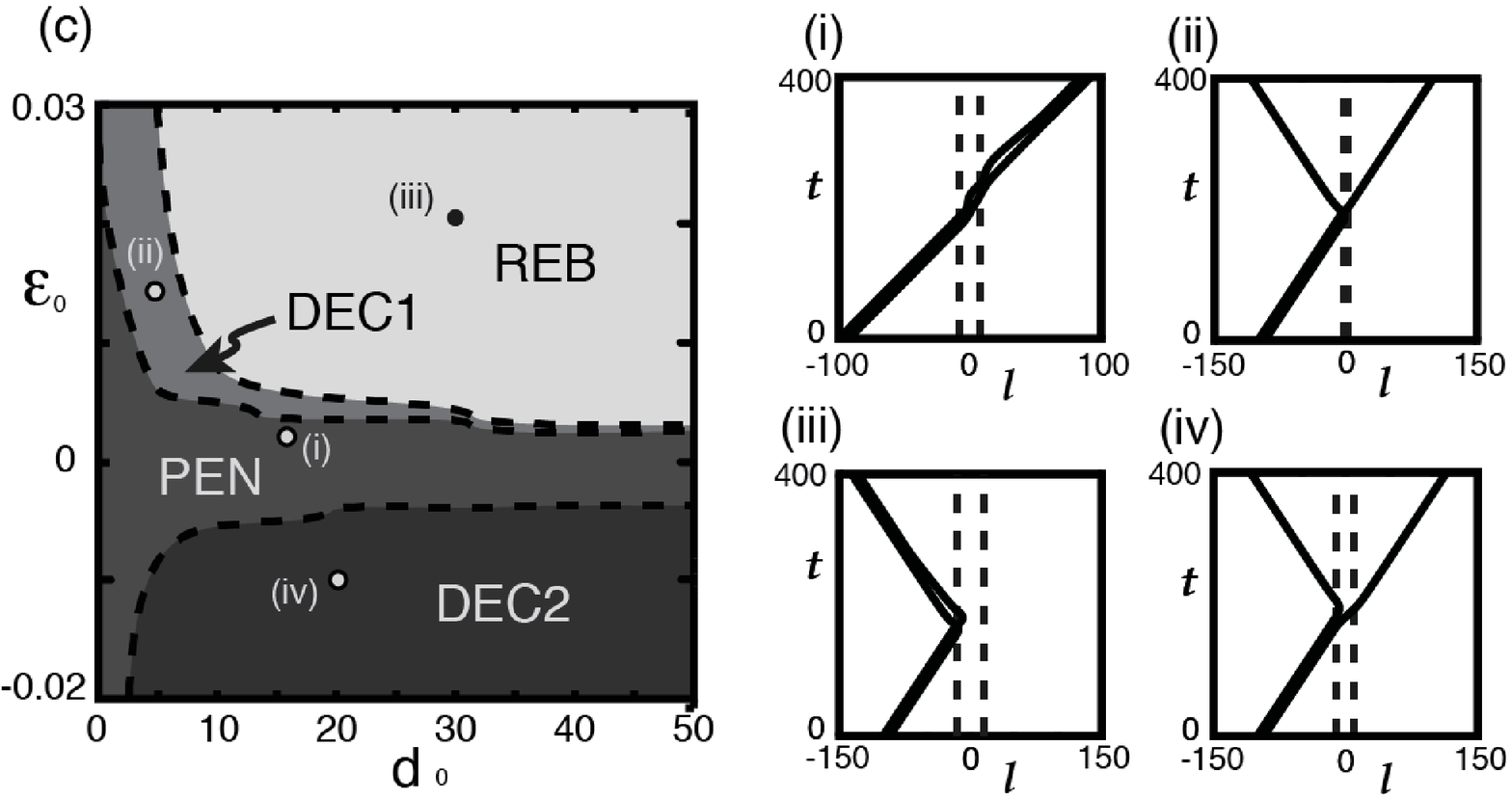}
 \caption{({\it continued}\hspace{0.2mm})
   (c) Phase diagram obtained numerically for the ODEs in (\ref{eq_ODE_renormalized_expanded1_2}) with the function (\ref{eq_Delta0}). The right side shows typical spatiotemporal plots for the four kinds of behavior where $l=l_2(t)$ and $l=l_1(t)$ are plotted, and the two vertical broken lines denote the edges of the bump heterogeneity at $l=\pm \hspace{0.5mm}d_0/2$. The parameters for each spatiotemporal plot correspond to the points labeled (i)--(iv) in the phase diagram.
   For the numerical simulation, we set $\tau=0.1700$, $\delta_0=0.001$, and the coefficients were computed as $m_0=3/(\hspace{0.5mm}16\sqrt{D}\hspace{0.5mm})$, $\tau_c=1/(\hspace{0.5mm}4\sqrt{2D}\hspace{0.5mm})$, $g_3=1/(\hspace{0.5mm}32(\sqrt{D})^3\hspace{0.5mm})$, $G_0=1/2$, and $G_1=1/(\hspace{0.5mm}4\sqrt{D}\hspace{0.5mm})$ with $D=1.0$.
 }
\end{figure}
\begin{figure}[ht!]
   \centering
 \includegraphics[width=12cm]{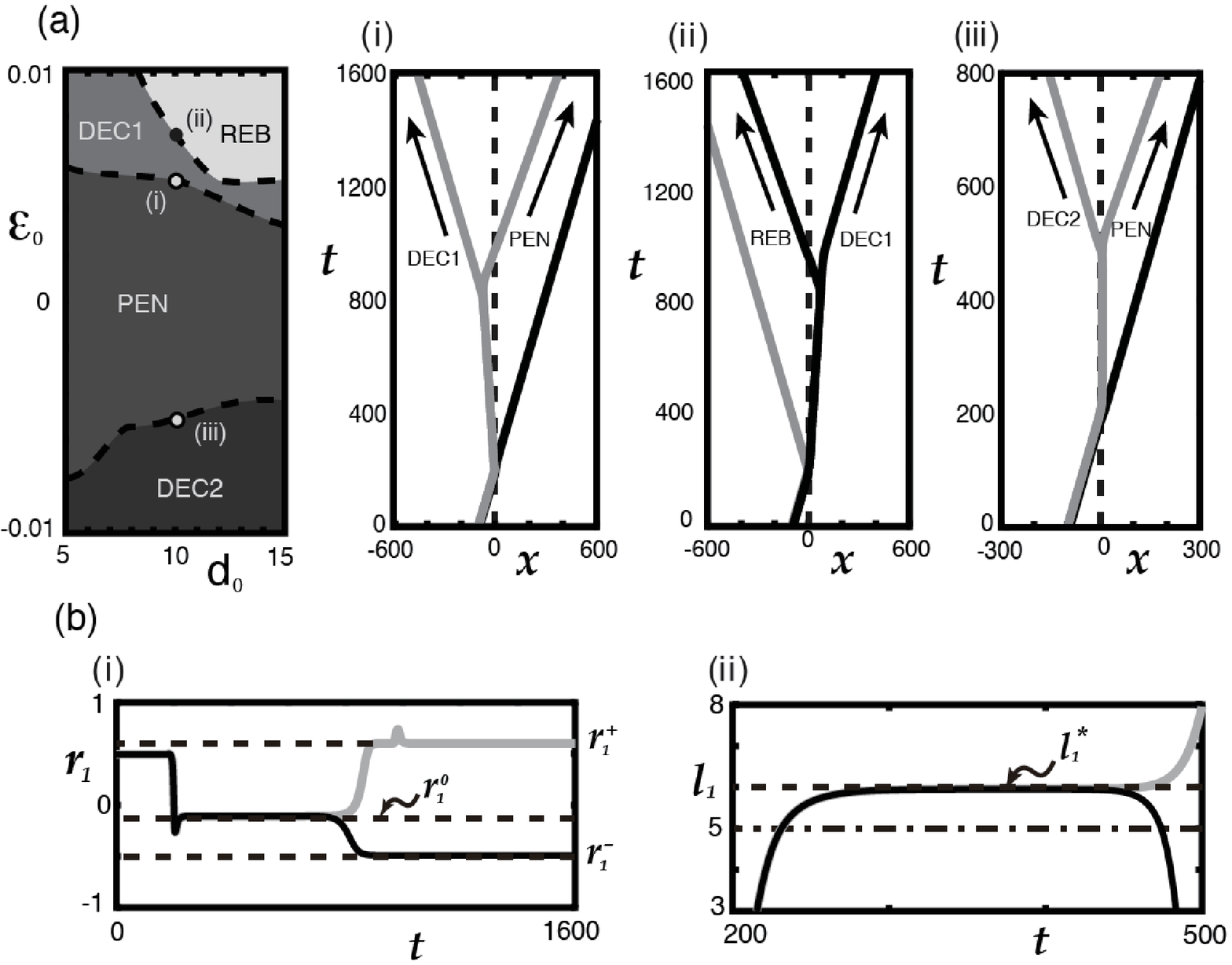}
 \caption{
   (a) Magnification of the phase diagram in Figure~\ref{phase_diagram_bump}(c) for the ODEs in (\ref{eq_ODE_renormalized_expanded1}) and the spatiotemporal plots of $l_2(t)$ (black) and $l_1(t)$ (gray) for parameters corresponding to the points labeled (i)--(iii) in the phase diagram. For (i), the parameters are close to the PEN-DEC1 phase boundary with $d_0=10$, and the two spatiotemporal plots for PEN and DEC1 are superposed. The same holds for cases (ii) and (iii). The other parameters were the same as in Figure~\ref{phase_diagram_bump}(c).
   (b) (i) Time series of $r_1(t)$ for the ODEs in (\ref{eq_ODE_renormalized_expanded1}) with the same parameters as in (a-i) for the PEN-DEC1 boundary. The black and gray lines correspond to the DEC1 and PEN behavior, respectively. The three horizontal broken lines denote $r^{+}_1=0.5985$, $r^{-}_1=-0.4942$, and $r^{0}_1=-0.1043$ computed from Eq.~(\ref{eq_front_velocity}).
   (ii) Time series of $l_1(t)$ for the ODEs in (\ref{eq_ODE_renormalized_expanded1}) with the same parameters as in (a-iii) for the PEN-DEC2 boundary. The black and gray lines correspond to the DEC2 and PEN behavior, respectively. The horizontal broken and dash-dot lines denote $l^{*}_1=5.9565$ computed from Eq.~(\ref{eq_front_position}) and the right edge of the bump heterogeneity at $d_0/2=5.0$, respectively. 
For (i) and (ii), the parameters were the same as in Figure~\ref{phase_diagram_bump}.
 }
   \label{phase_diagram_bump_scattor}
\end{figure}
\clearpage

%
\section{Summary and Discussion}
\subsection{Comparison between four dimensional ODEs in (\ref{eq_Ei_ODE2}) and (\ref{eq_ODE_renormalized_expanded2_0})}
In this study, the dynamics of pulse solutions have been examined both numerically and analytically for a bistable reaction-diffusion system (\ref{eq_PDE}). Applying the multiscale method to the hybrid system (\ref{eq_hybrid}), we formally derived four-dimensional ODEs. These were found to successfully reproduce the pulse dynamics observed in the original PDE system, both for the homogeneous and heterogeneous cases. In particular, they correctly preserved the order of the pitchfork and Hopf bifurcations of the SP solution, which the previously derived ODEs had failed to do within the framework of weak interaction theory (Fig.~\ref{fig_comparison}).
\begin{figure}[ht!]
   \centering
   \includegraphics[width=10cm]{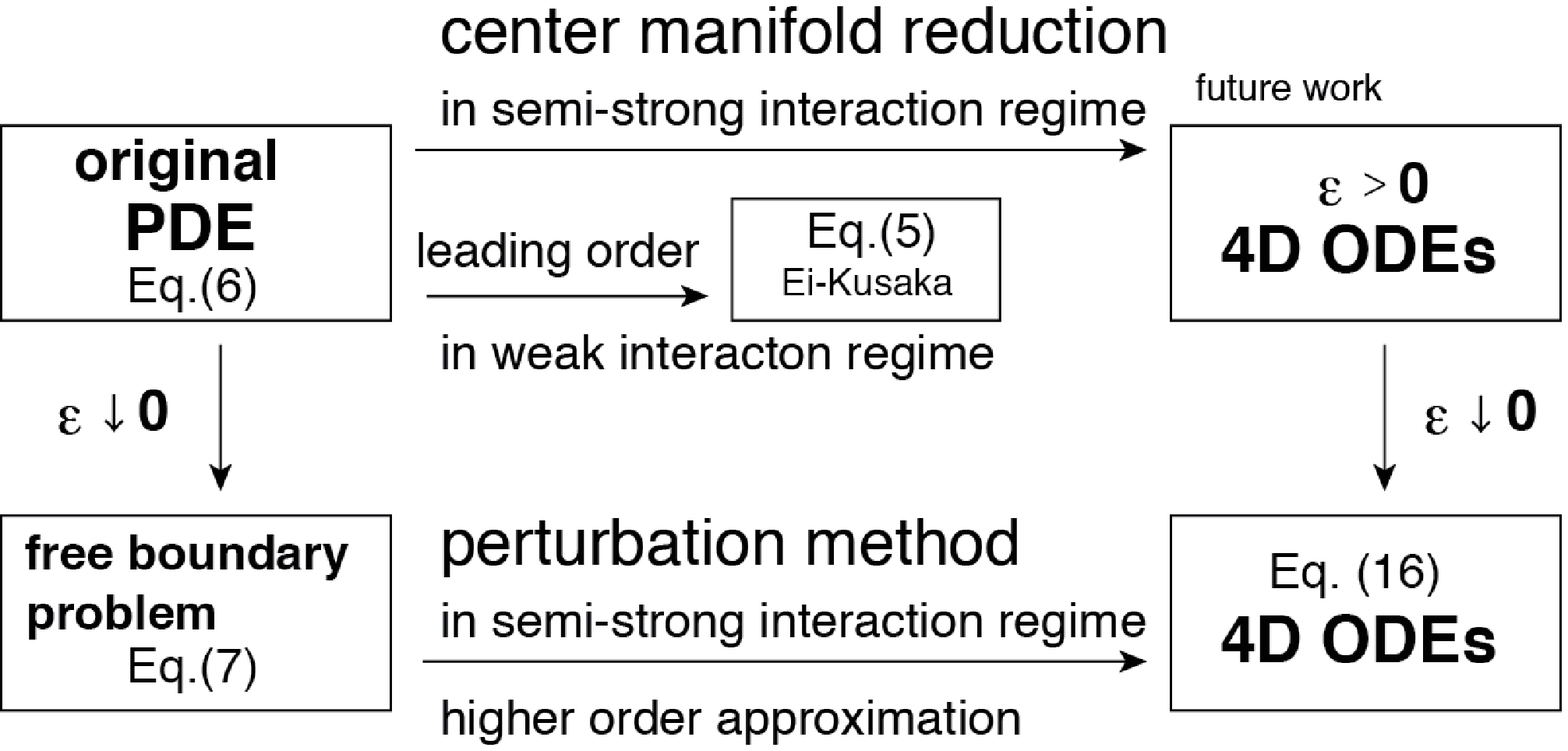}
 \caption{ %
   Comparison between the two kinds of reduction approach to the four-dimensional ODEs. In the present study, we made a detour via the limit system (\ref{eq_hybrid}) to the four-dimensional ODE system (\ref{eq_ODE_renormalized_expanded2_0}).}
   \label{fig_comparison}
\end{figure}

The pulse dynamics considered in this paper fall into the category of a semistrong front--front interaction \cite{Doelman_Eckhaus_Kaper_2000}\cite{Doelman_Eckhaus_Kaper_2001}\cite{Doelman_Kaper_2003}.
Away from the interfaces, the pulse solution comes very close to the equilibrium value for the $u$ component, but not for the $v$ component, which is largely deformed in between the two interfaces, as shown in Figure \ref{pulse_dynamics_pde_homo}(a) for a typical profile of the pulse solution. Hence, the front interaction is essentially determined by the slowly varying component of $v$.
In contrast, Ei and Kusaka \cite{Kusaka_2007} constructed the pulse solution by combining the front solutions in the framework of weak interaction theory. In this case, the fronts were assumed to be so far apart that they could only interact through their exponentially decaying tails, both for the $u$ and $v$ components.
Although the two ODEs (\ref{eq_Ei_ODE2}) and (\ref{eq_ODE_renormalized_expanded2_0}) were obtained by different approaches, they are not only similar in appearance, but also give excellent agreement in terms of their coefficients: constants such as $M_0$, $M_1$, and $M_2$ in Eq.~(\ref{eq_Ei_ODE2}) coincide with their counterparts in Eq.~(\ref{eq_ODE_renormalized_expanded2_0}) as $\epsilon \to 0$.
One of the major differences between (\ref{eq_Ei_ODE2}) and (\ref{eq_ODE_renormalized_expanded2_0}) lies in the interaction terms: Ei et al. only considered the lowest interaction terms $\pm \hspace{0.2mm}M_0 \hspace{0.2mm}e^{-h/\sqrt{D}}$, whereas we incorporated the higher-order ones of the form $(\hspace{0.5mm}G_0-G_1 \hspace{0.5mm}r_1\hspace{0.5mm}) \hspace{0.5mm}e^{-\frac{\;\;r_1+\phi(r_1)}{2 D}h}$ and $-(\hspace{0.5mm}G_0+G_1 \hspace{0.5mm}r_2\hspace{0.5mm}) \hspace{0.5mm}e^{-\frac{-r_2+\phi(r_2)}{2 D}h}$.
These observations suggest that the reduction from the original PDEs attempted by Ei et al. using center manifold reduction may work well if the higher-order interaction terms are appropriately taken into account, which is left as future work.
One may also be inclined to apply a renormalization group method for the reduction of pulse behavior, as in \cite{Heijster_Doelman_Kaper_Promislow_2010}, in which the authors not only derived ODEs for the motion of $N$-fronts by the geometric singular perturbation technique, but also ensured their validity based on the idea of the renormalization group method. As the resulting ODEs only describe monotonic attracting/repelling motion of the fronts, it would be a challenge to extend their rigorous approach to admit time-periodic and more complicated motions, as shown in Figures \ref{pulse_dynamics_ode_homo} and \ref{fig_square_well}.

\subsection{Sliding motion of standing breather in heterogeneous media}
As an application of the reduced ODE system (\ref{eq_ODE_renormalized_expanded1}), we also considered the propagation manner of the traveling pulse that encountered the bump-type heterogeneity, and revealed the influence of parameter variations on the pulse behavior. In particular, we focused on the behavioral changes that occurred near the phase boundaries to identify the unstable solutions, called scattors, which played a central role in determining the pulse behavior. The related unstable traveling and stationary solutions were explicitly given based on the reduced system.

In this study, the bump height $\epsilon_0$ and the bump width $d_0$ were varied as bifurcation parameters. Note that the bump width becomes infinitely large as $d_0 \to \infty$. In this limit, the bump-type heterogeneity can be regarded as a jump, which was studied in our previous paper \cite{Nishi_Nishiura_Teramoto_2013}.
As in the bump case, the PEN and REB behavior were also observed for the jump-type heterogeneity, for which the traveling pulse coming from the left eventually converges to another stable traveling pulse that exists for the homogeneous system either in the right infinity $x \to \infty$ (PEN) or the left infinity $x \to -\infty$ (REB).
However, the DEC1 regime appeared between the PEN and REB regimes (Figs.~\ref{phase_diagram_bump}(b),(c)); this was not observed for the jump-type case. The DEC1 regime occupies a relatively large parameter space when the bump width $d_0$ is comparable to the pulse size, and becomes narrower as $d_0$ increases. Eventually, the DEC1 regime seems to disappear for large values of $d_0$, which indicates that the PEN behavior directly changes to REB behavior as $\epsilon_0$ increases.
We numerically confirmed this for the reduced ODE system (Figs.~\ref{fig_wide_bump}(a)(b)). After the pulse entered the bump region, it started to oscillate, and its amplitude gradually grew to that of a stable oscillatory solution for the homogeneous system $\delta(x)=\delta_0+\epsilon_0$.
When the bump width was not sufficiently large, the pulse could not sustain the oscillation inside the bump region. After crossing the bump edges, the pulse interfaces decomposed into two counter-propagating front solutions.
For sufficiently large values of $d_0$, however, the oscillatory pulse did not immediately decompose into front solutions, but rather exhibited slow sliding motion that lasted exponentially long as the bump height $\epsilon_0$ remained close to the PEN-REB boundary at $\epsilon^{\scalebox{0.5}{PEN-REB}}_0$ (Fig.~\ref{fig_wide_bump}(c)). The mechanism for this sliding motion was analyzed for the jump case in \cite{Nishi_Nishiura_Teramoto_2013}. As a result of the sliding motion, the oscillatory pulse slowly approached one of the bump edges, and its interface eventually crossed the bump edge, leading to either the PEN or REB behavior.

These numerical observations strongly suggest that there exists an unstable oscillatory pulse solution located at the center of the bump, playing the role of a scattor for the PEN and REB behavior.
As we showed in Section 3, the reduced ODEs (\ref{eq_ODE_renormalized_expanded1}) can also be used to study unstable oscillatory pulse solutions. As these oscillatory pulse solutions often appear following the Hopf bifurcation of a stationary pulse solution, we first seek the corresponding stationary solution induced by the heterogeneity.\vspace{0.2cm}

\noindent
{\bf Lemma 4.}
{\it
For $\epsilon_0 < (G_0 - \delta_0)/(1-e^{-\frac{d_0}{2\sqrt{D}}})$, there exists a unique stationary solution to Eq.~(\ref{eq_ODE_renormalized_expanded1}) with (\ref{eq_Delta0}) for each value of $\epsilon_0$.
}
\vspace{0.2cm}\\
{\bf Proof.}\hspace{0.1cm}
{\it
Setting $\dot{l}_2=\dot{l}_1=0, \hspace{0.8mm} \dot{r}_2=\dot{r}_1=0$ in Eq.~(\ref{eq_ODE_renormalized_expanded1}) with (\ref{eq_Delta0}) yields
\begin{eqnarray}
\displaystyle{ 
\quad G_0 \hspace{0.2mm}e^{-\frac{h}{\sqrt{D}}}=\Delta_0(h/2)\hspace{0.5mm}, 
}
 \label{eq_HISP}
\end{eqnarray}
where $h=l_2-l_1>0$. By introducing a variable $z:=e^{-\frac{h}{2\sqrt{D}}}$ and a constant $z_0:=e^{-\frac{d_0}{2\sqrt{D}}}$, (\ref{eq_HISP}) with (\ref{eq_Delta0}) can be transformed as 
\begin{eqnarray}
 \epsilon_0=f(z)\hspace{0.5mm}, \label{eq_HISP2}
\end{eqnarray}
where
\begin{eqnarray}
\hspace{0.15cm}f(z)=\left \{
\begin{array}{l}
\displaystyle{ 
  \quad \frac{G_0 \hspace{0.2mm}z^{2} -\delta_0}{1-z_0 \left(\hspace{0.5mm} z+z^{-1} \hspace{0.5mm}\right)/2}\hspace{0.5mm}, \hspace{0.57cm}  (\hspace{0.5mm}z_0 \le z \hspace{0.5mm}) \vspace{0.5mm} } \vspace{0.2cm}\\
 \displaystyle{  \quad \quad \frac{G_0 \hspace{0.2mm}z^{2} -\delta_0}{z \left(\hspace{0.5mm} z^{-1}_0-z_0 \hspace{0.5mm}\right)/2}\hspace{0.5mm}. \hspace{0.8cm} (\hspace{0.5mm}0 < z \le z_0\hspace{0.5mm})}
\end{array}
\right. \label{eq_HISP3}
\end{eqnarray}
Note that $0< h < \infty$ is equivalent to $0 < z <1$. The function $f(z)$ in (\ref{eq_HISP3}) is continuous and monotonically increasing for $z > 0$, which, together with the monotonicity of $z=e^{-\frac{h}{2\sqrt{D}}}$, results in the unique existence of a stationary solution to Eq.~(\ref{eq_ODE_renormalized_expanded1}) with (\ref{eq_Delta0}) with a finite value of $h$ for $\epsilon_0 < f(1)= (G_0 - \delta_0)/(1-e^{-\frac{d_0}{2\sqrt{D}}})$.
} \hspace{0.0cm} \scalebox{1.0}{$\Box$} \vspace{0.2cm}
This lemma indicates the unique existence of a stationary pulse solution of width $h=-2\sqrt{D} \hspace{0.2mm}\log f^{-1}(\epsilon_0)$, centered at $x=0$, for the range of $\epsilon_0$ we now consider.
However, the unstable oscillatory pulse solution branch that originated from the stationary pulse solution does not reach the parameter range under consideration.
In fact, AUTO \cite{Doedel_2008} numerically revealed the stability of the stationary pulse solution, as shown in Figure \ref{fig_unstable_SB}.
Fixing $\epsilon_0$ and varying $\tau$ as a bifurcation parameter, the stationary pulse solution underwent a Hopf bifurcation at $\tau_{H}$, from which an unstable oscillatory pulse solution branch emerged (Fig.~\ref{fig_unstable_SB}(a)). The amplitude of the oscillation increased monotonically as $\tau$ decreased, and seemed to diverge at $\tau_{\infty}$, indicating that the oscillatory pulse solution exists for $\tau_{\infty}<\tau<\tau_{H}$.
By varying $\epsilon_0$ as well, we obtained the $\tau$\hspace{0.5mm}-\hspace{0.5mm}$\epsilon_0$ diagram for the existence of the oscillatory pulse solution when $d_0=40$ (Fig.~\ref{fig_unstable_SB}(b)). In particular, the diagram indicates that the unstable oscillatory pulse solution exists for $\epsilon_0 > \tilde{\epsilon}_0 \approx 0.0038$ when $\tau=0.1700$, the value used in the numerical simulations illustrated in Figures~\ref{phase_diagram_bump}(c) and \ref{fig_wide_bump}.
In contrast, the numerical results in Figure~\ref{phase_diagram_bump}(c) show that, for a bump width of $d_0=40$, the PEN-REB transition occurs around $\epsilon_0= 0.002$, which is outside the aforementioned range of $\epsilon_0$ for the existence of the unstable oscillatory solution.

This implies that the unstable oscillatory pulse solution described above is not responsible for the PEN-REB transition. It remains to be elucidated whether there is another unstable oscillatory solution that plays the role of the scattor for the PEN-REB transition.

\subsection{Designing trapped motion by square-well-type heterogeneity}
Finally, we remark that the reduction method presented here is not limited to the dynamics of the two interacting fronts along the whole line, but can readily be extended to the pulse dynamics in finite domains \cite{Suzuki_Ohta_Mimura_Sakaguchi_1995}\cite{Tsyganov_Ivanitsky_Zemskov_2014}, as well as to dynamics that involve more than two interfaces \cite{Ohta_Ito_Tetsuka_1990}\cite{Ohta_Nakazawa_1992}.
Another application of our study is the design of spatial heterogeneity based on the reduced system. In the present paper, we have dealt with the bump-type heterogeneity, partly because this is one of the most fundamental and ubiquitous heterogeneities seen in nature. In principle, however, the function $\delta(x)$ in Eq.~(\ref{eq_PDE}) can take any type of spatial heterogeneity, and may be designed to give the desired dynamics of moving localized patterns.

For instance, the pinning and trapping of moving localized patterns, for which pulses or spots become trapped in a finite region of space, are one of the most dramatic effects of heterogeneity. They play a central role in controlling the motion of localized patterns, together with the reversal and change of propagation direction \cite{Lober_Engel_2014}\cite{Mikhailov_Showalter_2006}.
Such trapped motion of traveling pulses and fronts has been observed for reaction-diffusion systems with jump- or bump-type heterogeneities \cite{Heijster_Doelman_Kaper_Nishiura_Ueda_2011}\cite{Heijster_Chen_Nishiura_Teramoto_2017}\cite{Ikeda_Ei_2010}\cite{Yuan_Teramoto_Nishiura_2007}, where, unlike in our case, the heterogeneity was introduced to the activator ($u$ component) rather than the inhibitor ($v$ component). In fact, trapped motion was not observed in our system, where the bump-type heterogeneity was introduced to the $v$ component.

Regardless, we can realize trapped motion by manipulating the form of $\delta(x)$ in Eq.~(\ref{eq_PDE}), based on the information obtained for the bump-type heterogeneity.
To this end, we note that the traveling pulse exhibited the PEN and REB behavior depending on the bump height (Fig.~\ref{phase_diagram_bump}), which readily leads to the idea that the traveling pulse may by trapped in between the two bumps by appropriately adjusting the bump height.
Thus, we applied the square-well-type heterogeneity for $\delta(x)$ in (\ref{eq_delta}):
\begin{eqnarray}
\delta(x) =\left \{
\begin{array}{l}
\displaystyle{ 
  \quad \delta_0 \hspace{0.5cm}(\hspace{0.5mm}x<-d_0/2-d_1\hspace{0.5mm}), } \vspace{0.2cm}\\
  \displaystyle{ 
  \quad \delta_0+\epsilon_1 \hspace{0.5cm}(\hspace{0.5mm}-d_0/2-d_1 \le x<-d_0/2\hspace{0.5mm}), } \vspace{0.2cm}\\
  \displaystyle{ 
  \quad \delta_0+\epsilon_2 \hspace{0.5cm}(\hspace{0.5mm}-d_0/2 \le x< d_0/2\hspace{0.5mm}), } \vspace{0.2cm}\\
    \displaystyle{ 
  \quad \delta_0+\epsilon_1 \hspace{0.5cm}(\hspace{0.5mm}d_0/2 \le x< d_0/2+d_1\hspace{0.5mm}), } \vspace{0.2cm}\\
  \displaystyle{ 
  \quad \delta_0 \hspace{0.5cm}(d_0/2+d_1 \le x),
  } \vspace{0.2cm}
\end{array}
\right. \label{eq_square_well}
\end{eqnarray}
and observed trapped oscillatory motions of the traveling pulse for the ODE system (Fig.~\ref{fig_square_well}) by changing the values of $\epsilon_1$ and $\epsilon_2$. These motions were also found for the original PDE system. Note that the swaying motion in Figure~\ref{fig_square_well}(b-iv) does not occur when $\epsilon_2=0$, where the two identical bumps are simply placed side by side.
Thus, the pulse dynamics for the original PDE system may be explored through the reduced ODEs, which facilitates numerical and analytical study of more complex behavior for other types of heterogeneity, as well as the design of external perturbations for controlling the traveling pulse.
\begin{figure}
   \centering
 \includegraphics[width=8cm]{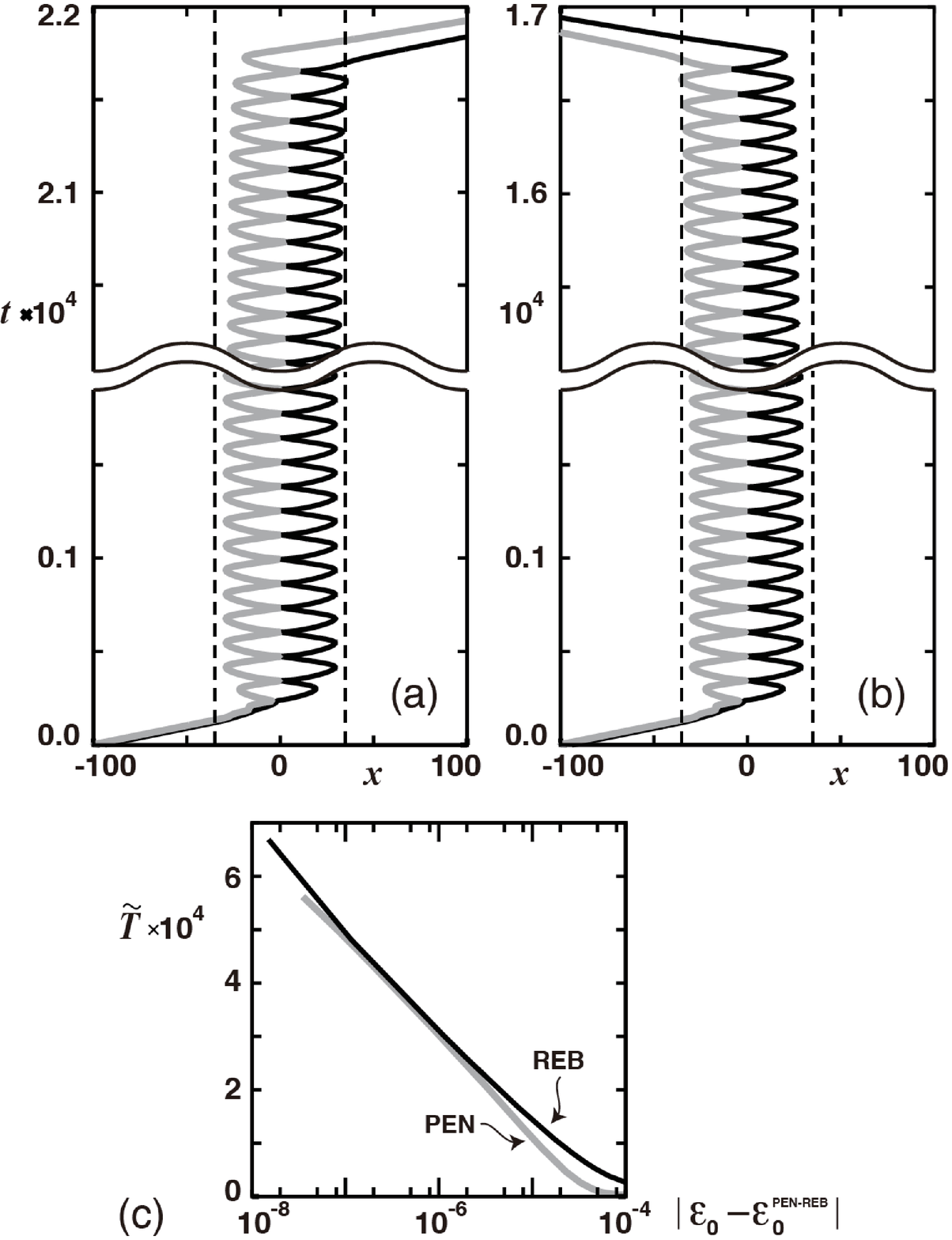}
 \caption{
   PEN-REB transition for the wide bump case $d_0=70$ observed for the ODEs in (\ref{eq_ODE_renormalized_expanded1}). 
   The PEN-REB transition occurred around $\epsilon^{\rm PEN-REB}_0 \approx 0.00230289$.
   In (a) and (b), typical spatiotemporal plots for PEN (\hspace{0.5mm}$\epsilon_0 = 0.00230$\hspace{0.5mm})
   and REB (\hspace{0.5mm}$\epsilon_0 = 0.00231$\hspace{0.5mm}) behavior are shown.
   The other parameters were the same as in Figure~\ref{phase_diagram_bump}.
   The vertical dashed lines indicate the right and left edges of a bump domain.
   The center of the oscillatory pulse is sliding slowly to the left or to the right. 
   (c) The $\epsilon_0$-dependence of the time interval $\tilde{T}$ during which the pulse stayed inside the bump domain.
   Note the log scale of the horizontal axis.
   The oscillation lasted longer as $\epsilon_0$ became closer to the PEN-REB boundary at $\epsilon_0^{\rm PEN-REB}$.
 }
   \label{fig_wide_bump}
\end{figure}
\begin{figure}
   \centering
 \includegraphics[width=12cm]{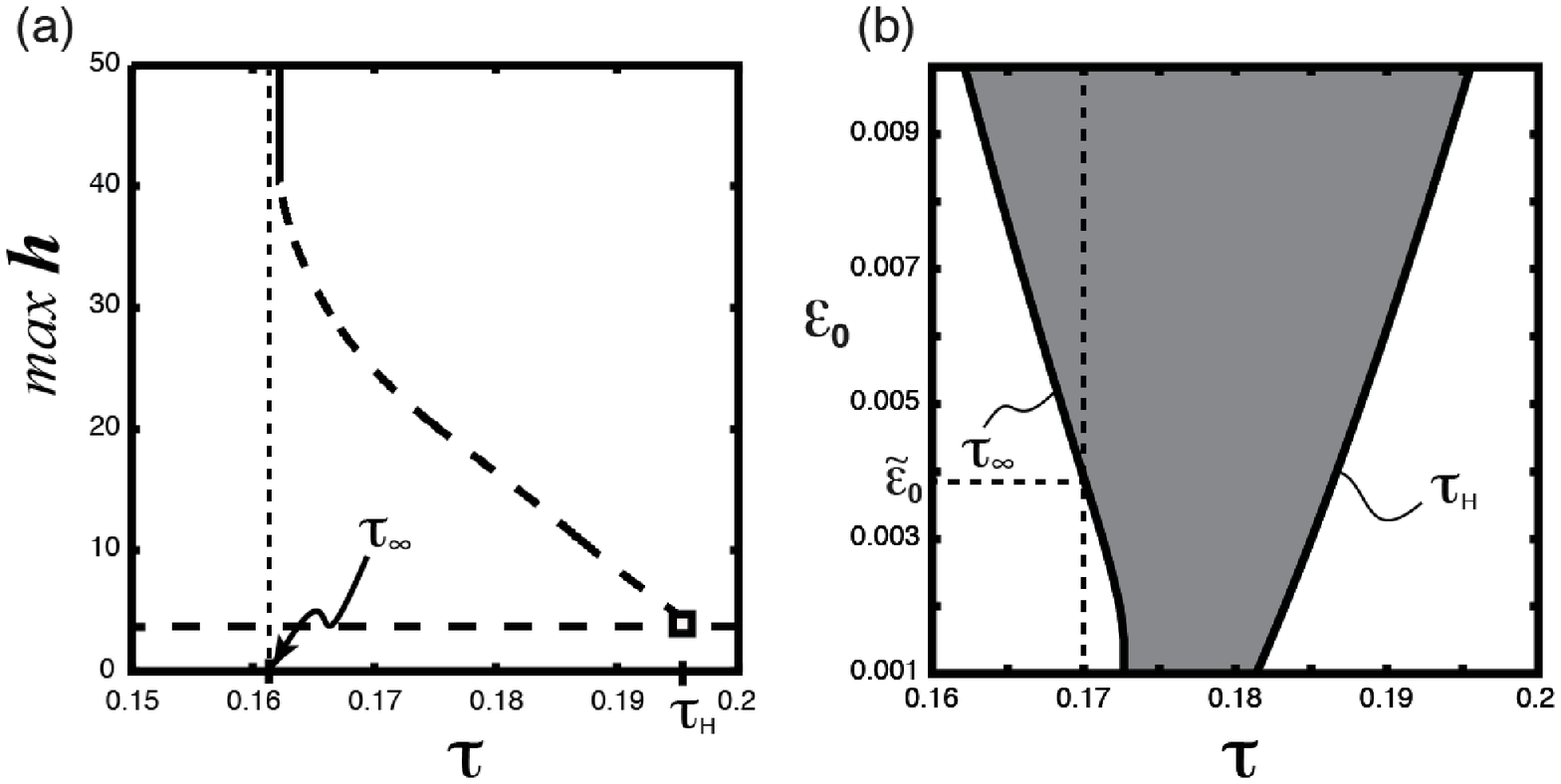}
 \caption{
   Existence of the heterogeneity-induced oscillatory pulse solution to the ODEs in (\ref{eq_ODE_renormalized_expanded1}), as found numerically by means of AUTO \cite{Doedel_2008}. The bump width was fixed at $d_0=40$. (a) Bifurcation branch of the stationary pulse and the oscillatory pulse solutions for $\epsilon_{0}=0.01$. Broken lines denote unstable solutions. When $\tau$ was varied as a bifurcation parameter, the stationary pulse solution underwent a Hopf bifurcation at $\tau_{H}=0.1955$, from which an unstable oscillatory pulse solution emerged. The maximal amplitude of the oscillatory pulse solution increased monotonically and seemed to diverge at $\tau_{\infty}=0.16216$. (b) By varying $\epsilon_{0}$ and plotting the values of $\tau_{H}$ and $\tau_{\infty}$ in (a), we obtained a $\tau$\hspace{0.5mm}-\hspace{0.5mm}$\epsilon_{0}$ diagram, for which the unstable oscillatory pulse solution existed in the shaded region. The value of $\epsilon_{0}$ at which the $\tau_{\infty}$ curve intersected with $\tau=0.170$ was indicated by $\tilde{\epsilon}_{0} \approx 0.0038$. For (a) and (b), the other parameters were the same as in Figure~\ref{phase_diagram_bump}.}
   \label{fig_unstable_SB}
\end{figure}
\begin{figure}
   \centering
 \includegraphics[width=12cm]{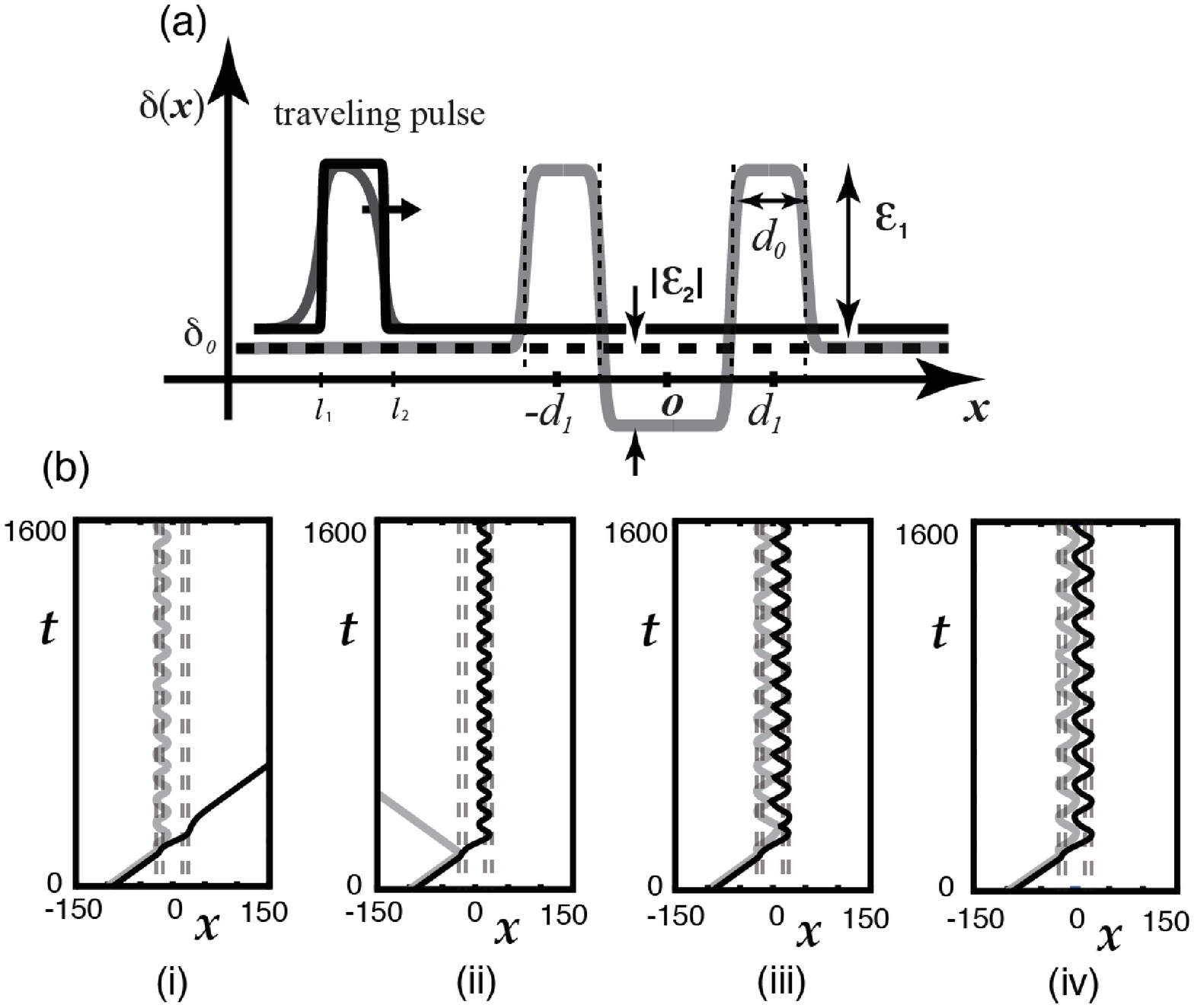}
 \caption{ (a) Schematic figure for the traveling pulse that encounters the square-well-type heterogeneity.
   (b) Trapped motions of the pulse interfaces observed for the ODE system in (\ref{eq_ODE_renormalized_expanded1}) with the square-well-type heterogeneity (\ref{eq_square_well}). The vertical broken lines denote the edges of the square-well at $x=\pm \hspace{0.5mm}d_0/2$, $d_0/2+d_1$, and $-d_0/2-d_1$, where $d_0$ and $d_1$ were fixed to $d_0=10$ and $d_1=30$. The values of $\epsilon_1$ and $\epsilon_2$ were set to (i) $\epsilon_1=0.0048$, $\epsilon_2=-0.0080$, (ii) $\epsilon_1=0.0070$, $\epsilon_2=-0.0070$, (iii) $\epsilon_1=0.0048$, $\epsilon_2=-0.0040$, and (iv) $\epsilon_1=0.0048$, $\epsilon_2=-0.0048$. The other parameters were the same as in Figure~\ref{phase_diagram_bump}.
 }
   \label{fig_square_well}
\end{figure}
\clearpage

%

%
\appendix
\section{Derivation of reduced ODE system (\ref{eq_ODE_renormalized_full})}
We introduce a slow time $T = \mu \hspace{0.5mm}t$ and consider a perturbation series for $v$ in Eq.~(\ref{eq_hybrid}) of the form  
\begin{eqnarray}
  \begin{array}{ll}
    \hspace{0.8cm}v(x,t,T) = v_0(x,t,T) + \mu \hspace{0.5mm}v_1 (x,t,T)+ \mu^2 \hspace{0.5mm}v_2 (x,t,T) \\
    \hspace{6.8cm} \displaystyle + \mu^3 \hspace{0.5mm}v_3 (x,t,T)+ \cdots ,
  \label{eq_expansion1}
  \end{array}
\end{eqnarray}
where $0<\hspace{0.5mm}\mu \ll \hspace{0.5mm}1$ is some infinitesimally small parameter. Then, the total time derivative can be rewritten as $\hspace{1.0mm}\partial/\partial t \hspace{1.0mm}\rightarrow \hspace{1.0mm}\partial/\partial \hspace{0.5mm}t+\mu \hspace{0.5mm}\partial/\partial \hspace{0.5mm}T\hspace{0.5mm}$.
Substituting this into the third equation in (\ref{eq_hybrid}), we obtain 
\begin{eqnarray}
\begin{array}{l}
  \displaystyle \hspace{1.5cm}\left(\hspace{0.5mm} \frac{\partial}{\partial \hspace{0.5mm}t} + \mu \hspace{0.5mm}\frac{\partial}{\partial \hspace{0.5mm}T} \hspace{0.5mm}\right) 
\left(\hspace{0.5mm}v_0 + \mu \hspace{0.5mm}v_1+ \mu^2 \hspace{0.5mm}v_2 + \cdots \right) \vspace{0.25cm}\\
 = \displaystyle \hspace{0.2cm} D \hspace{0.5mm}\frac{\partial^2}{\partial x^2} \left(\hspace{0.5mm}v_0 + \mu \hspace{0.5mm}v_1+ \mu^2 \hspace{0.5mm}v_2 + \cdots \right) + u -\left(\hspace{0.5mm}v_0 + \mu \hspace{0.5mm}v_1+ \mu^2 \hspace{0.5mm}v_2 + \cdots \right) + \delta(x).
\end{array} \nonumber
\end{eqnarray}
\begin{equation}
\hspace{5.0cm} \label{eq_expansion2}
\end{equation}
Collecting terms with equal powers of $\mu$ yields equations for $v_i \hspace{1.5mm}(\hspace{0.5mm}i=0,1,2,\cdots\hspace{0.5mm})$:
\begin{eqnarray}
\begin{array}{l}
\displaystyle{ 
\mathcal{O}(\; \mu^0 \;): \quad \frac{\partial v_0}{\partial t}=\left( D \frac{\partial^2}{\partial x^2} -1 \right) v_0(x,T)+ u(\hspace{0.5mm}x;l_2,l_1\hspace{0.5mm})+\delta(x), } \vspace{0.2cm} \\
\displaystyle{
\mathcal{O}(\; \mu^1 \;): \quad \frac{\partial v_1}{\partial t}=\left( D \frac{\partial^2}{\partial x^2} -1 \right) v_1(x,T)-\frac{\partial v_0}{\partial T}, } \vspace{0.2cm} \\
\displaystyle{
\mathcal{O}(\; \mu^2 \;): \quad \frac{\partial v_1}{\partial t}=\left( D \frac{\partial^2}{\partial x^2} -1 \right) v_2(x,T)-\frac{\partial v_1}{\partial T}, } \vspace{0.2cm} \\
\hspace{5.0cm} \vdots \vspace{0.2cm}\\
\displaystyle{ 
\mathcal{O}(\; \mu^n \;): \quad \frac{\partial v_n}{\partial t}=\left( D \frac{\partial^2}{\partial x^2} -1 \right) v_n(x,T)-\frac{\partial v_{n-1}}{\partial T}. }
\end{array} \label{eq_order}
\end{eqnarray}
Here, we focus on the asymptotic pulse behavior governed by the slow timescale $T$, and assume that $v$, $l_2$, and $l_1$ are independent of $t$ as $t \to \infty$. Hence, $\partial \hspace{0.5mm}v_i \hspace{0.5mm}/ \hspace{0.5mm}\partial \hspace{0.5mm}t=0 \hspace{1.5mm}(\hspace{0.5mm}i=0,1,2,\cdots\hspace{0.5mm})$ in Eq.~(\ref{eq_order}). Each equation in (\ref{eq_order}) then reduces to an ODE with respect to $x$, and is solved iteratively as follows.

For the zeroth order in Eq.~(\ref{eq_order}), the equation reads
\begin{eqnarray}
\begin{array}{rcl}
\displaystyle \left( D \frac{\partial^2}{\partial x^2} - 1 \right) v_0(x,T_1)
+ u(x \hspace{0.5mm};l_2, l_1) +\delta(x)= 0,
\end{array} \label{eq_0th_equation}
\end{eqnarray}
where $u(x\hspace{0.5mm};l_2,l_1)= F(\hspace{0.2mm}x-l_1\hspace{0.2mm}) - F(\hspace{0.2mm}x - l_2\hspace{0.2mm}) - 1/\hspace{0.2mm}2$, as shown in (\ref{eq_function_u}). This is readily solved as
\begin{eqnarray}
\begin{array}{c}
  \displaystyle v_0(x,T_1) = - \frac{1}{2} + \Delta_0(x) + \overline{v}_0(\hspace{0.5mm}x-l_2(T)\hspace{0.5mm}) - \overline{v}_0(\hspace{0.5mm}x-l_1(T)\hspace{0.5mm}),
\end{array} \label{eq_0th_solution}
\end{eqnarray}
where $\Delta_0(x)$ and $\overline{v}_0(x)$ satisfy
\begin{eqnarray}
\begin{array}{ll}
  \displaystyle L \hspace{0.5mm}\overline{v}_0(x)+ F(x) =0\hspace{0.5mm}, \hspace{1.0cm}L \hspace{0.5mm}\Delta_0(x)+ \delta(x) =0\hspace{0.5mm},
\end{array} \label{eq_0th_function}
\end{eqnarray}
with $L := D \hspace{0.5mm}d^2/d x^2 - 1$.
Substituting (\ref{eq_0th_solution}) into the second-order equation in (\ref{eq_order}) yields the equation for $v_1$:
\begin{eqnarray}
\begin{array}{ll}
\displaystyle \left( D \frac{\partial^2}{\partial x^2} - 1 \right) \hspace{0.5mm}v_1(x,T) 
\displaystyle =\frac{d \hspace{0.5mm}l_2}{d \hspace{0.5mm}T}\hspace{1.0mm}\frac{\partial \overline{v}_0}{\partial x}(x-l_2)- \frac{d \hspace{0.5mm}l_1}{d \hspace{0.5mm}T}\hspace{1.0mm} \frac{\partial \overline{v}_0}{\partial x}(x-l_1)  \hspace{0.5mm},
\end{array} \label{eq_1st_equation}
\end{eqnarray}
which is solved by
\begin{eqnarray}
\begin{array}{ll}
\displaystyle  v_1(x,T) = \frac{d \hspace{0.5mm}l_2}{d \hspace{0.5mm}T}\hspace{1.0mm}\overline{v}_1(x-l_2(T)) - \frac{d \hspace{0.5mm}l_1}{d \hspace{0.5mm}T}\hspace{1.0mm}\overline{v}_1(x-l_1(T)) \hspace{0.5mm},
\end{array} \label{eq_1st_solution}
\end{eqnarray}
where $\overline{v}_1$ satisfies
\begin{eqnarray}
\displaystyle{ L \hspace{0.5mm}\overline{v}_1(x) = - \frac{d \overline{v}_0}{dx}}\hspace{0.5mm}. \label{eq_1st_function}
\end{eqnarray}
Similarly, the second-order equation in Eq.~(\ref{eq_order}) can be rewritten as 
\begin{eqnarray}
\begin{array}{ll}
\displaystyle \left( D \frac{\partial^2}{\partial x^2} - 1 \right) \hspace{0.5mm}v_2(x,T) =\frac{d^2 \hspace{0.5mm}l_2}{d \hspace{0.5mm}T^2}\hspace{1.0mm}\frac{\partial \overline{v}_1}{\partial x}(x-l_2)-\left( \hspace{0.5mm} \frac{d \hspace{0.5mm}l_2}{d \hspace{0.5mm}T}  \hspace{0.5mm}\right)^2 \hspace{1.0mm}\frac{d \overline{v}_1}{d x}(x-l_2) \vspace{0.2cm}\\
\displaystyle \hspace{4.0cm}-\frac{d^2 \hspace{0.5mm}l_1}{d \hspace{0.5mm}T^2}\hspace{1.0mm}\frac{\partial \overline{v}_1}{\partial x}(x-l_1)+\left( \hspace{0.5mm} \frac{d \hspace{0.5mm}l_1}{d \hspace{0.5mm}T}  \hspace{0.5mm}\right)^2 \hspace{1.0mm}\frac{d \overline{v}_1}{d x}(x-l_1),
\end{array} \label{eq_2nd_equation}
\end{eqnarray}
which is solved by
\begin{eqnarray}
\begin{array}{ll}
  \displaystyle
 v_2(x,T)=\frac{d^2 \hspace{0.5mm}l_2}{d \hspace{0.5mm}T^2}\hspace{1.0mm} \overline{v}_{20}(x-l_2(T))-\left( \hspace{0.5mm} \frac{d \hspace{0.5mm}l_2}{d \hspace{0.5mm}T}  \hspace{0.5mm}\right)^2 \hspace{1.0mm}\overline{v}_{11}(x-l_2(T)) \vspace{0.2cm}\\
\displaystyle \hspace{2.5cm}-\frac{d^2 \hspace{0.5mm}l_1}{d \hspace{0.5mm}T^2}\hspace{1.0mm} \overline{v}_{20}(x-l_1(T))+\left( \hspace{0.5mm} \frac{d \hspace{0.5mm}l_1}{d \hspace{0.5mm}T}  \hspace{0.5mm}\right)^2 \hspace{1.0mm} \overline{v}_{11}(x-l_1),
\end{array} \label{eq_2nd_solution}
\end{eqnarray}
where $\overline{v}_{20}$ and $\overline{v}_{11}$ satisfy
\begin{eqnarray}
  \displaystyle L \hspace{0.5mm}\overline{v}_{20}(x) = \overline{v}_1 \hspace{0.5mm}, \hspace{0.5cm}L \hspace{0.5mm}\overline{v}_{11}(x) = \frac{d \overline{v}_1}{dx}\hspace{0.5mm}. \label{eq_2nd_function}
\end{eqnarray}
Thus, calculating up to $\mathcal{O}(\hspace{0.5mm} \mu^4 \hspace{0.5mm})$, we obtain
%
\begin{eqnarray}
\begin{array}{ll}
  \mathcal{O}(\mu^{0}):
  & \displaystyle{ v_0(x,T) = -\frac{1\hspace{0.2mm}}{\hspace{0.2mm}2}+\Delta_0(x)+\overline{v}_0(x-l_2)
      - [\hspace{0.5mm} \mbox{terms of } \hspace{1.0mm}l_1 \hspace{0.5mm}]\hspace{0.5mm}, } \vspace{0.3cm}\\
  & \displaystyle{ \hspace{0.2cm} L \hspace{0.5mm}\overline{v}_0 + F(x) = 0 \hspace{0.5mm}, \hspace{0.5cm} L \hspace{0.5mm}\Delta_0(x) + \delta(x) = 0 \hspace{0.5mm},} \vspace{0.4cm}\\
  \mathcal{O}(\mu^1): 
  & \displaystyle{ v_1(x,T) = \left(\hspace{0.2mm} \frac{d \hspace{0.2mm}l_2}{d \hspace{0.2mm}T}\hspace{0.2mm} \right) \hspace{1.0mm}\overline{v}_1(x-l_2) - [\hspace{0.5mm} \mbox{terms of } \hspace{1.0mm}l_1 \hspace{0.5mm}]\hspace{0.5mm}, } \vspace{0.2cm}\\
  & \displaystyle{ \hspace{0.2cm} L \hspace{0.5mm}\overline{v}_1 =- \frac{d \hspace{0.2mm}\overline{v}_0}{d \hspace{0.2mm}x} \hspace{0.5mm}, } \vspace{0.4cm}\\
  \mathcal{O}(\mu^2):
 & \displaystyle{ v_2(x,T) = \left(\frac{d^2 \hspace{0.2mm}l_2}{d \hspace{0.2mm}T^2} \right) \hspace{1.0mm}\overline{v}_{20}(x-l_2) - \left( \frac{d \hspace{0.2mm}l_2}{d \hspace{0.2mm}T}\right)^2 \hspace{1.0mm} \overline{v}_{11}(x-l_2) } \vspace{0.2cm}\\
  & \displaystyle{ \hspace{6.5cm} - [\hspace{0.5mm} \mbox{terms of } \hspace{1.0mm}l_1 \hspace{0.5mm}]\hspace{0.5mm}, } \vspace{0.0cm}\\
  & \displaystyle{ \hspace{0.2cm} L \hspace{0.5mm}\overline{v}_{20} = \overline{v}_1\hspace{0.5mm}, \hspace{0.5cm} L \hspace{0.5mm}\overline{v}_{11} =\frac{d \hspace{0.2mm}\overline{v}_1}{d \hspace{0.2mm}x} \hspace{0.5mm}, } \vspace{0.4cm}\\
  \mathcal{O}(\mu^3):
  & \displaystyle{ v_3(x,T) = \left( \frac{d^3 \hspace{0.2mm}l_2}{d \hspace{0.2mm}T^3}\right) \hspace{1.0mm}\overline{v}_{300}(x-l_2) - \left(\frac{d^2 \hspace{0.2mm}l_2}{d \hspace{0.2mm}T^2}\right)\hspace{0.2mm}\left(\frac{d \hspace{0.2mm}l_2}{d \hspace{0.2mm}T}\right) \hspace{1.0mm} \overline{v}_{210}(x-l_2) } \vspace{0.2cm}\\
  & \displaystyle{ \hspace{1.8cm}+ \left(\frac{d \hspace{0.2mm}l_2}{d \hspace{0.2mm}T}\right)^3 \hspace{1.0mm} \overline{v}_{111}(x-l_2) - [\hspace{0.5mm} \mbox{terms of } \hspace{1.0mm}l_1 \hspace{0.5mm}]\hspace{0.5mm}, } \vspace{0.3cm}\\
  & \displaystyle{ \hspace{0.2cm} L \hspace{0.5mm}\overline{v}_{300} = \overline{v}_{20}\hspace{0.5mm}, \hspace{0.5cm} L \hspace{0.5mm}\overline{v}_{210} =\frac{d \hspace{0.2mm}\overline{v}_{20}}{d \hspace{0.2mm}x} + 2 \hspace{0.5mm}\overline{v}_{11}\hspace{0.5mm}, \hspace{0.5cm} L \hspace{0.5mm}\overline{v}_{111} =\frac{d \hspace{0.2mm}\overline{v}_{11}}{d \hspace{0.2mm}x} \hspace{0.5mm}, } \vspace{0.4cm}\\
  \mathcal{O}(\mu^4):
  & \displaystyle{ v_4(x,T) = \left(\frac{d^4 \hspace{0.2mm}l_2}{d \hspace{0.2mm}T^4}\right) \hspace{1.0mm}\overline{v}_{4000}(x-l_2) - \left(\frac{d^3 \hspace{0.2mm}l_2}{d \hspace{0.2mm}T^3}\right)\hspace{0.2mm}\left(\frac{d \hspace{0.2mm}l_2}{d \hspace{0.2mm}T}\right) \hspace{1.0mm} \overline{v}_{3100}(x-l_2) } \vspace{0.2cm}\\
 & \displaystyle{ \hspace{1.8cm} - \left(\frac{d^2 \hspace{0.2mm}l_2}{d \hspace{0.2mm}T^2}\right)^2 \hspace{1.0mm} \overline{v}_{2200}(x-l_2) + \left(\frac{d^2 \hspace{0.2mm}l_2}{d \hspace{0.2mm}T^2}\right)\hspace{0.2mm}\left(\frac{d \hspace{0.2mm}l_2}{d \hspace{0.2mm}T}\right)^2 \hspace{1.0mm} \overline{v}_{2110}(x-l_2) } \vspace{0.2cm}\\
  & \displaystyle{ \hspace{2.3cm} - \left(\frac{d \hspace{0.2mm}l_2}{d \hspace{0.2mm}T}\right)^4 \hspace{1.0mm} \overline{v}_{1111}(x-l_2) - [\hspace{0.5mm} \mbox{terms of } \hspace{1.0mm}l_1 \hspace{0.5mm}] \hspace{0.5mm}, } \vspace{0.3cm}\\
  & \left \{
  \begin{array}{l}
    \displaystyle{ L \hspace{0.5mm}\overline{v}_{4000} = \overline{v}_{300}, \hspace{0.2cm}
      L \hspace{0.5mm}\overline{v}_{3100} = \frac{d \hspace{0.2mm}\overline{v}_{300} }{d \hspace{0.2mm}x} + \overline{v}_{210} = 0, \hspace{0.2cm} L \hspace{0.5mm}\overline{v}_{2200} = \overline{v}_{210},} \vspace{0.2cm}\\
    \displaystyle{ L \hspace{0.5mm}\overline{v}_{2110} = \frac{d \hspace{0.2mm}\overline{v}_{210} }{d \hspace{0.2mm}x} + 3 \hspace{0.5mm}\overline{v}_{111}, \hspace{0.2cm}
      L \hspace{0.5mm}\overline{v}_{1111} = \frac{d \hspace{0.2mm}\overline{v}_{111}}{d \hspace{0.2mm}x}.
    }
  \end{array}
  \right. \vspace{0.3cm}\\
\end{array} \nonumber
\end{eqnarray}
\begin{eqnarray}
  \label{eq_formulas_summary}
\end{eqnarray}
Here, to avoid redundancy, the notation $[\hspace{0.5mm} \mbox{terms of } \hspace{1.0mm}l_1 \hspace{0.5mm}]$ is introduced to denote the preceding $\overline{v}$\hspace{0.5mm}-\hspace{0.5mm}terms with $l_2$ replaced by $l_1$. Note that each equation for $\overline{v}_i$ is a nonhomogeneous linear ODE with constant coefficients. In particular, the nonhomogeneous term $F(x)$ for $\mathcal{O}(\mu^0)$ is piecewise linear in our case, which allows us to derive explicit formulas for not only  $\overline{v}_0$, but all $\overline{v}_i \hspace{1.5mm}(\hspace{0.5mm}i \ge 1\hspace{0.5mm})$. Substituting the $v_i$ in (\ref{eq_formulas_summary}) into the expansion (\ref{eq_expansion1}) and rewriting the differentiation with the overdot $\hspace{1.5mm}\dot{} \hspace{1.5mm}:= \hspace{0.5mm}d\hspace{0.5mm}/d\hspace{0.5mm}T$, we obtain
\begin{eqnarray}
  \begin{array}{lll}
    & v(x,T) = \displaystyle v_0(x,T) + \mu v_1 (x,T)+ \mu^2 v_2 (x,T) + \mu^3 v_3 (x,T)+ \cdots \vspace{0.2cm}\\
    & \displaystyle \hspace{0.2cm} =  -1\hspace{0.2mm}/\hspace{0.2mm}2+\Delta_0(x)+\overline{v}_0(x-l_2)+\hspace{0.5mm}\mu \hspace{0.5mm} \dot{l}_2 \hspace{1.0mm}\overline{v}_1(x-l_2) \vspace{0.1cm} \\
    & \displaystyle \hspace{0.6cm} +\hspace{0.5mm}\mu^2 \hspace{0.5mm} \left \{ \hspace{0.5mm} \ddot{l}_2 \hspace{1.0mm}\overline{v}_{20}(x-l_2)  - \dot{l}^{\hspace{0.5mm}2}_2 \hspace{1.0mm} \overline{v}_{11}(x-l_2) \hspace{0.5mm} \right \} \vspace{0.1cm} \\
    & \displaystyle \hspace{0.8cm} +\hspace{0.5mm}\mu^3 \hspace{0.5mm} \left \{ \hspace{0.5mm} \dddot{l}_2 \hspace{1.0mm}\overline{v}_{300}(x-l_2) - \ddot{l}_2 \hspace{0.2mm}\dot{l}_2 \hspace{1.0mm} \overline{v}_{210}(x-l_2) + \dot{l}^{\hspace{0.5mm}3}_2 \hspace{1.0mm} \overline{v}_{111}(x-l_2) \hspace{0.5mm} \right \} \vspace{0.1cm} \\
    & \displaystyle \hspace{1.0cm} +\hspace{0.5mm}\mu^4 \hspace{0.5mm} \left \{ \hspace{0.5mm} \ddddot{l}_2 \hspace{1.0mm}\overline{v}_{4000}(x-l_2) - \dddot{l}_2\hspace{0.2mm}\dot{l}_2  \hspace{1.0mm} \overline{v}_{3100}(x-l_2) - \ddot{l}^{\hspace{0.5mm}2}_2 \hspace{1.0mm} \overline{v}_{2200}(x-l_2) \right. \vspace{0.1cm}\\
    & \left. \displaystyle \hspace{3.8cm} + \ddot{l}_2\hspace{0.2mm}\dot{l}^{\hspace{0.5mm}2}_2 \hspace{1.0mm} \overline{v}_{2110}(x-l_2) - \dot{l}^{\hspace{0.5mm}4}_2 \hspace{1.0mm} \overline{v}_{1111}(x-l_2) \hspace{0.5mm} \right \} +\cdots \vspace{0.1cm}\\
   & \hspace{8.2cm} - [\hspace{0.5mm} \mbox{terms of } \hspace{1.0mm}l_1 \hspace{0.5mm}] \hspace{0.5mm}.
  \end{array} \label{eq_expansion3} 
\end{eqnarray}
Among the many terms in Eq.~(\ref{eq_expansion3}), we choose the following and truncate the others as higher-order terms:
\begin{eqnarray}
  \begin{array}{lll}
    & v(x,T) \vspace{0.2cm}\\
    & \displaystyle \hspace{0.1cm} =  -1\hspace{0.2mm}/\hspace{0.2mm}2+\Delta_0(x)+\overline{v}_0(x-l_2)+\hspace{0.5mm}\mu \hspace{0.5mm} \dot{l}_2 \hspace{1.0mm}\overline{v}_1(x-l_2) -\hspace{0.5mm}\mu^2 \hspace{0.5mm} \dot{l}^{\hspace{0.5mm}2}_2 \hspace{1.0mm} \overline{v}_{11}(x-l_2) \vspace{0.1cm}\\
    & \hspace{4.2cm}+\hspace{0.5mm}\mu^3 \hspace{0.5mm}\dot{l}^{\hspace{0.5mm}3}_2 \hspace{1.0mm} \overline{v}_{111}(x-l_2) -\hspace{0.5mm}\mu^4 \hspace{0.5mm} \dot{l}^{\hspace{0.5mm}4}_2 \hspace{1.0mm} \overline{v}_{1111}(x-l_2) \hspace{0.5mm} +\cdots \vspace{0.1cm}\\
   & \hspace{1.0cm} +\hspace{0.5mm}\mu^2 \hspace{0.5mm} \ddot{l}_2 \hspace{1.0mm}\overline{v}_{20}(x-l_2) - \hspace{0.5mm}\mu^3 \hspace{0.5mm} \ddot{l}_2\hspace{0.2mm}\dot{l}_2 \hspace{1.0mm} \overline{v}_{210}(x-l_2) + \hspace{0.5mm}\mu^4 \hspace{0.5mm} \ddot{l}_2\hspace{0.2mm}\dot{l}^{\hspace{0.5mm}2}_2 \hspace{1.0mm} \overline{v}_{2110}(x-l_2)  \hspace{0.5mm} +\cdots  \vspace{0.1cm}\\
   & \hspace{9.5cm} - [\hspace{0.5mm} \mbox{terms of } \hspace{1.0mm}l_1 \hspace{0.5mm}] \hspace{0.5mm} \vspace{0.2cm}\\
   & \displaystyle \hspace{0.2cm} =  -1\hspace{0.2mm}/\hspace{0.2mm}2+\Delta_0(x)+\overline{U}_2^{(1)}(x-l_2,T) +\hspace{0.5mm}\mu^2 \hspace{0.5mm} \ddot{l}_2 \hspace{1.0mm} \overline{U}_2^{(2)}(x-l_2,T) \vspace{0.2cm}\\
   & \hspace{5.0cm} -\hspace{0.5mm}\overline{U}_1^{(1)}(x-l_1,T) -\hspace{0.5mm}\mu^2 \hspace{0.5mm} \ddot{l}_1 \hspace{1.0mm} \overline{U}_1^{(2)}(x-l_1,T)\hspace{0.5mm}, \vspace{0.2cm}
  \end{array} \nonumber 
\end{eqnarray}
\begin{eqnarray}
  \label{eq_expansion4}
\end{eqnarray}
where we have defined
\begin{eqnarray}
    \begin{array}{lll}
      & \overline{U}_i^{(1)}(x,T) \hspace{0.5mm}=\hspace{0.5mm} \overline{v}_0(x)+\hspace{0.5mm}\mu \hspace{0.5mm} \dot{l}_i \hspace{1.0mm}\overline{v}_1(x) -\hspace{0.5mm}\mu^2 \hspace{0.5mm} \dot{l}^{\hspace{0.5mm}2}_i \hspace{1.0mm} \overline{v}_{11}(x) +\hspace{0.5mm}\mu^3 \hspace{0.5mm}\dot{l}^{\hspace{0.5mm}3}_i \hspace{1.0mm} \overline{v}_{111}(x) \vspace{0.2cm}\\
    & \hspace{7.2cm}-\hspace{0.5mm}\mu^4 \hspace{0.5mm} \dot{l}^{\hspace{0.5mm}4}_i \hspace{1.0mm} \overline{v}_{1111}(x) \hspace{0.5mm} +\cdots  \hspace{0.5mm}, \vspace{0.3cm}\\
      & \overline{U}_i^{(2)}(x,T) \hspace{0.5mm}=\hspace{0.5mm} \overline{v}_{20}(x) - \hspace{0.5mm}\mu \hspace{0.5mm} \dot{l}_i \hspace{1.0mm} \overline{v}_{210}(x) + \hspace{0.5mm}\mu^2 \hspace{0.5mm} \dot{l}^{\hspace{0.5mm}2}_i \hspace{1.0mm} \overline{v}_{2110}(x)  \hspace{0.5mm} +\cdots \hspace{0.5mm}, \vspace{0.35cm}\\
      & \displaystyle \hspace{9.0cm} (\hspace{0.5mm}i=1,\hspace{0.5mm}2\hspace{0.5mm}) 
  \end{array} \label{eq_function_U}
\end{eqnarray}
Note that both $\overline{U}_i^{(1)}(x,T)$ and $\overline{U}_i^{(2)}(x,T)$ are expanded in power series of $\mu \hspace{0.5mm}l_i$, which can be written as follows in renormalized form. First, applying the operator $L$ to $\overline{U}_i^{(1)}$ and using the relations between $\overline{v}$ given by the lemma in Appendix B, we obtain (in the following, the argument $x$ for $\overline{v}_i$ is sometimes omitted, unless this may cause confusion):
\begin{eqnarray}
    \begin{array}{lll}
      & L \hspace{0.5mm} \overline{U}_i^{(1)}(x,T) \vspace{0.1cm}\\
      & \displaystyle \hspace{0.2cm}\hspace{0.5mm}=\hspace{0.5mm} L \hspace{0.5mm} \overline{v}_0(x)-(\hspace{0.5mm}-\hspace{0.5mm}\mu \hspace{0.5mm}l_i \hspace{0.5mm}) \hspace{1.0mm} L \hspace{0.5mm}\overline{v}_1(x) -(\hspace{0.5mm}-\hspace{0.5mm}\mu \hspace{0.5mm}l_i \hspace{0.5mm})^2\hspace{0.5mm} \hspace{1.0mm} L \hspace{0.5mm} \overline{v}_{11}(x) \vspace{0.1cm}\\
      & \hspace{2.5cm} -(\hspace{0.5mm}-\hspace{0.5mm}\mu \hspace{0.5mm}l_i \hspace{0.5mm})^3 \hspace{1.0mm} L \hspace{0.5mm} \overline{v}_{111}(x) -(\hspace{0.5mm}-\hspace{0.5mm}\mu \hspace{0.5mm}l_i \hspace{0.5mm})^4 \hspace{1.0mm} L \hspace{0.5mm} \overline{v}_{1111}(x) \hspace{0.5mm} +\cdots \vspace{0.3cm}\\
      & \displaystyle \hspace{0.2cm} = \hspace{0.5mm}-F(x)+(\hspace{0.5mm}-\hspace{0.5mm}\mu \hspace{0.5mm}l_i \hspace{0.5mm}) \hspace{1.0mm} \left(\hspace{0.5mm}\frac{d \hspace{0.2mm}\overline{v}_{0}}{d \hspace{0.2mm}x} \hspace{0.5mm}\right) -(\hspace{0.5mm}-\hspace{0.5mm}\mu \hspace{0.5mm}l_i \hspace{0.5mm})^2\hspace{0.5mm} \hspace{1.0mm} \left(\hspace{0.5mm}\frac{d \hspace{0.2mm}\overline{v}_{1}}{d \hspace{0.2mm}x} \hspace{0.5mm}\right) \vspace{0.2cm}\\
      & \displaystyle \hspace{2.2cm} -(\hspace{0.5mm} -\hspace{0.5mm}\mu \hspace{0.5mm}l_i \hspace{0.5mm})^3 \hspace{1.0mm} \left(\hspace{0.5mm}\frac{d \hspace{0.2mm}\overline{v}_{11}}{d \hspace{0.2mm}x} \hspace{0.5mm}\right) -(\hspace{0.5mm}-\hspace{0.5mm}\mu \hspace{0.5mm}l_i \hspace{0.5mm})^4 \hspace{1.0mm} \left(\hspace{0.5mm}\frac{d \hspace{0.2mm}\overline{v}_{111}}{d \hspace{0.2mm}x} \hspace{0.5mm}\right) \hspace{0.5mm} +\cdots \vspace{0.2cm}\\
      & \displaystyle \hspace{0.2cm} = \hspace{0.5mm}-F(x)-\hspace{0.5mm}\mu \hspace{0.5mm}l_i \hspace{1.0mm} \left(\hspace{0.5mm}\frac{d \hspace{0.2mm}\overline{U}^{(1)}_{i}}{\hspace{0.5mm}d \hspace{0.2mm}x} \hspace{0.5mm}\right) \hspace{0.5mm} 
   \end{array} \label{eq_renormalization1-1}
\end{eqnarray}
or, equivalently,
\begin{eqnarray}
    \begin{array}{l}
      \tilde{L}\hspace{1.0mm} \overline{U}_i^{(1)}(x,T) \hspace{0.5mm} = \hspace{0.5mm}-F(x) \hspace{0.5mm}, 
   \end{array} \label{eq_renormalization1-2}
\end{eqnarray}
where
\begin{eqnarray}
    \begin{array}{l}
      \displaystyle \tilde{L} \hspace{1.0mm}:= \hspace{1.0mm}D \hspace{0.5mm} \frac{d^2}{\hspace{0.5mm}d \hspace{0.2mm}x^2} + \hspace{0.5mm}\mu \hspace{0.5mm}l_i \hspace{0.5mm} \frac{d}{\hspace{0.5mm}d \hspace{0.2mm}x} -1 \hspace{0.5mm}. 
   \end{array} \label{eq_linear_operator2}
\end{eqnarray}
This is again a nonhomogeneous linear differential equation with constant coefficients, which can be solved explicitly as
\begin{eqnarray}
\begin{array}{rcl}
\overline{U}_i^{(1)} & = & \displaystyle{
\left \{
\begin{array}{cc}
\displaystyle
\hspace{1.2mm}\frac{1}{2} - \frac{\mu \hspace{0.5mm}\dot{l}_i + \phi (\dot{l}_i) }{2 \phi (\dot{l}_i)}
\exp \left( \frac{ -\mu \hspace{0.5mm}\dot{l}_i + \phi (\dot{l}_i) }{2D} \hspace{0.5mm}x \right) \hspace{0.5mm}, &\hspace{0.3cm}(\hspace{0.5mm} x \le 0\hspace{0.5mm}) \vspace{0.3cm}\\
\displaystyle
\hspace{0.2mm}-\frac{1}{2 } + \frac{-\mu \hspace{0.5mm}\dot{l}_i + \phi(\dot{l}_i)}{2 \phi(\dot{l}_i)}
\exp \left( -\frac{\mu \hspace{0.5mm}\dot{l}_i + \phi(\dot{l}_i) }{2D} \hspace{0.5mm}x \right) \hspace{0.5mm}, &\hspace{0.3cm}(\hspace{0.5mm} x > 0\hspace{0.5mm})
\end{array} \label{eq_renormalized_U1}
\right.
}
\end{array}
\end{eqnarray}
where we have defined $\phi(\dot{l}_i) := \sqrt{(\mu \hspace{0.5mm}\dot{l}_i)^{\hspace{0.5mm}2} + 4D}$.

Similarly, applying the operator $L$ to $\overline{U}_i^{(2)}$ and using the relations given by the lemma in Appendix B, we obtain 
\begin{eqnarray}
    \begin{array}{lll}
      & \displaystyle L \hspace{0.5mm} \overline{U}_i^{(2)}(x,T) \vspace{0.1cm}\\
      & \displaystyle \hspace{0.2cm}\hspace{0.5mm}=\hspace{0.5mm} L \hspace{0.5mm} \overline{v}_{20}(x)+(\hspace{0.5mm}-\hspace{0.5mm}\mu \hspace{0.5mm}l_i \hspace{0.5mm}) \hspace{1.0mm} L \hspace{0.5mm}\overline{v}_{210}(x) +(\hspace{0.5mm}-\hspace{0.5mm}\mu \hspace{0.5mm}l_i \hspace{0.5mm})^2\hspace{0.5mm} \hspace{1.0mm} L \hspace{0.5mm} \overline{v}_{2110}(x) \vspace{0.1cm}\\
      & \displaystyle \hspace{5.5cm} +(\hspace{0.5mm}-\hspace{0.5mm}\mu \hspace{0.5mm}l_i \hspace{0.5mm})^3\hspace{0.5mm} \hspace{1.0mm} L \hspace{0.5mm} \overline{v}_{21110}(x) \hspace{0.5mm} +\cdots \vspace{0.2cm}\\
      & \displaystyle \hspace{0.2cm} = \hspace{0.5mm}(\hspace{0.5mm}-\hspace{0.5mm}\mu \hspace{0.5mm}l_i \hspace{0.5mm}) \hspace{1.0mm} \left(\hspace{0.5mm}\frac{d \hspace{0.2mm}\overline{v}_{20}}{d \hspace{0.2mm}x} \hspace{0.5mm}\right) +(\hspace{0.5mm}-\hspace{0.5mm}\mu \hspace{0.5mm}l_i \hspace{0.5mm})^2\hspace{0.5mm} \hspace{1.0mm} \left(\hspace{0.5mm}\frac{d \hspace{0.2mm}\overline{v}_{210}}{d \hspace{0.2mm}x} \hspace{0.5mm}\right) \vspace{0.2cm}\\
      & \displaystyle \hspace{4.5cm} -(\hspace{0.5mm} -\hspace{0.5mm}\mu \hspace{0.5mm}l_i \hspace{0.5mm})^3 \hspace{1.0mm} \left(\hspace{0.5mm}\frac{d \hspace{0.2mm}\overline{v}_{2110}}{d \hspace{0.2mm}x} \hspace{0.5mm}\right) \hspace{0.5mm} +\cdots \vspace{0.2cm}\\      
      & \displaystyle \hspace{0.5cm} + \left \{\hspace{0.5mm} \overline{v}_{1}+(\hspace{0.5mm}-\hspace{0.5mm}\mu \hspace{0.5mm}l_i \hspace{0.5mm}) \hspace{1.0mm} 2 \hspace{0.5mm} \overline{v}_{11} +(\hspace{0.5mm}-\hspace{0.5mm}\mu \hspace{0.5mm}l_i \hspace{0.5mm})^2\hspace{0.5mm} \hspace{1.0mm} 3 \hspace{0.5mm} \overline{v}_{111} +(\hspace{0.5mm} -\hspace{0.5mm}\mu \hspace{0.5mm}l_i \hspace{0.5mm})^3 \hspace{1.0mm} 4 \hspace{0.5mm} \overline{v}_{1111} \hspace{0.5mm} +\cdots \hspace{0.5mm} \right \} \vspace{0.2cm}\\
      & \displaystyle \hspace{0.2cm} = (\hspace{0.5mm}-\hspace{0.5mm}\mu \hspace{0.5mm}l_i\hspace{0.5mm}) \hspace{1.0mm} \left(\hspace{0.5mm}\frac{d \hspace{0.2mm}\overline{U}^{(2)}_{i}}{\hspace{0.5mm}d \hspace{0.2mm}x} \hspace{0.5mm}\right) +\overline{V}^{(2)}_{i} \hspace{0.5mm}, 
   \end{array} \label{eq_renormalization2-1}
\end{eqnarray}
where $\overline{V}_i^{(2)} := \overline{v}_{1}+(\hspace{0.5mm}-\hspace{0.5mm}\mu \hspace{0.5mm}l_i \hspace{0.5mm}) \hspace{1.0mm} 2 \hspace{0.5mm} \overline{v}_{11} +(\hspace{0.5mm}-\hspace{0.5mm}\mu \hspace{0.5mm}l_i \hspace{0.5mm})^2\hspace{0.5mm} \hspace{1.0mm} 3 \hspace{0.5mm} \overline{v}_{111} +(\hspace{0.5mm} -\hspace{0.5mm}\mu \hspace{0.5mm}l_i \hspace{0.5mm})^3 \hspace{1.0mm} 4 \hspace{0.5mm} \overline{v}_{1111} \hspace{0.5mm} +\cdots \hspace{0.5mm}$.
Applying the operator $L$ to $\overline{V}_i^{(2)}$ yields
\begin{eqnarray}
\begin{array}{lll}
 \displaystyle L  \hspace{0.5mm}\overline{V}_i^{(2)} (x,T)= \displaystyle (\hspace{0.5mm}-\hspace{0.5mm}\mu \hspace{0.5mm}l_i\hspace{0.5mm}) \hspace{1.0mm} \left(\hspace{0.5mm}\frac{d \hspace{0.2mm}\overline{V}^{(2)}_{i}}{\hspace{0.5mm}d \hspace{0.2mm}x} \hspace{0.5mm}\right) \hspace{0.1mm}-\frac{d \hspace{0.2mm}\overline{W}^{(2)}_{i}}{\hspace{0.5mm}d \hspace{0.2mm}x}  \hspace{0.5mm},
\end{array} \label{eq_renormalization2-2}
\end{eqnarray}
where $\overline{W}_i^{(2)} := \overline{v}_0 - (\hspace{0.5mm}-\hspace{0.5mm}\mu \hspace{0.5mm}l_i\hspace{0.5mm}) \hspace{1.0mm}\overline{v}_1 - (\hspace{0.5mm}-\hspace{0.5mm}\mu \hspace{0.5mm}l_i\hspace{0.5mm})^2 \hspace{1.0mm}\overline{v}_{11} - (\hspace{0.5mm}-\hspace{0.5mm}\mu \hspace{0.5mm}l_i\hspace{0.5mm})^3 \hspace{1.0mm}\overline{v}_{111} + \cdots $.
Applying $L$ to $\overline{W}_i^{(2)}$ yields
\begin{eqnarray}
\begin{array}{lll}
 \displaystyle L \hspace{0.5mm}\overline{W}_i^{(2)}(x,T)  = \displaystyle (\hspace{0.5mm}-\hspace{0.5mm}\mu \hspace{0.5mm}l_i\hspace{0.5mm})  \hspace{0.1mm}\left(\frac{\hspace{0.5mm} d \hspace{0.2mm}\overline{W}^{(2)}_{i}}{\hspace{0.5mm}d \hspace{0.2mm}x}\hspace{0.5mm}\right) - F(x). 
\end{array} \label{eq_renormalization2-2}
\end{eqnarray}
Thus, we have three linear ODEs
\begin{eqnarray}
\tilde{L} \hspace{0.5mm}\overline{U}_i^{(2)}(x,T) = \overline{V}_i^{(2)}, \hspace{0.1cm}
\tilde{L} \hspace{0.5mm}\overline{V}_i^{(2)}(x,T) = - \frac{d \overline{W}_i^{(2)}}{d x}, \hspace{0.1cm}
\tilde{L} \hspace{0.5mm}\overline{W}_i^{(2)}(x,T) = - F(x),
\end{eqnarray}
which are simultaneously solved to yield \vspace{0.2cm}
\begin{eqnarray} 
\begin{array}{lll}
\overline{U}_i^{(2)} =  \displaystyle{
\left \{
\begin{array}{cc}
\displaystyle
\hspace{0.1cm} \left(\hspace{0.5mm}
\frac{6D^2}{\phi(\dot{l}_i)^5}
-\frac{3D x}{\phi(\dot{l}_i)^4}
+\frac{x^2}{2\phi(\dot{l}_i)^3}
\hspace{0.5mm} \right) \hspace{1.0mm}
\exp \left( \frac{ -\mu \hspace{0.5mm}\dot{l}_i + \phi (\dot{l}_i) }{2D} \hspace{0.5mm}x \right) \hspace{0.5mm}, \hspace{0.5cm}(\hspace{0.5mm} x \le 0\hspace{0.5mm}) \vspace{0.3cm}\\
\displaystyle
\hspace{0.1cm} \left(\hspace{0.5mm}
\frac{6D^2}{\phi(\dot{l}_i)^5}
+\frac{3D x}{\phi(\dot{l}_i)^4}
+\frac{x^2}{2\phi(\dot{l}_i)^3}
\hspace{0.5mm} \right) \hspace{1.0mm}
\exp \left( -\frac{\mu \hspace{0.5mm}\dot{l}_i + \phi(\dot{l}_i) }{2D} \hspace{0.5mm}x \right)  \hspace{0.5mm}, \hspace{0.5cm}(\hspace{0.5mm} x > 0\hspace{0.5mm})
\end{array}
\right. \nonumber
} \vspace{0.2cm}
\end{array}
\end{eqnarray}
\begin{eqnarray}
  \label{eq_renormalized_U2}
\end{eqnarray}
%

Now that we have perturbatively solved for $v(x,T)$, we derive equations of motion from Eqs.~(\ref{eq_hybrid}) and (\ref{eq_function_U}). In terms of the slow timescale, the first equation in (\ref{eq_hybrid}) can be rewritten as
\begin{eqnarray}
\begin{array}{l}
  \displaystyle \mu \hspace{0.5mm}\frac{d \hspace{0.5mm}l_2(T)}{d \hspace{0.5mm}T}=-\frac{v(l_2,T)}{\sqrt{2} \hspace{0.5mm}\tau} \hspace{0.5mm}. \\
\end{array} \label{eq_reduced_ODE1}
\end{eqnarray}
Substituting Eq.~(\ref{eq_expansion4}) on the right-hand side of Eq.~(\ref{eq_reduced_ODE1}), we have
\begin{eqnarray}
  \begin{array}{lll}
    &\displaystyle -\sqrt{2} \hspace{0.5mm}\tau \hspace{0.3mm}\mu \hspace{0.5mm}\dot{l}_2 = v(l_2,T) \\
    &\hspace{1.7cm} = -1\hspace{0.2mm}/\hspace{0.2mm}2+\Delta_0(l_2)+\overline{U}_2^{(1)}(0,T) +\hspace{0.5mm}\mu^2 \hspace{0.5mm} \ddot{l}_2 \hspace{1.0mm} \overline{U}_2^{(2)}(0,T) \vspace{0.2cm}\\
    &\hspace{4.8cm} -\hspace{0.5mm}\overline{U}_1^{(1)}(h,T) -\hspace{0.5mm}\mu^2 \hspace{0.5mm} \ddot{l}_1 \hspace{1.0mm} \overline{U}_1^{(2)}(h,T)\hspace{0.5mm}, 
\end{array} \label{eq_reduced_ODE2}
\end{eqnarray}
where $h:=l_2-l_1$. The last four terms in Eq.~(\ref{eq_reduced_ODE2}) are computed from Eqs.~(\ref{eq_renormalized_U1}) and (\ref{eq_renormalized_U2}) as
\begin{eqnarray}
\begin{array}{lll}
  \displaystyle \overline{U}_2^{(1)}(0,T) = - \frac{\mu \hspace{0.5mm}\dot{l}_2}{2 \phi(\dot{l}_2)}, \hspace{0.5cm}
    \overline{U}_2^{(2)}(0,T) = \frac{6 D^2}{\phi(\dot{l}_2)^5}, \vspace{0.3cm}\\
  \displaystyle \overline{U}_1^{(1)}(h,T) = - \frac{1}{2} + \frac{-\mu \hspace{0.5mm}\dot{l}_1 + \phi(\dot{l}_1)}{2 \phi(\dot{l}_1)}
     \exp \left(\hspace{0.5mm} - \frac{\mu \hspace{0.5mm}\dot{l}_1 + \phi(\dot{l}_1) }{2D} h \hspace{0.5mm}\right), \vspace{0.3cm} \\
  \displaystyle \overline{U}_1^{(2)}(h,T)  =  \left(\hspace{0.5mm} \frac{6 D^2}{\phi(\dot{l}_1)^5}+ \frac{3 D h}{\phi(\dot{l}_1)^4} + \frac{h^2}{2 \phi(\dot{l}_1)^3}  \hspace{0.5mm}\right)
     \hspace{0.5mm} \exp \left( \hspace{0.5mm}- \frac{\mu \hspace{0.5mm}\dot{l}_1 + \phi(\dot{l}_1) }{2D} h \hspace{0.5mm} \right). 
\end{array} \label{eq_reduced_ODE3}
\end{eqnarray}
Substituting these into Eq.~(\ref{eq_reduced_ODE2}) and using the variables $r_2 := \mu \hspace{0.5mm}\dot{l}_2$ and $r_1 := \mu \hspace{0.5mm}\dot{l}_1$ leads to \vspace{0.35cm}
\begin{eqnarray}
  \begin{array}{lll}
   \displaystyle \frac{6 D^2}{\phi(r_2)^5} \hspace{0.5mm} \mu \hspace{0.5mm}\dot{r}_2
    -\left(\hspace{0.5mm} \frac{6 D^2}{\phi(r_1)^5}+ \frac{3 D h}{\phi(r_1)^4} + \frac{h^2}{2 \phi(r_1)^3}  \hspace{0.5mm}\right)
      \hspace{0.5mm} \exp \left( \hspace{0.5mm}- \frac{ r_1 + \phi(r_1) }{2D} h \hspace{0.5mm} \right) \hspace{0.5mm} \mu \hspace{0.5mm} \dot{r}_1 \vspace{0.35cm}\\
      \displaystyle = -\sqrt{2} \hspace{0.5mm}\tau \hspace{0.3mm}r_2+\frac{r_2}{2 \phi(r_2)}+ \frac{-r_1 + \phi(r_1)}{2 \phi(r_1)}
    \exp \left(\hspace{0.5mm} - \frac{r_1 + \phi(r_1) }{2D} h \hspace{0.5mm}\right)-\Delta_0(l_2) \hspace{0.5mm}. 
\end{array} \label{eq_reduced_ODE4} 
\end{eqnarray} \vspace{0.35cm}
Finally, by redefining the time and the time derivative as $t:=T/\hspace{0.5mm}\mu$ and $\hspace{1.0mm}\dot{}\hspace{1.0mm}:=d/d\hspace{0.5mm}t$, (\ref{eq_reduced_ODE4}) reduces to
\begin{eqnarray}
  \begin{array}{ll}
   &\hspace{0.3cm} \displaystyle m(r_2) \hspace{0.5mm}\dot{r}_2
    -M(r_2,h) \hspace{0.5mm} \exp \left( \hspace{0.5mm}- \frac{ r_1 + \phi(r_1) }{2D} h \hspace{0.5mm} \right) \hspace{0.5mm} \dot{r}_1 \vspace{0.3cm}\\
   &\displaystyle = g(r_2)+ G(-r_1) \hspace{0.5mm} \exp \left(\hspace{0.5mm} - \frac{r_1 + \phi(r_1) }{2D} h \hspace{0.5mm}\right)-\Delta_0(l_2) \hspace{0.5mm}, 
\end{array} \label{eq_reduced_ODE5} 
\end{eqnarray} 
where the functions are defined as
\begin{eqnarray}
\begin{array}{cc}
M(r,h) := \displaystyle \left( \hspace{0.5mm} \frac{6D^2}{\phi(r)^5} + \frac{3Dh}{\phi(r)^4} + \frac{h^2}{2\phi(r)^3} \hspace{0.5mm}\right), \vspace{0.2cm} \\
\displaystyle \hspace{1.0cm} m(r) := \frac{6 D^2}{\phi(r)^5}, \hspace{0.5cm} g(r) :=  -\sqrt{2}\hspace{0.2mm} \tau \hspace{0.2mm}r + \frac{r}{2 \phi(r)}, \hspace{0.5cm} G(r) := \frac{r + \phi(r)}{2 \phi(r)}.
\end{array} \label{eq_prefactor_function_appendix}
\end{eqnarray}
Similarly, the second equation in (\ref{eq_hybrid}) can also be transformed; together with the aforementioned $\dot{l}_2=r_2$ and $\dot{l}_1=r_1$, this gives the reduced ODE system in (\ref{eq_ODE_renormalized_full}).
%

%
\clearpage
\section{Relations between function $\overline{v}$ in Eq.~(\ref{eq_formulas_summary})}
By setting $d\hspace{0.2mm}v_i/d \hspace{0.2mm}t \hspace{0.2mm}=0$ for Eq.~(\ref{eq_order}) in Appendix A, we find that
\begin{eqnarray}
  \begin{array}{l}
    L \hspace{0.5mm}v_{n+1}(x,T)=\partial \hspace{0.5mm}v_{n}/\partial \hspace{0.5mm}T\hspace{0.5mm}
\end{array} \label{eq_order_appendix_B}
\end{eqnarray}
holds for $n \ge 1$, where $L:= D (\hspace{0.5mm}d^2 \hspace{0.2mm}/d \hspace{0.2mm}x^2 \hspace{0.5mm})-1$. We can show that terms of the form $(d \hspace{0.2mm}l_2/d \hspace{0.2mm}T)^n$ and $(d^2 \hspace{0.2mm}l_2/d \hspace{0.2mm}T^2)\hspace{0.5mm}(d \hspace{0.2mm}l_2/d \hspace{0.2mm}T)^{n-2}$ always appear in the expressions for $v_n$ when $n \ge 2$, and that their associated functions satisfy particular relations. \vspace{0.5cm}\\
\noindent {\bf Lemma \hspace{0.2mm}B\hspace{0.2mm}1.}\\
{\it 
For $n\ge 2$, the $n$-th order function $v_n(x,T)$ in Eq.~(\ref{eq_formulas_summary}) includes terms of the form $(-1)^{n+1} \hspace{0.5mm}(d \hspace{0.2mm}l_2/d \hspace{0.2mm}T)^n$ and $(-1)^{n} \hspace{0.5mm}(d^2 \hspace{0.2mm}l_2/d \hspace{0.2mm}T^2)\hspace{0.5mm}(d \hspace{0.2mm}l_2/d \hspace{0.2mm}T)^{n-2}$.
Furthermore, let their associated functions be denoted by $\displaystyle{ \overline{v}_{\scalebox{0.7}{$\underbrace{111 \cdots 11}_{\hspace{-0.3cm} \scalebox{1.8}{$n$}}$}} }$ and $\displaystyle{ \overline{v}_{2\scalebox{0.7}{$\underbrace{11 \cdots 1}_{\scalebox{1.8}{$n-2$}}$}0} }$, respectively. That is, 
\begin{eqnarray}
\begin{array}{ll}
  \mathcal{O}(\hspace{0.5mm} \mu^n \hspace{0.5mm}): \\
  \displaystyle \hspace{0.5cm} v_n(x,T) = \Pi_{n}(x-l_2,T) +(-1)^{n+1}\hspace{0.5mm} \left(d \hspace{0.2mm}l_2/d \hspace{0.2mm}T \right)^n \hspace{1.0mm} \overline{v}_{\scalebox{0.7}{$\underbrace{111 \cdots 11}_{\hspace{-0.3cm} \scalebox{1.8}{$n$}}$}}(x-l_2) \\
  \displaystyle \hspace{2.5cm} +(-1)^{n}\hspace{0.5mm}\left(d^2 \hspace{0.2mm}l_2/d \hspace{0.2mm}T^2 \right)\hspace{0.5mm}\left(d \hspace{0.2mm}l_2/d \hspace{0.2mm}T \right)^{n-2} \hspace{1.0mm}\overline{v}_{2\scalebox{0.7}{$\underbrace{11 \cdots 1}_{\scalebox{1.8}{$n-2$}}$}0}(x-l_2) \\
  \displaystyle \hspace{8.5cm} - [\hspace{0.5mm} \mbox{terms of } \hspace{1.0mm}l_1 \hspace{0.5mm}]\hspace{0.5mm},
\end{array} \label{eq_statement1}
\end{eqnarray}
where $\Pi_{n}(x-l_2,T)$ denotes all the other terms. The following relations then hold between the two functions: }
\begin{eqnarray}
\begin{array}{ll}
  \hspace{0.2cm} L \hspace{0.5mm} \overline{v}_{\scalebox{0.7}{$\underbrace{111 \cdots 11}_{\scalebox{1.8}{$n+1$}}$}}(x)= \frac{d}{d\hspace{0.3mm}x}  \hspace{0.5mm}\overline{v}_{\scalebox{0.7}{$\underbrace{111 \cdots 11}_{\hspace{-0.2cm}\scalebox{2.0}{$n$}}$}}(x) \hspace{0.5mm}, \vspace{0.2cm} \\
\hspace{0.2cm} L \hspace{0.5mm} \overline{v}_{2\scalebox{0.7}{$\underbrace{11\cdots 1}_{ \scalebox{1.8}{$n-1$}}$ }0}(x) = \frac{d}{d\hspace{0.3mm}x}  \hspace{0.5mm} \overline{v}_{2\scalebox{0.7}{$\underbrace{11\cdots 1}_{ \scalebox{1.8}{$n-2$} }$}0}(x) + n \hspace{0.5mm}\overline{v}_{\scalebox{0.7}{$\underbrace{111 \cdots 11}_{\hspace{-0.3cm} \scalebox{1.8}{$n$}}$}}(x) \hspace{0.5mm}. \vspace{0.2cm} \\
\end{array} \label{eq_statement1}
\end{eqnarray}
\vspace{0.2cm}\\
{\bf Proof.}\hspace{0.1cm}
{\it We prove the lemma using mathematical induction. From Eq.~(\ref{eq_formulas_summary}) in Appendix A, the statement is true for $n=2$. Assume that, for $n=k$ \hspace{0.1cm}$(\hspace{0.5mm}k \ge 2\hspace{0.5mm})$, the solution $v_k$ to Eq.~(\ref{eq_order_appendix_B}) is represented by
\begin{eqnarray}
\begin{array}{ll}
  \displaystyle \hspace{0.5cm} v_{k}(x,T) = \Pi_{k}(x-l_2,T) +(-1)^{k+1}\hspace{0.5mm} \left(d \hspace{0.2mm}l_2/d \hspace{0.2mm}T\right)^k \hspace{1.0mm} \overline{v}_{\scalebox{0.7}{$\underbrace{111 \cdots 11}_{\hspace{-0.3cm} \scalebox{1.8}{$k$}}$}}(x-l_2) \\
  \displaystyle \hspace{2.5cm} +(-1)^{k}\hspace{0.5mm}\left(d^2 \hspace{0.2mm}l_2/d \hspace{0.2mm}T^2 \right)\hspace{0.5mm}\left(d \hspace{0.2mm}l_2/d \hspace{0.2mm}T \right)^{k-2} \hspace{1.0mm}\overline{v}_{2\scalebox{0.7}{$\underbrace{11 \cdots 1}_{\scalebox{1.8}{$k-2$}}$}0}(x-l_2) \\
  \displaystyle \hspace{8.5cm} - [\hspace{0.5mm} \mbox{terms of } \hspace{1.0mm}l_1 \hspace{0.5mm}]\hspace{0.5mm}.
\end{array} \label{eq_proof1}
\end{eqnarray}
Differentiating both sides of Eq.~(\ref{eq_proof1}) with respect to $T$, we have
\begin{eqnarray}
 \begin{array}{ll}
  \displaystyle \hspace{0.5cm} d \hspace{0.5mm}v_{k}(x,T)/dT  \vspace{0.2cm}\\
  \displaystyle \hspace{0.7cm}= \tilde{\Pi}_{k}(x-l_2,T) \vspace{0.2cm}\\
  \displaystyle \hspace{1.4cm} +(-1)^{k+1}\hspace{0.5mm}\left(d^2 \hspace{0.2mm}l_2/d \hspace{0.2mm}T^2\right) \hspace{1.0mm}k \hspace{1.0mm}\left(d \hspace{0.2mm}l_2/d \hspace{0.2mm}T\right)^{k-1} \hspace{1.0mm}\overline{v}_{\scalebox{0.7}{$\underbrace{111 \cdots 11}_{\hspace{-0.3cm} \scalebox{1.8}{$k$}}$}}(x-l_2) \\
  \displaystyle \hspace{1.7cm} +(-1)^{k+2}\hspace{0.5mm} \left(d \hspace{0.2mm}l_2/d \hspace{0.2mm}T\right)^{k+1} \hspace{1.0mm} d \hspace{0.5mm}\overline{v}_{\scalebox{0.7}{$\underbrace{111 \cdots 11}_{\hspace{-0.3cm} \scalebox{1.8}{$k$}}$}}(x-l_2)/d\hspace{0.2mm}x \\
  \displaystyle \hspace{2.0cm} +(-1)^{k+1}\hspace{0.5mm}\left(d^2 \hspace{0.2mm}l_2/d \hspace{0.2mm}T^2\right)\hspace{0.5mm}\left(d \hspace{0.5mm}\hspace{0.2mm}l_2/d \hspace{0.2mm}T\right)^{k-1} \hspace{1.0mm}d \hspace{0.5mm}\overline{v}_{2\scalebox{0.7}{$\underbrace{11 \cdots 1}_{\scalebox{1.8}{$k-2$}}$}0}(x-l_2)/d\hspace{0.2mm}x \\
  \displaystyle \hspace{8.5cm} - [\hspace{0.5mm} \mbox{terms of } \hspace{1.0mm}l_1 \hspace{0.5mm}] \vspace{0.2cm}\\
  \displaystyle \hspace{0.7cm} =\tilde{\Pi}_{k}(x-l_2,T) +(-1)^{k+2}\hspace{0.5mm} \left(d \hspace{0.2mm}l_2/d \hspace{0.2mm}T \right)^{k+1} \hspace{1.0mm} d \hspace{0.5mm}\overline{v}_{\scalebox{0.7}{$\underbrace{111 \cdots 11}_{\hspace{-0.3cm} \scalebox{1.8}{$k$}}$}}(x-l_2)/d\hspace{0.2mm}x \\
  \displaystyle \hspace{2.3cm} +(-1)^{k+1}\hspace{0.5mm}\left(d^2 \hspace{0.2mm}l_2/d \hspace{0.2mm}T^2 \right)\hspace{0.5mm}(d \hspace{0.2mm}l_2/d \hspace{0.2mm}T)^{k-1} \vspace{0.2cm}\\
  \displaystyle \hspace{3.0cm} \times\{\hspace{0.5mm} d \hspace{0.5mm}\overline{v}_{2\scalebox{0.7}{$\underbrace{11 \cdots 1}_{\scalebox{1.8}{$k-2$}}$}0}(x-l_2)/d\hspace{0.2mm}x +\hspace{0.5mm}k \hspace{1.0mm}\overline{v}_{\scalebox{0.7}{$\underbrace{111 \cdots 11}_{\hspace{-0.3cm} \scalebox{1.8}{$k$}}$}}(x-l_2) \hspace{0.5mm}\} \vspace{0.2cm}\\
  \displaystyle \hspace{8.5cm} - [\hspace{0.5mm} \mbox{terms of } \hspace{1.0mm}l_1 \hspace{0.5mm}]\hspace{0.5mm},  
\end{array} \label{eq_proof2}
\end{eqnarray} \vspace{0.2cm}
where the irrelevant terms are included in $\tilde{\Pi}_{k}(x-l_2,T)$. However, it follows from $L\hspace{0.5mm} v_{k+1}=d\hspace{0.2mm}v_{k}/d\hspace{0.2mm}T$ that
\begin{eqnarray}
 \begin{array}{ll}
  \displaystyle \hspace{0.2cm}v_{k+1}(x,T) \vspace{0.2cm}\\
  \displaystyle \hspace{0.5cm} = L^{-1}\hspace{0.5mm} \tilde{\Pi}_{k}(x-l_2,T) +(-1)^{k+2}\hspace{0.5mm} (d \hspace{0.2mm}l_2/d \hspace{0.2mm}T)^{k+1} \hspace{1.0mm}L^{-1}\hspace{0.5mm} d \hspace{0.5mm}\overline{v}_{\scalebox{0.7}{$\underbrace{111 \cdots 11}_{\hspace{-0.3cm} \scalebox{1.8}{$k$}}$}}(x-l_2)/d\hspace{0.2mm}x \\
  \displaystyle \hspace{2.4cm} +(-1)^{k+1}\hspace{0.5mm}(d^2 \hspace{0.2mm}l_2/d \hspace{0.2mm}T^2)\hspace{0.5mm}(d \hspace{0.2mm}l_2/d \hspace{0.2mm}T)^{k-1} \vspace{0.2cm}\\
  \displaystyle \hspace{3.0cm} \times L^{-1}\hspace{0.5mm} \{\hspace{0.5mm} d \hspace{0.5mm}\overline{v}_{2\scalebox{0.7}{$\underbrace{11 \cdots 1}_{\scalebox{1.8}{$k-2$}}$}0}(x-l_2)/d\hspace{0.2mm}x +\hspace{0.5mm}k \hspace{1.0mm}\overline{v}_{\scalebox{0.7}{$\underbrace{111 \cdots 11}_{\hspace{-0.3cm} \scalebox{1.8}{$k$}}$}}(x-l_2) \hspace{0.5mm}\} \vspace{0.2cm}\\
  \displaystyle \hspace{8.5cm} - [\hspace{0.5mm} \mbox{terms of } \hspace{1.0mm}l_1 \hspace{0.5mm}]\hspace{0.5mm},  
\end{array} \label{eq_proof3}
\end{eqnarray}
where $L^{-1}$ is the inverse operator of $L= D (\hspace{0.5mm}d^2 \hspace{0.2mm}/d \hspace{0.2mm}x^2\hspace{0.5mm})-1$. We define the functions in Eq.~(\ref{eq_proof3}) as
\begin{eqnarray}
 \begin{array}{ll}
   \displaystyle \hspace{0.5cm}v_{k+1}(x,T) \vspace{0.2cm}\\
  \displaystyle \hspace{0.5cm} = \Pi_{k+1}(x-l_2,T) +(-1)^{k+2}\hspace{0.5mm} (d \hspace{0.2mm}l_2/d \hspace{0.2mm}T)^{k+1} \hspace{1.0mm} \overline{v}_{\scalebox{0.7}{$\underbrace{111 \cdots 11}_{\hspace{-0.3cm} \scalebox{1.8}{$k+1$}}$}}(x-l_2) \vspace{0.2cm}\\
  \displaystyle \hspace{2.5cm} +(-1)^{k+1}\hspace{0.5mm}(d^2 \hspace{0.2mm}l_2/d \hspace{0.2mm}T^2)\hspace{0.5mm}(d \hspace{0.2mm}l_2/d \hspace{0.2mm}T)^{k-1} \hspace{1.0mm} \overline{v}_{2\scalebox{0.7}{$\underbrace{11 \cdots 1}_{\hspace{-0.3cm} \scalebox{1.8}{$k-1$}}$}0}(x-l_2) \vspace{0.2cm}\\
  \displaystyle \hspace{8.5cm} - [\hspace{0.5mm} \mbox{terms of } \hspace{1.0mm}l_1 \hspace{0.5mm}]\hspace{0.5mm},  
\end{array} \label{eq_proof2}
\end{eqnarray}
where
\begin{eqnarray}
 \begin{array}{ll}
  \displaystyle \Pi_{k+1}(x,T) := \hspace{1.0mm}L^{-1}\hspace{0.5mm} \tilde{\Pi}_{k}(x,T)\hspace{0.5mm}, \vspace{0.2cm}\\
  \displaystyle \overline{v}_{\scalebox{0.7}{$\underbrace{111 \cdots 11}_{\hspace{-0.3cm} \scalebox{1.8}{$k+1$}}$}}(x):=\hspace{1.0mm}L^{-1}\hspace{0.5mm} d \hspace{0.5mm}\overline{v}_{\scalebox{0.7}{$\underbrace{111 \cdots 11}_{\hspace{-0.3cm} \scalebox{1.8}{$k$}}$}}(x)/d\hspace{0.2mm}x\hspace{0.5mm}, \vspace{0.2cm}\\
  \displaystyle \overline{v}_{2\scalebox{0.7}{$\underbrace{11 \cdots 1}_{\scalebox{1.8}{$k-1$}}$}0}(x):=L^{-1}\hspace{0.5mm} \{\hspace{0.5mm} d \hspace{0.5mm}\overline{v}_{2\scalebox{0.7}{$\underbrace{11 \cdots 1}_{\scalebox{1.8}{$k-2$}}$}0}(x)/d\hspace{0.2mm}x +\hspace{0.5mm}k \hspace{1.0mm}\overline{v}_{\scalebox{0.7}{$\underbrace{111 \cdots 11}_{\hspace{-0.3cm} \scalebox{1.8}{$k$}}$}}(x) \hspace{0.5mm}\}\hspace{0.5mm}
\end{array} \label{eq_proof4}
\end{eqnarray}
or, equivalently,
\begin{eqnarray}
 \begin{array}{ll}
  \displaystyle  L \hspace{0.5mm}\Pi_{k+1}(x,T) = \tilde{\Pi}_{k}(x,T)\hspace{0.5mm}, \vspace{0.2cm}\\
  \displaystyle L \hspace{0.5mm}\overline{v}_{\scalebox{0.7}{$\underbrace{111 \cdots 11}_{\hspace{-0.3cm} \scalebox{1.8}{$k+1$}}$}}(x) \hspace{1.0mm}=\hspace{1.0mm} d \hspace{0.5mm}\overline{v}_{\scalebox{0.7}{$\underbrace{111 \cdots 11}_{\hspace{-0.3cm} \scalebox{1.8}{$k$}}$}}(x)/d\hspace{0.2mm}x\hspace{0.5mm}, \vspace{0.2cm}\\
  \displaystyle L \hspace{0.5mm}\overline{v}_{2\scalebox{0.7}{$\underbrace{11 \cdots 1}_{\scalebox{1.8}{$k-1$}}$}0}(x) \hspace{1.0mm}= \hspace{1.0mm}d \hspace{0.5mm}\overline{v}_{2\scalebox{0.7}{$\underbrace{11 \cdots 1}_{\scalebox{1.8}{$k-2$}}$}0}(x)/d\hspace{0.2mm}x +\hspace{0.5mm}k \hspace{1.0mm}\overline{v}_{\scalebox{0.7}{$\underbrace{111 \cdots 11}_{\hspace{-0.3cm} \scalebox{1.8}{$k$}}$}}(x) \hspace{0.5mm}.
\end{array} \label{eq_proof5}
\end{eqnarray}
Hence, the statement also holds for $n=k+1$. } \hspace{0.0cm} \scalebox{1.0}{$\Box$} \vspace{0.2cm}

\clearpage

%
\section*{Acknowledgement}
This work was supported in part by JSPS KAKENHI Nos. 26247015, 15KT0100, and 17K05355, and by the JSPS A3 Foresight program (YN). The authors are grateful to Shin-Ichiro Ei (Hokkaido University, Japan), who introduced the thesis \cite{Kusaka_2007} of Tomomi Kusaka he supervised and had fruitful discussions with the authors.

%
%
%

\bibliography{mybibfile}

\end{document}